%% file: rsga_final.tex
\documentclass[lecture]{pcms-l}
\input pcmslmod-modified.tex 
\input pcmslmod14fr.tex

\usepackage{graphics,amssymb,amscd,amsthm,verbatim,ifthen}
\input epsf                 

\newcounter{colorversion}
\usepackage{color}
\definecolor{gray}{rgb}{0.7,0.7,0.7}  
\definecolor{lightgray}{rgb}{0.8,0.8,0.8}  
\newcommand{\light}[1]{\color{gray}{#1}\color{black}}
\newcommand{\verylight}[1]{\color{lightgray}{#1}\color{black}}
\newcommand{\dark}[1]{\color{black}{#1}\color{black}}

\ifthenelse{\thecolorversion=1}
{
\definecolor{light}{rgb}{1,0.1,0.1}        
\definecolor{dark}{cmyk}{1,0.4,0,0}  
\definecolor{bulletcolor}{cmyk}{1,0.3,1,0}  
\renewcommand{\light}[1]{\color{light}{#1}\color{black}}
\renewcommand{\dark}[1]{\color{dark}{#1}\color{black}}

}{}

\newtheorem{proposition}{Proposition}[chapter]
\newtheorem{theorem}[proposition]{Theorem}
\newtheorem{corollary}[proposition]{Corollary}
\newtheorem{lemma}[proposition]{Lemma}
\newtheorem{conjecture}[proposition]{Conjecture}

\newtheorem{definition}[proposition]{Definition}
\newtheorem{remark}[proposition]{Remark}

\newtheorem{example}[proposition]{Example}

\newcommand{\xx}{{\bf x}}

\newcommand{\reals}{\mathbb R}
\newcommand{\complexes}{\mathbb C}

\newcommand{\Hy}{{\mathcal H}}

\newcommand{\set}[1]{{\left\lbrace #1 \right\rbrace}}

\newcommand{\br}[1]{\langle #1 \rangle}
\newcommand{\ep}{\varepsilon}

\def\cc{\mathbf{c}}

\def\ZZ{\mathbb{Z}}
\def\AA{\mathcal{A}}

\def\RR{\mathbb{R}}

\begin{document}

\setcounter{page}{1}
\part*{Root Systems and\\ Generalized Associahedra}
\pauth{Sergey Fomin and Nathan Reading}

\tableofcontents

\mainmatter
\setcounter{page}{3}

\LogoOn

\lectureseries[Roots and Associahedra]{Root Systems and\\
    Generalized Associahedra}

\auth[Fomin and Reading]{Sergey Fomin and Nathan Reading}
\address{Department of Mathematics, 
University of Michigan, 
Ann Arbor, MI 48109-1109,
USA}
\email{fomin@umich.edu, nreading@umich.edu}
\thanks{This work was partially supported by NSF grants DMS-0245385 (S.F.) 
and DMS-0202430 (N.R.).}

\setaddress


These lecture notes provide an overview of root systems, 
generalized associahedra, and the combinatorics of clusters. 
Lectures~\ref{lec1}-\ref{lec2} cover classical material:
root systems, finite reflection groups, and the Cartan-Killing
classification. 
Lectures~\ref{lec3}--\ref{lec:cluster} provide an introduction to
cluster algebras from a combinatorial perspective. 
Lecture~\ref{lec:num} is devoted to related topics in enumerative
combinatorics.


\medskip

There are essentially no proofs but an abundance of examples.
We label unproven assertions as either ``lemma'' or ``theorem''
depending on whether they are easy or difficult to prove. 
We encourage the reader to try proving the lemmas, or at least get an idea 
of why they are true.

\medskip

For additional information on root systems, reflection groups and Coxeter 
groups, the reader is referred to \cite{Bourbaki,Fu-Ha,Humphreys}.
For basic definitions related to convex polytopes and lattice theory, 
see~\cite{Ziegler} and~\cite{Gr}, respectively.
Primary sources on generalized associahedra and cluster combinatorics 
are~\cite{gaPoly,ga,ca2}. 
Introductory surveys on cluster algebras were given in \cite{cdm,
  camblec, zel-korea}. 

\medskip

\emph{Note added in press (February 2007):}
Since these lecture notes were written, there has been much progress in
the general area of cluster algebras and Catalan combinatorics of Coxeter
groups and root systems.  We have not attempted to update the text to
reflect these most recent advances. Instead, we refer the reader to the
online Cluster Algebras Portal, maintained by the first author. 


\newpage


\noindent\textbf{Acknowledgments}

\bigskip

We thank Christos Athanasiadis, Jim Stasheff and Andrei Zelevinsky for careful
readings of earlier versions of these notes and for a number of
editorial suggestions, which led to the improvement of the paper. 

\bigskip

\noindent
S.F.: 
I am grateful to the organizers of the 2004 Graduate Summer School at
Park City (Ezra Miller, Vic Reiner, and Bernd Sturmfels)
for the invitation to deliver these lectures, and for their support,
understanding, and technical help. 

Sections~\ref{sec:matrix-mut}-\ref{sec:exchange-rel} 
and Lecture~\ref{lec:cluster} present results of an ongoing joint
project with Andrei Zelevinsky centered around cluster algebras. 

\bigskip

\noindent
N.R.: I would like to thank Vic Reiner for teaching the course which 
sparked my interest in Coxeter groups; Anders Bj\"{o}rner and Francesco 
Brenti for making a preliminary version of their forthcoming book available 
to the students in Reiner's course; and John Stembridge, whose course 
and lecture notes have deepened my knowledge of Coxeter groups and root 
systems.

\bigskip

\noindent
Some of the figures in these notes are inspired by figures produced by Satyan 
Devadoss, Vic Reiner and Rodica Simion.
Several figures were borrowed from \cite{gaPoly,ga,ca1,ca2, tptp}.

\newpage

\lecture{Reflections and Roots}
\label{lec1}

\section{The pentagon recurrence} 
\label{abel}

Consider a sequence $f_1, f_2, f_3,\dots$
defined recursively by $f_1=x$, \ $f_2=y$, and 
\begin{equation}
\label{eq:abel}
f_{n+1}=\frac{f_n+1}{f_{n-1}}\,. 
\end{equation}
Thus, the first five entries are 
\begin{equation}
\label{eq:abel2}
x,\ y,\ \frac{y+1}{x},\ \frac{x+y+1}{xy},\ \frac{x+1}{y}.
\end{equation}
Unexpectedly, the sixth and seventh entries are $x$ and $y$, 
respectively, so the sequence is periodic with period five!
We will call \eqref{eq:abel} the \emph{pentagon recurrence}.\footnote{
The discovery of this recurrence and its $5$-periodicity are 
  sometimes attributed to R.~C.~Lyness (1942); see,
  e.g.,~\cite{csornyei-laczkovich}. 
It was probably already known to N.~H.~Abel. 
This recurrence is closely related to (and easily deduced from) 
the famous ``pentagonal identity'' 
for the dilogarithm function, 
first obtained by W.~Spence (1809)  
and rediscovered by Abel (1830) and C.~H.~Hill (1830). 
See, e.g.,~\cite{lewin}. 
}

This sequence has another important property. 
{\em A priori}, we can only expect its terms to be rational functions
of $x$ and~$y$. In fact, each $f_i$ is a 
\emph{Laurent polynomial} (actually, with nonnegative integer
coefficients). 
This is an instance of what is called the {\em Laurent phenomenon}.

It will be helpful to represent this recurrence as the
evolution of a ``moving window'' consisting of two consecutive terms
$f_i$ and~$f_{i+1}$:
\[
\begin{bmatrix}
f_1\\ f_2
\end{bmatrix}
\stackrel{\textstyle\tau_1}{\longrightarrow}
\begin{bmatrix}
f_3\\ f_2
\end{bmatrix}
\stackrel{\textstyle\tau_2}{\longrightarrow}
\begin{bmatrix}
f_3\\ f_4
\end{bmatrix}
\stackrel{\textstyle\tau_1}{\longrightarrow}
\begin{bmatrix}
f_5\\ f_4
\end{bmatrix}
\stackrel{\textstyle\tau_2}{\longrightarrow}
\begin{bmatrix}
f_5\\ f_6
\end{bmatrix}
\longrightarrow\,\cdots,
\]
where the maps $\tau_1$ and $\tau_2$ are defined by 
\begin{equation}
\label{eq:tau1tau2}
\tau_1:\begin{bmatrix}f\\ g\end{bmatrix}\longmapsto
       \begin{bmatrix}\frac{g+1}{f}\\g\end{bmatrix}
\quad\text{and}\quad
\tau_2:\begin{bmatrix}f\\ g\end{bmatrix}\longmapsto
       \begin{bmatrix}f\\\frac{f+1}{g}\end{bmatrix}. 
\end{equation} 
Both $\tau_1$ and $\tau_2$ are involutions:
$\tau_1^2=\tau_2^2=1$, where $1$ denotes the identity map.
The $5$-periodicity of the recurrence~\eqref{eq:abel} translates into
the identity $(\tau_2\tau_1)^5=1$. 
That is, the group generated by $\tau_1$ and~$\tau_2$ is a dihedral
group with 10 elements. 

Let us now consider a
similar but simpler pair of maps. 
Throw away the $+1$'s that occur in the definitions of $\tau_1$ and
$\tau_2$, and take logarithms.
We then obtain a pair of linear maps 
\[
s_1:\begin{bmatrix}x\\ y\end{bmatrix}\longmapsto
       \begin{bmatrix}y-x\\y\end{bmatrix}
\quad\text{and}\quad
s_2:\begin{bmatrix}x\\ y\end{bmatrix}\longmapsto
       \begin{bmatrix}x\\x-y\end{bmatrix}. 
\]

A (linear) {\em hyperplane} in a vector space $V$ is a linear subspace
of codimension~$1$. 
A (linear) {\em reflection} is a map that fixes
all the points in some linear hyperplane, and has an eigenvalue of $-1$. 
The maps $s_1$ and $s_2$ are linear reflections 
satisfying $(s_2s_1)^3=1$.
Thus, the group $\br{s_1,s_2}$ is a dihedral group with 6 elements.

We are led to wonder if the dihedral behavior of $\br{\tau_1,\tau_2}$ is 
related to, or even explained by the dihedral behavior of $\br{s_1,s_2}$.
To test this unlikely-sounding hypothesis, let us try to find similar 
examples.
What other pairs $(s,s')$ of linear reflections generate 
finite dihedral groups? 
To keep things simple, we set $s=s_1$ and confine the choice
of~$s'$ to maps of the form 
\[
s':\begin{bmatrix}x\\ y\end{bmatrix}\longmapsto
       \begin{bmatrix}x\\L(x,y)\end{bmatrix},
\]
where $L$ is a linear function.
Keeping in mind that $s_1$ and $s_2$ arose as logarithms, we 
require that $L$ have integer coefficients.

After some work, one determines that besides $x-y$, 
the functions $2x-y$ and $3x-y$ are
the only good choices for~$L$. 
More specifically, define 
\[
s_3:\begin{bmatrix}x\\ y\end{bmatrix}\longmapsto
       \begin{bmatrix}x\\2x-y\end{bmatrix}
\quad\text{and}\quad
s_4:\begin{bmatrix}x\\ y\end{bmatrix}\longmapsto
       \begin{bmatrix}x\\3x-y\end{bmatrix}. 
\]
Then $(s_3s_1)^4=1$ and $(s_4s_1)^6=1$.
Thus, $\br{s_1,s_3}$  and $\br{s_1,s_4}$ are dihedral groups with 8
and 12 elements, respectively.

By analogy with \eqref{eq:tau1tau2}, we next define 
\[
\tau_3:\begin{bmatrix}f\\ g\end{bmatrix}\longmapsto
       \begin{bmatrix}f \\ \frac{f^2+1}{g}\end{bmatrix}
\quad\text{and}\quad
\tau_4:\begin{bmatrix}f\\ g\end{bmatrix}\longmapsto
       \begin{bmatrix}f\\\frac{f^3+1}{g}\end{bmatrix}. 
\]
Calculations show that $(\tau_3\tau_1)^6=1$, and the group 
$\br{\tau_1,\tau_3}$ is dihedral with 12 elements.  
We can think of $\tau_1$ and $\tau_3$ as defining a
``moving window'' for the sequence 
\begin{equation}
\label{eq:abel6}
x,\ y,\ \frac{y+1}{x},\ \frac{x^2+(y+1)^2}{x^2y},\
\frac{x^2+y+1}{xy},\ \frac{x^2+1}{y},\ x ,\ y,\ \ldots
\end{equation}
Notice that the Laurent phenomenon holds: 
these rational functions are Laurent polynomials---again,
with nonnegative integer coefficients. 

Likewise, $(\tau_4\tau_1)^8=1$, the group $\br{\tau_1,\tau_4}$ is dihedral with
16 elements, and $\tau_1$ and $\tau_4$ define an $8$-periodic sequence
of Laurent polynomials. 



In the first two lectures, we will develop the basic theory of finite
reflection groups that will include their complete classification. 
This theory will later help explain the periodicity and Laurentness of
the sequences discussed above, and provide appropriate algebraic
and combinatorial tools for the study of other similar recurrences. 


\section{Reflection groups}
\label{ref gp}

Our first goal will be to understand the finite groups generated by linear 
reflections in a vector space~$V$.
It turns out that for such a group, it is always possible to define a
Euclidean structure on $V$ so that all of the reflections in the group
are ordinary {\em orthogonal reflections}. 
The study of groups generated by orthogonal reflections is a classical subject, 
which goes back to the classification of Platonic solids by the 
ancient Greeks.

Let $V$ be a Euclidean space. 
In what follows, all reflecting hyperplanes pass through the origin,
and all reflections are orthogonal. 
A finite {\em reflection group} is a finite group generated by
some reflections in~$V$.
In other words, we choose a collection of hyperplanes such that the group of
orthogonal transformations generated by the corresponding reflections is finite.
Infinite reflection groups are also interesting, but in these
lectures, ``reflection group'' will always mean a finite one. 

\begin{figure}[ht]
\centerline{\begin{picture}(0,0)(-150,-95)
                \put(100,5){$ststs=tstst$}
                \put(85,50){$tst$}
                \put(37,90){$t$}
                \put(-23,90){$s$}
                \put(-75,60){$sts$}
        \put(1,60){$1$}
        \put(37,48){$t$}
        \put(-34,48){$s$}
        \put(55,16){$ts$}
        \put(-57,16){$st$}
        \put(53,-21){$tst$}
        \put(-59,-21){$sts$}
        \put(35,-53){$tsts$}
        \put(-43,-53){$stst$}
        \put(-7,-60){$ststs$}
        \put(1,-68){$\shortparallel$}
        \put(-7,-76){$tstst$}
        \end{picture}
        \epsfbox{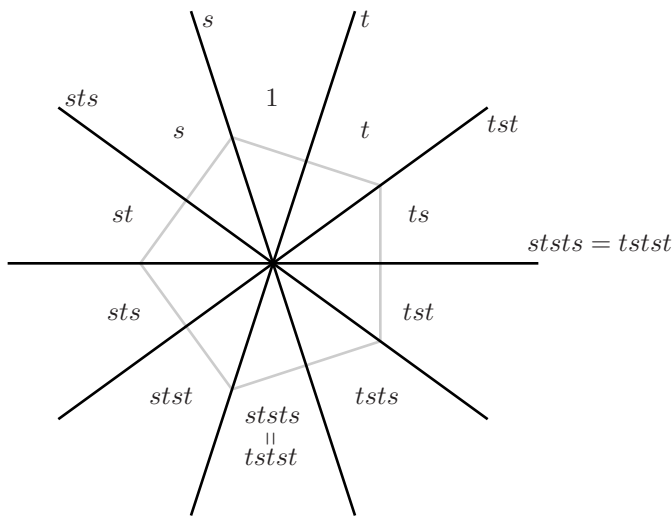}
}   
\caption{The reflection group $I_2(5)$.}
\label{I25}
\end{figure}

The set of reflections in a reflection group $W$ is typically larger
than a minimal set of reflections generating~$W$.
This is illustrated in Figure~\ref{I25}, 
where $W$ is the group of symmetries of a regular pentagon. 
This $10$-element group is generated by
two reflections $s$ and~$t$ whose reflecting lines make an angle
of~$\pi/5$. 
It consists of $5$ reflections, $4$ rotations, and the
identity element. 
In Figure~\ref{I25}, each of the $5$~lines is labeled by the corresponding
reflection.

\begin{lemma}
\label{conjug1}
If $t$ is the reflection fixing a hyperplane $H$ and $w$ an 
orthogonal transformation, then 
$wtw^{-1}$ is the reflection fixing the hyperplane $wH$.
\end{lemma}

\begin{lemma}
\label{conjug2}
Let $W$ be a finite group generated by a finite set $T$ of
 reflections.
Then the set of all 
reflections in $W$ is $\set{wtw^{-1}:w\in W,t\in T}$.
\end{lemma}

The set $\Hy$ of all reflecting hyperplanes of a reflection group $W$ is called
a {\em Coxeter arrangement}.
In light of Lemmas~\ref{conjug1} and~\ref{conjug2}, one can give an alternate 
definition of a Coxeter arrangement:
A Coxeter arrangement is a collection $\Hy$ of hyperplanes which is closed 
under reflections in the hyperplanes.
Like any hyperplane arrangement in $V$, a Coxeter arrangement cuts $V$ into 
connected components called {\em regions}.
That is, the regions are the connected components of the complement to
the union of all hyperplanes in~$\Hy$. 

The regions are in one-to-one correspondence with the elements of~$W$,
as follows. 
Once and for all, fix an arbitrary region $R_1$ to represent the
identity element. 

\begin{lemma}
The map $w\mapsto R_w\stackrel{\rm def}{=}w(R_1)$
is a bijection between a reflection group~$W$ and the set of regions
of the corresponding Coxeter arrangement. 
\end{lemma}

To illustrate, each of the $10$ regions in Figure~\ref{I25}
is labeled by the corresponding element of the group. 

The choice of a region representing the identity element leads to a 
distinguished choice of a minimal set of generating reflections.
The {\em facet hyperplanes} of $R_1$ are the hyperplanes in $\Hy$ whose 
intersection with the closure of $R_1$ has dimension $n-1$.

\begin{lemma}
\label{gen}
The reflections in the facet hyperplanes of $R_1$ generate~$W$.
This generating set is minimal by inclusion.
\end{lemma}

\section{Symmetries of regular polytopes} 
\label{sec:regular-polytopes}

A \emph{regular polytope} in a Euclidean space is a convex polytope whose
symmetry group (i.e., the group of isometries of the space that leave the
polytope invariant) acts transitively on \emph{complete flags of
  faces}, i.e., on nested collections of the form
\[
\text{vertex}\subset\text{edge}\subset\text{$2$-dim.\ face}\subset\cdots
\]

\begin{theorem}
\label{regular}
The symmetry group of any regular polytope is a reflection group.
\end{theorem}

The converse is false---see Remark~\ref{rem:dn}. 

We illustrate Theorem~\ref{regular} with several concrete examples.

\begin{example}\rm
\label{mgon}
Consider a regular $m$-gon on a Euclidean plane, centered at the origin.
The symmetry group of the $m$-gon is denoted by~$I_2(m)$.
This group contains (and is generated by)
$m$ reflections, which correspond to 
the $m$ lines of reflective symmetry of the $m$-gon. 

The group $I_2(m)$ is a dihedral group with $2m$ elements. 
It is generated by two reflections $s$ and~$t$ satisfying 
$(st)^m=1$. 
To define $s$ and~$t$, we use the construction of Lemma~\ref{gen}. 
Pick a side of the polygon, and consider two
reflecting lines: one perpendicular to the side and another passing
through one of its endpoints. 
The case $m=5$ is shown in Figure~\ref{I25}. 
\end{example}

\begin{figure}[ht]
\centerline{\begin{picture}(0,0)(-100,-100)
                \put(7,-65){$(0,0,1)$}
                \put(57,20){$(1,0,0)$}
                \put(-83,20){$(0,1,0)$}
                \put(7,65){$x = y$}
                \put(70,-33){$x = z$}
                \put(-58,-44){$y = z$}
        \end{picture}
        \epsfbox{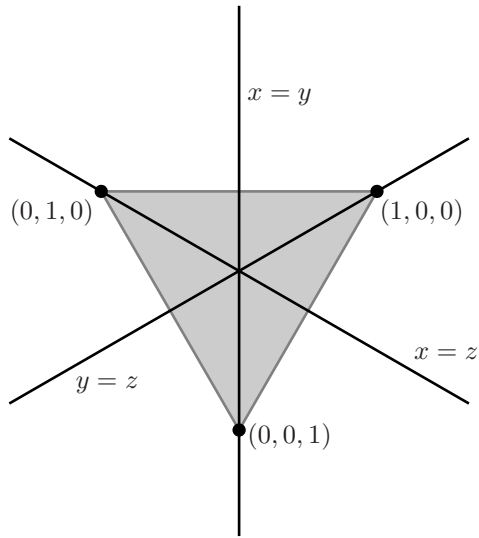}
}   
\caption{The reflection group $A_2$.
}
\label{A2simp}
\end{figure}

\begin{example}\rm
\label{An}
Take a regular tetrahedron in $3$-space, with the vertices labeled 1,
 2, 3, and~4.
Its symmetry group is obviously isomorphic to 
the {\em symmetric group}~$\mathcal{S}_4$, which consists of the permutations of
the set $\set{1,2,3,4}$.
For each edge of the tetra\-hedron, choose a plane which is
perpendicular to the edge and contains the other two vertices.
Reflections in these six hyperplanes generate the symmetry group.

In general, the symmetry group of a regular simplex can be described
as follows. 
Let $(e_1,\dots,e_{n+1})$ be the standard basis in~$\reals^{n+1}$.
The {\em standard $n$-dimensional simplex} (or $n$-simplex) is the convex hull of
the endpoints of the vectors $e_1,\ldots,e_{n+1}$.
Thus the standard $1$-simplex is a line segment in~$\reals^2$, the standard 
$2$-simplex is an equilateral triangle in~$\reals^3$, and the 
standard $3$-simplex is the regular tetrahedron described above,
sitting in~$\reals^4$. 
The symmetry group $A_n$ of the standard $n$-simplex is
canonically isomorphic to~$\mathcal{S}_{n+1}$, the symmetric group of permutations of the set 
$[n+1]\stackrel{\rm def}{=}\set{1,2,\ldots,n+1}$. 

For each edge $[e_i,e_j]$ of the standard simplex, there is a
hyperplane $x_i-x_j=0$ 
perpendicular to the edge and containing all the other vertices.
Reflection through this hyperplane interchanges the endpoints of the edge and 
fixes the rest of the vertices.
These $\binom{n+1}{2}$ reflections generate~$A_n$. 

To construct a minimal generating set of reflections,
we again use Lemma~\ref{gen}.  
Let $R_1$ be the connected component of the complement to the
$\binom{n+1}{2}$ reflecting hyperplanes defined by 
\begin{equation}
\label{eq:x1<x2<...}
R_1=\{x_1<x_2<\cdots < x_{n+1}\}.
\end{equation}
The facet hyperplanes of $R_1$ are given by the equations
\[
x_i-x_{i+1}=0, \quad \text{for $i=1,\dots,n$}.
\]
Then Lemma~\ref{gen} reduces to the well-known fact that 
the symmetric group $\mathcal{S}_{n+1}$ is generated by the adjacent transpositions
$s_1,\dots,s_n$. 
(Here each $s_i$ exchanges $i$ and~$i+1$, keeping everything else in
its place.) 

Figure~\ref{A2simp} illustrates the special case $n\!=\!2$, 
the symmetry group of the standard $2$-simplex (shaded).
The plane of the page represents the plane \hbox{$x+y+z=1$} in~$\reals^3$.
\end{example}

\begin{example}\rm
\label{Bn}
The {\em $n$-crosspolytope} is the convex hull
of (the endpoints of) the vectors $\pm e_1,\pm e_2,\ldots,\pm e_n$ in~$\reals^n$. 
For example, the $3$-crosspolytope is the regular octahedron.
The symmetry group of this polytope is the {\em hyperoctahedral
group}~$B_n$.
As in the previous examples, it is generated by the reflections it
contains. 

The special case $n=3$ 
(the symmetry group $B_3$ of a regular octahedron) 
is shown in Figure~\ref{B3pic}.
The dotted lines show the intersections of reflecting hyperplanes with
the front surface of the octahedron.
Each edge of the octahedron is also contained in a reflecting plane.

\begin{figure}[ht]
\centerline{\epsfbox{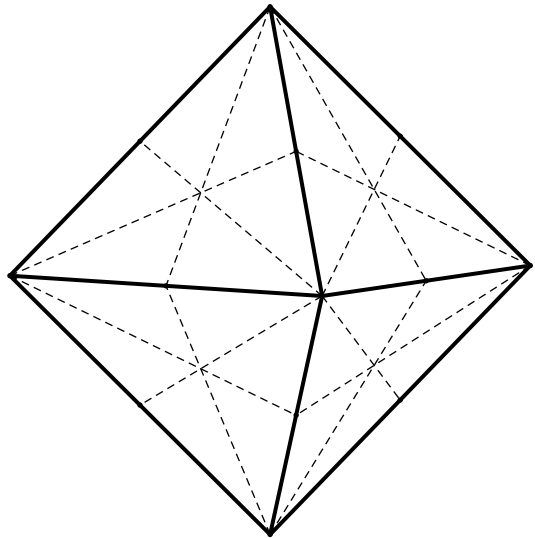}}   
\caption{The reflection group $B_3$}
\label{B3pic}
\end{figure}

There are two types of reflections in the symmetry group of the crosspolytope.
One type of reflection transposes a vertex with its negative and fixes all 
other vertices.
Also, for each pair $i\neq j$, there is a reflection 
which transposes $e_i$ and $e_j$, transposes
$-e_i$ and $-e_j$, and fixes all other vertices.

To construct a minimal set of reflections generating~$B_n$, take the 
minimal generating set for $A_{n-1}$ given in Example~\ref{An} 
and adjoin the reflection that interchanges $e_1$ and~$-e_1$.

The group $B_n$ is also the symmetry group of the $n$-dimensional cube.
\end{example}

\begin{figure}[ht]
\centerline{\epsfbox{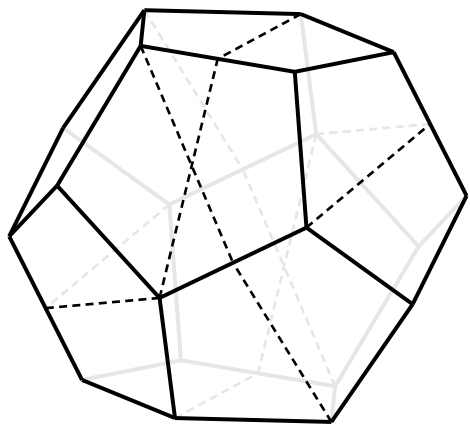}}   
\caption{The reflection group $H_3$}
\label{H3pic}
\end{figure}

\begin{example}\rm
\label{H3}
The symmetry group of a regular dodecahedron (or a regular icosahedron) is
the reflection group~$H_3$.
Figure~\ref{H3pic} shows the dodecahedron and a minimal set of three 
reflections generating its symmetry group.
The dotted lines show the intersections of the corresponding three
hyperplanes with the front surface of the dodecahedron. 
\end{example}

\begin{example}\rm
\label{H4}
In $4$-space, there are six types of regular polytopes.  
The obvious three are the $4$-simplex, the $4$-cube, and the $4$-crosspolytope.
There are two regular polytopes whose symmetry group is the reflection group 
called~$H_4$.
One of these, the {\em $120$-cell}, has 600 vertices and 120
dodecahedral faces; 
the other, the {\em $600$-cell}, has 120 vertices and 600 tetrahedral faces.
The remaining regular $4$-dimensional polytope is the {\em $24$-cell},
with 24 vertices and 24 octahedral faces.
Its symmetry group is a reflection group denoted by~$F_4$.
\end{example}

Not every reflection group is the 
symmetry group of a regular polytope.
A counterexample is constructed as follows. 

\begin{example}\rm
\label{Dn}
Let $n\geq 3$. 
Returning to the crosspolytope, ignore the reflections which transpose an 
opposite pair of vertices.
The remaining reflections generate a reflection group called~$D_n$, which is a 
proper subgroup of~$B_n$.
The reflections of $D_3$ are represented by the dotted lines in 
Figure~\ref{B3pic}.  
We note that the Coxeter arrangements of types $A_3$ and~$D_3$ are
related by an orthogonal transformation, 
so the reflection groups $A_3$ and $D_3$ are isomorphic to each other.
\end{example}


\begin{remark}
\label{rem:dn}
{\rm
It can be shown that, for $n\geq 4$, the group $D_n$ is not a 
symmetry group of a regular polytope.
See Section~\ref{sec:other} for further details. 
}
\end{remark}

\section{Root systems}

Root systems are configurations of vectors obtained by
replacing each reflecting hyperplane of a reflection group by a pair
of opposite normal vectors; 
the resulting configuration should be invariant under the action of
the group. Here is a formal definition. 
A finite {\em root system} is a finite non-empty collection
$\Phi$ of nonzero vectors in $V$ called {\em roots} with the following
properties: 
\begin{enumerate}
\item[(i)] Each one-dimensional subspace of $V$ either contains no roots, or 
contains two roots $\pm\alpha$.
\item[(ii)] For each $\alpha\in\Phi$, the reflection $\sigma_\alpha$ permutes 
$\Phi$.  
\end{enumerate}
The following lemma shows that the study of root systems is
essentially equivalent to the study of reflection groups. 

\begin{lemma}
\label{root ref}
For a finite root system~$\Phi$, 
the group generated by the reflections
$\{\sigma_\alpha:\alpha\!\in\!\Phi\}$ is finite. 
The corresponding reflecting hyperplanes 
form a Coxeter arrangement.
Conversely, for any reflection group~$W$, there is a root system~$\Phi$ such
that the orthogonal reflections $\{\sigma_\alpha\}_{\alpha\in\Phi}$  
are precisely the reflections in~$W$. 
\end{lemma}

In Section~\ref{ref gp}, we fixed a region $R_1$ of the associated Coxeter 
arrangement $\Hy$.
The {\em simple roots} in~$\Phi$ 
are the roots normal to the facet hyperplanes of~$R_1$ and
pointing into the half-space containing~$R_1$.
The {\em rank} of $\Phi$ is the cardinality $n$ of the set of simple
roots~$\Pi$. 
Since $W$ acts transitively on the regions of $\Hy$, the rank of $\Phi$ does 
not depend on the choice of~$\Pi$, and is equal to the dimension of
the linear span of~$\Phi$. 
It will be convenient to fix an indexing set $I$ so that
$\Pi=\set{\alpha_i:i\in I}$.
The standard choice is $I=[n]=\{1,\dots,n\}$.

For any $\alpha\in\Phi$, the coefficients $c_i$ in the expansion $\alpha=\sum_{i\in
  I}c_i\alpha_i$ are called the {\em simple root coordinates} of~$\alpha$. 
The set $\Phi_+$ of {\em positive roots} consists of all roots whose
simple root coordinates are all non-negative. 
The \emph{negative roots} $\Phi_-$ are those with non-positive
simple root coordinates.

\begin{lemma}
\label{pm}
$\Phi$ is the disjoint union of $\Phi_+$ and $\Phi_-$.
\end{lemma}

\begin{figure}[ht]
\centerline{
        \begin{picture}(0,0)(-60,-30)
                \put(30,-27){$\beta$}
                \put(0,64){$\alpha$}
                \put(21,28){$\sigma_\alpha(\beta)$}
                \put(66,5){$a_{\alpha\beta}\alpha$}
        \end{picture}
        \epsfbox{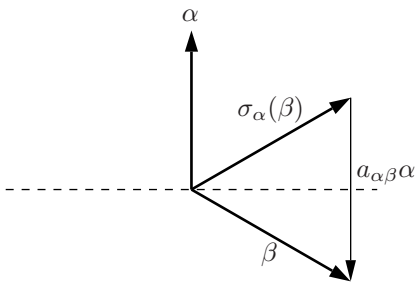}
}
\caption{Reflecting $\beta$ in the hyperplane perpendicular to~$\alpha$.}
\label{crys}
\end{figure}

In these lectures, we focus on the study of the important class
of finite \emph{crystallographic} root systems. 
These are the finite non-empty collections of vectors that, 
in addition to the axioms~(i)--(ii) above, satisfy the ``crystallographic
condition'' 
\begin{enumerate}
\item[(iii)] For any $\alpha,\beta\!\in\!\Phi$, 
we have $\sigma_\alpha(\beta)=\beta-a_{\alpha\beta}\alpha$ with 
$a_{\alpha\beta}\!\in\!\mathbb{Z}$. (See Figure~\ref{crys}.) 
\end{enumerate}
Equivalently, the simple root coordinates of any root are integers. 

For the rest of these lectures, a ``root system'' will always
be presumed finite and crystallographic. 

\begin{figure}[ht]
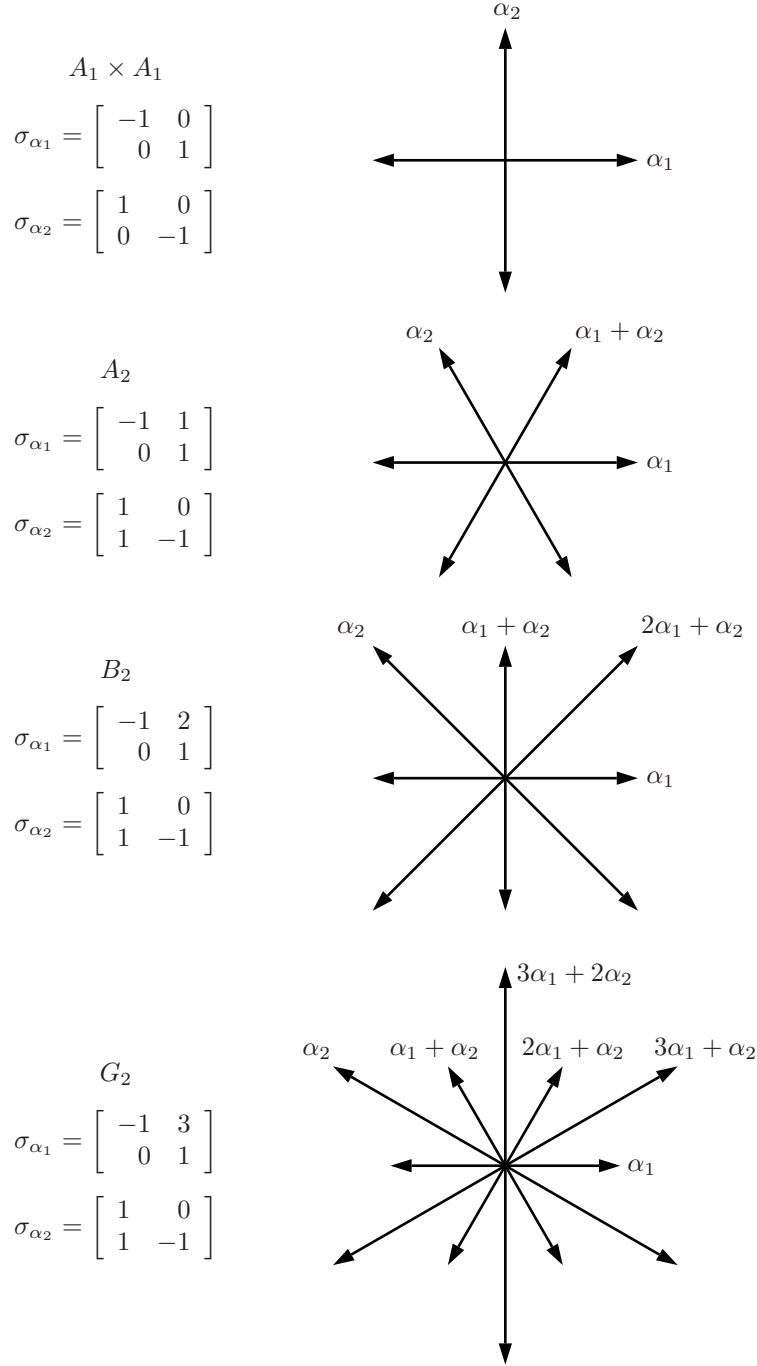

\begin{tabular}{cc}
\hbox{$\begin{array}{c}
A_1\times A_1\\[.1in]
\sigma_{\alpha_1}=\left[\begin{array}{rr}-1&0\\0&1\end{array}\right]\\[.2in]
\sigma_{\alpha_2}=\left[\begin{array}{rr}1&0\\0&-1\end{array}\right]
\end{array}$}&
\hbox{\begin{picture}(0,0)(-60,0)
                \put(57,-2){$\alpha_1$}
                \put(-1,55){$\alpha_2$}
        \end{picture}
        \epsfbox{A1A1root.ps}
}\\[.7in]
$\begin{array}{c}A_2\\[.1in]
\sigma_{\alpha_1}=\left[\begin{array}{rr}-1&1\\0&1\end{array}\right]\\[.2in]
\sigma_{\alpha_2}=\left[\begin{array}{rr}1&0\\1&-1\end{array}\right]
\end{array}$&
\hbox{\begin{picture}(0,0)(-60,0)
                \put(57,-2){$\alpha_1$}
                \put(-34,47){$\alpha_2$}
                \put(30,47){$\alpha_1+\alpha_2$}
        \end{picture}
\epsfbox{A2root.ps}
}\\[.7in]
$\begin{array}{c}B_2\\[.1in]
\sigma_{\alpha_1}=\left[\begin{array}{rr}-1&2\\0&1\end{array}\right]\\[.2in]
\sigma_{\alpha_2}=\left[\begin{array}{rr}1&0\\1&-1\end{array}\right]\\\\
\end{array}$&
\hbox{\begin{picture}(0,0)(-60,0)
                \put(57,-2){$\alpha_1$}
                \put(-60,55){$\alpha_2$}
                \put(55,55){$2\alpha_1+\alpha_2$}
                \put(-13,55){$\alpha_1+\alpha_2$}
        \end{picture}
\epsfbox{B2root.ps}
}\\[.8in]
$\begin{array}{c}G_2\\[.1in]
\sigma_{\alpha_1}=\left[\begin{array}{rr}-1&3\\0&1\end{array}\right]\\[.2in]
\sigma_{\alpha_2}=\left[\begin{array}{rr}1&0\\1&-1\end{array}\right]
\end{array}$&
\hbox{\begin{picture}(0,0)(-90,0)
                \put(50,-2){$\alpha_1$}
                \put(60,42){$3\alpha_1+\alpha_2$}
                \put(-40,42){$\alpha_1+\alpha_2$}
                \put(10,42){$2\alpha_1+\alpha_2$}
                \put(8,70){$3\alpha_1+2\alpha_2$}
                \put(-73,42){$\alpha_2$}
        \end{picture}
\epsfbox{G2root.ps}
}\\\\\\\\
\end{tabular}
\caption{The finite crystallographic root systems of rank~$2$}
\label{rank2}
\end{figure}

\begin{example}\rm
\label{ranks12}
A root system of rank~$1$ is called~$A_1$; it consists of a pair of
vectors~$\pm\alpha$. 
There are four non-isomorphic (finite crystallographic)
root systems of rank~$2$, 
called $A_1\times A_1$, $A_2$, $B_2$ and~$G_2$; see
Figure~\ref{rank2}.
\end{example}
\vspace{-.1in}

For the root systems $A_2$, $B_2$ and $G_2$, the 
reflections $\sigma_{\alpha_1}$ and $\sigma_{\alpha_2}$ 
have appeared earlier  in Section~\ref{abel}. 
(The matrices of these reflections in the basis
$(\alpha_1,\alpha_2)$  of simple roots are shown in
Figure~\ref{rank2}.) 
In these three cases, 
the pair $(\sigma_{\alpha_1},\sigma_{\alpha_2})$ 
coincides with $(s_2,s_1)$, $(s_3,s_1)$, and $(s_4,s_1)$, 
respectively, in the notation of Section~\ref{abel}.

\section{Root systems of types A, B, C, and D}

Here we present four classical families
of root systems, traditionally denoted by~$A_n$,
$B_n$, $C_n$ and~$D_n$. 
The corresponding reflection groups have types $A_n$,
$B_n$, $B_n$ and~$D_n$  
(cf.\ Examples~\ref{An}, \ref{Bn}, and~\ref{Dn}).  
In each case, $n$ is the rank of a root system. 

We realize each root system inside a Euclidean space with a fixed 
orthonormal 
basis $(e_1,e_2,\dots)$,
and describe particular choices of the sets of simple and positive
roots. 
There is no ``canonical'' way to make these choices.
Our realizations of root systems coincide with those  in~\cite{Bourbaki, Humphreys},
but our choices of simple/positive roots 
 (which are motivated by notational convenience alone)
are different. 

\subsection*{The root system $A_n$}

The root system $A_n$ can be realized as the set of vectors $e_i-e_j$ in 
$\reals^{n+1}$ with $i\neq j$.
Let $R_1$ be given by~\eqref{eq:x1<x2<...}.
Then the $n$ simple roots are
$\alpha_i\stackrel{\rm def}{=}e_{i+1}-e_i$, for $i=1,\dots,n$,
and the positive roots are $e_i-e_j$, for $1\le j<i\le n+1$.

\begin{figure}[ht]
\centerline{
        \begin{picture}(0,0)(-100,-87)
                \put(108,-2){$\alpha_1$}
                \put(-60,93){$\alpha_2$}
                \put(-1,-68){$\alpha_3$}
                \put(58,30){$\alpha_1+\alpha_2+\alpha_3$}
                \put(40,93){$\alpha_1+\alpha_2$}
                \put(-83,30){$\alpha_2+\alpha_3$}
        \end{picture}
        \epsfbox{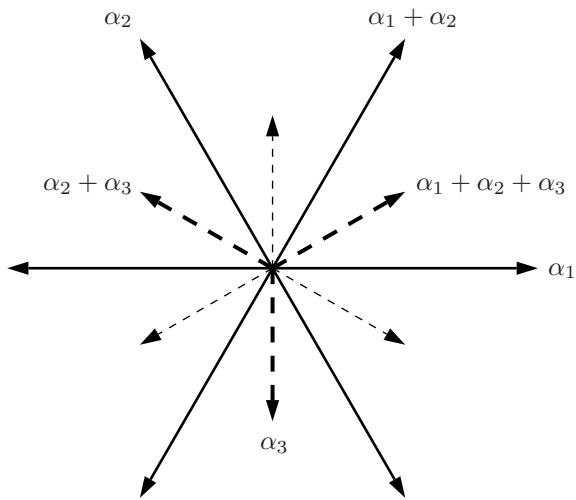}
}
\caption{The root system $A_3$
}
\label{A3root}
\end{figure}

Figure~\ref{A3root} shows a planar projection of the root system~$A_3$.
The positive roots are labeled by their simple root coordinates.
The solid lines are in the plane of the page.
Thick dotted lines are above the plane, while thin dotted lines are
below it.  

\subsection*{The root systems $B_n$ and $C_n$}

\begin{figure}[ht]
\centerline{
        \begin{picture}(0,0)(-110,-110)
                \put(29,-37){$\alpha_1$}
                \put(-83,38){$\alpha_2$}
                \put(-52,-106){$\alpha_3$}
        \end{picture}
\scalebox{0.75}{\epsfbox{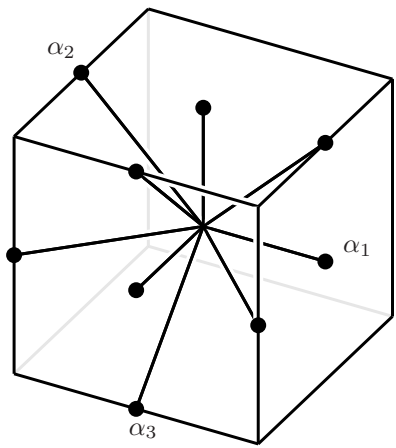}}
}
\caption{The root system $B_3$.
The endpoints of the $9$ positive roots are shown as black circles on
the cube's front.
The negative roots are not shown.
}
\label{B3root}
\end{figure}

The root system $B_n$ can be realized as the set of vectors in $\reals^n$ of 
the form $\pm e_i$ or $\pm e_i\pm e_j$ with $i\neq j$.
Choose $R_1=\{0<x_1<x_2<\cdots<x_n\}$.
Then the vectors $\alpha_0=e_1$ and 
$\alpha_i=e_{i+1}-e_i$ for $i\in[n-1]$ form a set of simple roots.
The positive roots are $e_i$ for $i\in [n]$ and $e_i\pm e_j$ for 
$1\le j<i\le n$.
See Figure~\ref{B3root}.

\begin{figure}[ht]
\centerline{
        \begin{picture}(0,0)(-120,-120)
                \put(70,-29){$\alpha_1$}
                \put(-83,22){$\alpha_2$}
                \put(-10,-80){$\alpha_3$}
        \end{picture}
\scalebox{0.75}{\epsfbox{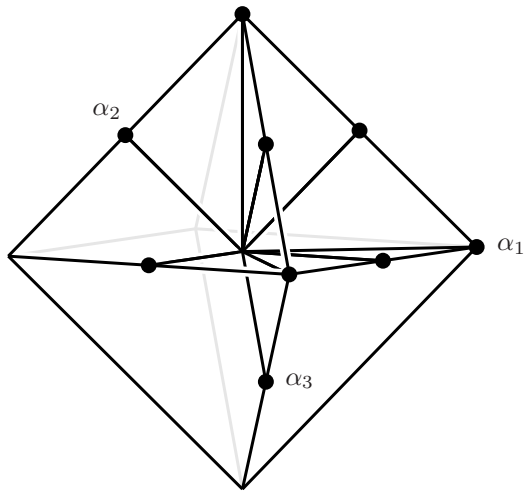}}
}
\caption{The root system $C_3$.
The endpoints of the $9$ positive roots are shown 
on the front of the octahedron.
The negative roots are not shown.
}
\label{C3root}
\end{figure}

The root system $C_n$ can be realized as the set of vectors in $\reals^n$ of 
the form $\pm 2 e_i$ or $\pm e_i\pm e_j$.
The vectors $\alpha_0=2 e_1$ and $\alpha_i=e_{i+1}-e_i$ form a set of simple 
roots. The positive roots are $2 e_i$ and $e_i\pm e_j$.
See Figure~\ref{C3root}.

The root system $C_n$ is a rescaling of $B_n$, so the corresponding reflection 
groups $W$ coincide.
In contrast to the type~$A_n$, 
the action of $W$ on the roots of $B_n$ or~$C_n$ is not transitive: 
there are two orbits, corresponding to two different lengths of roots.

\subsection*{The root system $D_n$}

The root system $D_n$ can be realized as the vectors $\pm e_i\pm e_j$ with 
$i\neq j$.
One choice of simple roots is $\alpha_0=e_2+e_1$ and
$\alpha_i=e_{i+1}-e_i$,  
giving the positive roots $e_i\pm e_j$ for $1\le j<i\le n$.
This comes from setting $R_1=\set{-x_2<x_1<x_2<\cdots<x_n}$.

\lecture{Dynkin Diagrams and Coxeter Groups}
\label{lec2}

\section{Finite type classification}

The most fundamental result in the theory of (finite crystallographic)
root systems is their complete classification, obtained by
W.~Killing and E.~Cartan in late nineteenth -- early twentieth century.
(See the historical notes in~\cite{Bourbaki}.) 
To present this classification, we will need a few preliminaries. 

First, we will need the notion of isomorphism.
The {\em ambient space} $Q_\reals=Q_\reals(\Phi)$ of a root
system~$\Phi$ is the real span of~$\Phi$. 
It inherits a Euclidean structure from~$V$.
Root systems $\Phi$ and $\Phi'$ are \emph{isomorphic} if there is an isometry 
map $Q_\reals(\Phi)\to Q_\reals(\Phi')$ of their ambient spaces
that sends $\Phi$ to some dilation $c\Phi'$ of~$\Phi'$.


The {\em Cartan matrix} of a root system $\Phi$ is the integer
matrix $[a_{ij}]_{i,j\in I}$, 
where $a_{ij}$ is such that 
$\sigma_{\alpha_i}(\alpha_j)=\alpha_j-a_{ij}\alpha_i$, as in
part~(iii) of the definition of a root system.
(This convention
  agrees with \cite{ca2,kac} but is 
``transposed'' to the one in \cite{Bourbaki,Humphreys}.)

\begin{lemma}
Root systems $\Phi$ and $\Phi'$ are isomorphic if and only if they
have the same Cartan matrix, up to simultaneous rearrangement of rows
and columns. 
\end{lemma}

\begin{example}
\label{example:cartan-rank2}
\rm
The Cartan matrices for the root systems of rank two are:
\[\begin{array}{rlcrl}
A_1\times A_1:&\left[\begin{array}{rr}2&0\\0&2\end{array}\right]&&
A_2:&\left[\begin{array}{rr}2&-1\\-1&2\end{array}\right]\\[.2in]
B_2:&\left[\begin{array}{rr}2&-2\\-1&2\end{array}\right]&&
G_2:&\left[\begin{array}{rr}2&-3\\-1&2\end{array}\right]\end{array}\]
\end{example}

\begin{example}
\label{example:a4b4c4d4}
The Cartan matrices for the root systems of type $A_4$, $B_4$, $C_4$,
and~$D_4$ are, respectively: 
\[
\begin{array}{ll}
A_4:\quad
\left[\begin{array}{rrrrr}
2&-1&0&0\\
-1&2&-1&0\\
0&-1&2&-1\\
0&0&-1&2\\
\end{array}\right] &\qquad
B_4:\quad 
\left[\begin{array}{rrrrr}
2&-2&0&0\\
-1&2&-1&0\\
0&-1&2&-1\\
0&0&-1&2\\
\end{array}\right]\\[.4in]
C_4:\quad 
\left[\begin{array}{rrrrr}
2&-1&0&0\\
-2&2&-1&0\\
0&-1&2&-1\\
0&0&-1&2\\
\end{array}\right] &\qquad
D_4:\quad 
\left[\begin{array}{rrrrr}
2&0&-1&0\\
0&2&-1&0\\
-1&-1&2&-1\\
0&0&-1&2\\
\end{array}\right]
\end{array}
\]
\end{example}

The Cartan matrices of (finite crystallographic) root systems 
are sometimes called \emph{Cartan matrices of finite type}. 
This class of matrices is completely characterized by several
elementary properties.

\begin{theorem}
\label{cartan thm}
An integer $n\times n$ matrix $[a_{ij}]$ is a Cartan matrix of a root
system if and only if 
\begin{enumerate}
\item[(i)] $a_{ii}=2$ for every $i$; 
\item[(ii)] $a_{ij}\leq 0$ for any $i\neq j$, with $a_{ij}=0$ if and only if $a_{ji}=0$;
\item[(iii)] there exists a diagonal matrix $D$ with positive diagonal
  entries such 
that $DAD^{-1}$ is symmetric and positive definite.
\end{enumerate}
\end{theorem}

\begin{remark}
{\rm
Condition~(iii) can be replaced by 
\begin{enumerate}
\item[(iii${}'$)] \emph{there exists a diagonal matrix $D'$ with positive
  integer diagonal
  entries such that $D'A$ is symmetric and positive definite.}
\end{enumerate}
}
\end{remark}

\begin{example}\rm
For the root systems $A_1\times A_1$ and $A_2$, the $2\times 2$
identity matrix serves as~$D$. 
For $B_2$ and~$G_2$, take $D=\Bigl[\begin{smallmatrix}1& 0\\0 &
  \sqrt{2}\end{smallmatrix}\Bigr]$ 
and $D=\Bigl[\begin{smallmatrix}1& 0\\0 &
  \sqrt{3}\end{smallmatrix}\Bigr]$, respectively. 
\end{example}

The characterization in Theorem~\ref{cartan thm} can be used to completely classify the
Cartan matrices of finite type, or the corresponding root systems. 
It turns out that each of those is built from blocks taken from a
certain relatively short list. Let us be more precise.

A root system $\Phi$ is called \emph{reducible}
if $\Phi$ is a disjoint union of root systems $\Phi_1$ and $\Phi_2$ such that 
every $\beta_1\in \Phi_1$ is normal to every $\beta_2\in \Phi_2$.
If such a decomposition does not exist,  $\Phi$ is called {\em irreducible}.
The parallel definition for Cartan matrices is that a Cartan matrix of
finite type is 
\emph{indecomposable} if its rows and columns cannot be simultaneously rearranged to
bring the matrix into block-diagonal form with more than one block.

The Cartan matrices of finite type can be encoded by their {\em Dynkin 
diagrams}.
The vertices of a Dynkin diagram are labeled by the elements of the
indexing set~$I$; thus they are in bijection with the simple roots. 
Each pair of vertices $i$ and $j$ is then connected as shown below
(with the vertex~$i$ on the left): 

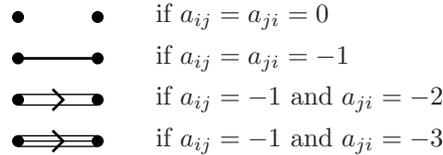
\centerline{\begin{tabular}{ccl}
\epsfbox{A1A1dynk.ps}& & if $a_{ij}=a_{ji}=0$\\
\epsfbox{A2dynk.ps}& & if $a_{ij}=a_{ji}=-1$\\
\epsfbox{B2dynk.ps}& & if $a_{ij}=-1$ and $a_{ji}=-2$\\
\epsfbox{G2dynk.ps}& & if $a_{ij}=-1$ and $a_{ji}=-3$\\[.1in]
\end{tabular}}

\noindent
(It follows from Theorem~\ref{cartan thm} that these are the only
possible pairs of values for $a_{ij}$ 
and~$a_{ji}$.
Cf.\ Example~\ref{example:cartan-rank2}.) 

\begin{lemma}
A Cartan matrix of finite type (resp., a root system)
is indecomposable (resp., irreducible) if and only if its Dynkin
diagram is connected. 
\end{lemma}

\begin{figure}[ht]
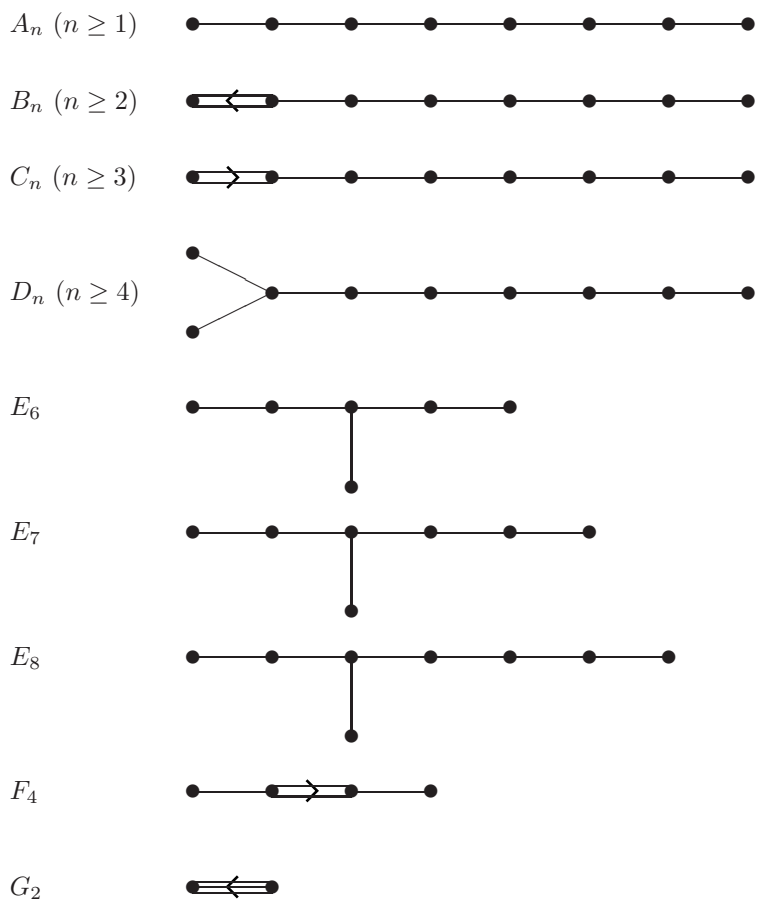

\vspace{-.2in} 
\[ 
\begin{array}{lcl} 
A_n\  (n\geq 1) && 
\setlength{\unitlength}{1.5pt} 
\begin{picture}(140,17)(0,-2) 
\put(0,-0.15){\line(1,0){140}} 
\multiput(0,0)(20,0){8}{\circle*{3}} 
\end{picture}\\
B_n\ (n\geq 2)
&& 
\setlength{\unitlength}{1.5pt} 
\begin{picture}(140,17)(0,-2) 
\put(0,0){\epsfbox{leftdir.ps}}
\put(20,-0.15){\line(1,0){120}}
\put(0,1.25){\line(1,0){20}} 
\put(0,-1.5){\line(1,0){20}} 
\multiput(0,0)(20,0){8}{\circle*{3}} 
\end{picture} \\
C_n\ (n\geq 3)
&& 
\setlength{\unitlength}{1.5pt} 
\begin{picture}(140,17)(0,-2) 
\put(0,0){\epsfbox{rightdir.ps}}
\put(20,-0.15){\line(1,0){120}}
\put(0,1.25){\line(1,0){20}} 
\put(0,-1.5){\line(1,0){20}} 
\multiput(0,0)(20,0){8}{\circle*{3}} 
\end{picture} 
\\[.2in] 
D_n\ (n\geq 4)
&& 
\setlength{\unitlength}{1.5pt} 
\begin{picture}(140,17)(0,-2) 
\put(20,-0.15){\line(1,0){120}} 
\put(0,10){\line(2,-1){20}} 
\put(0,-10){\line(2,1){20}} 
\multiput(20,0)(20,0){7}{\circle*{3}} 
\put(0,10){\circle*{3}} 
\put(0,-10){\circle*{3}} 
\end{picture} 
\\[.2in] 
E_6 
&& 
\setlength{\unitlength}{1.5pt} 
\begin{picture}(140,17)(0,-2) 
\put(0,-0.15){\line(1,0){80}} 
\put(40,0){\line(0,-1){20}} 
\put(40,-20){\circle*{3}} 
\multiput(0,0)(20,0){5}{\circle*{3}} 
\end{picture} 
\\[.25in] 
E_7 
&& 
\setlength{\unitlength}{1.5pt} 
\begin{picture}(140,17)(0,-2) 
\put(0,-0.15){\line(1,0){100}} 
\put(40,0){\line(0,-1){20}} 
\put(40,-20){\circle*{3}} 
\multiput(0,0)(20,0){6}{\circle*{3}} 
\end{picture} 
\\[.25in] 
E_8 
&& 
\setlength{\unitlength}{1.5pt} 
\begin{picture}(140,17)(0,-2) 
\put(0,-0.15){\line(1,0){120}} 
\put(40,0){\line(0,-1){20}} 
\put(40,-20){\circle*{3}} 
\multiput(0,0)(20,0){7}{\circle*{3}} 
\end{picture} 
\\[.3in] 
F_4 
&& 
\setlength{\unitlength}{1.5pt} 
\begin{picture}(140,17)(0,-2) 
\put(20,0){\epsfbox{rightdir.ps}}
\put(20,1.25){\line(1,0){20}} 
\put(20,-1.5){\line(1,0){20}} 
\put(0,-0.25){\line(1,0){20}} 
\put(40,-0.25){\line(1,0){20}} 
\multiput(0,0)(20,0){4}{\circle*{3}} 
\end{picture} 
\\[.1in] 
G_2 
&& 
\setlength{\unitlength}{1.5pt} 
\begin{picture}(140,17)(0,-2) 
\put(0,0){\epsfbox{leftdir.ps}}
\put(0,1.25){\line(1,0){20}} 
\put(0,-1.5){\line(1,0){20}} 
\put(0,-0.15){\line(1,0){20}} 
\multiput(0,0)(20,0){2}{\circle*{3}} 
\end{picture} 
\end{array} 
\] 
\vspace{-.1in} 
\caption{Dynkin diagrams of finite irreducible root systems.}
\label{diagrams} 
\end{figure}

\begin{theorem}[Cartan-Killing classification of irreducible root
  systems and Cartan matrices of finite type] 
The complete list of Dynkin diagrams of irreducible root systems is presented 
in Figure~\ref{diagrams}.
\end{theorem}

\pagebreak

Root systems are just one example among a large number of mathematical objects of
``finite type'' which are classified by (some class of) 
Dynkin diagrams.  
The appearance of the ubiquitous Dynkin diagrams in a variety of
seemingly unrelated classification problems has
fascinated several generations of mathematicians,
and helped establish nontrivial connections between different
areas of mathematics. 
See Section~\ref{sec:other} and references therein. 

\section{Coxeter groups} 

Let $\Phi$ be a (finite crystallographic)
root system and $\alpha\neq\beta$ a pair of roots in~$\Phi$. 
The angle between the 
corresponding reflecting hyperplanes is a rational multiple of $\pi$
with denominator 2, 3, 4 or~6.
Thus the rotation $\sigma_\alpha\sigma_\beta$ has order 2, 3, 4, or 6 as an 
element of the associated reflection group $W$.
The insight that the order of a product of reflections is directly related to 
the angle between the corresponding hyperplanes
leads to the definition of a Coxeter group.

\begin{definition}
{\rm
A \emph{Coxeter system} $(W,S)$ is a pair consisting of a group $W$ 
together with a finite subset $S\subset W$ satisfying the following
conditions:
\begin{itemize}
\item[(i)]
each $s\in S$ is an involution: $s^2=1$;
\item[(ii)]
some pairs $\{s,t\}\subset S$ satisfy relations of the form
$(st)^{m_{st}}\!=\!1$ with \hbox{$m_{st}\geq 2$;} 
\item[(iii)]
the relations in (i)--(ii) form a presentation of the group~$W$. 
\end{itemize}
In other words, $S$ generates~$W$, and any identity in~$W$ is a formal
consequence of (i)--(ii) and the axioms of a group. 

A group $W$ is called a \emph{Coxeter group} if it has a presentation
of the above form. 
}
\end{definition}

The following theorem demonstrates that the notion of a 
Coxeter group indeed captures the geometric essence 
of reflection groups.

\begin{theorem}
\label{Cox ref root}
Any finite Coxeter group is isomorphic to a reflection group. 
\end{theorem}

Conversely, a reflection group associated with a (finite
crystallographic) root system~$\Phi$ is a Coxeter group,
in the following sense. 
Let $\Pi$ be the set of simple roots in~$\Phi$. 
For each simple root $\alpha_i\in\Pi$, 
the associated {\em simple reflection} is  $s_i\stackrel{\rm
  def}{=}\sigma_{\alpha_i}$. 

\begin{theorem}
\label{th:refgp-is-coxeter}
Let $W$ be the group generated by the 
reflections $\set{\sigma_\beta:\beta\in\Phi}$. 
Let 
\[
S={\set{s_i}}_{i\in I}
=\set{\sigma_\alpha:\alpha\in\Pi}
\]  
be the set of simple reflections. 
Then $(W,S)$ is a Coxeter system. 
\end{theorem}


Furthermore, $W$ is a \emph{crystallographic 
Coxeter group}, where the adjective ``crystallographic'' refers to
restricting the integers $m_{st}$ to the set $\set{2,3,4,6}$.

\section{Other ``finite type'' classifications}
\label{sec:other}

The classification of root systems is similar or identical to several other
classifications of objects of ``finite type,'' 
briefly reviewed below. 

\subsection*{Non-crystallographic root systems}

Lifting the crystallographic restriction does not allow very many
additional root systems. 
The only non-crystallographic irreducible finite root systems are those of
types $H_3$, $H_4$ and $I_2(m)$ for $m=5$ or $m\ge 7$.
See~\cite{Humphreys}.

\subsection*{Coxeter groups and reflection groups}

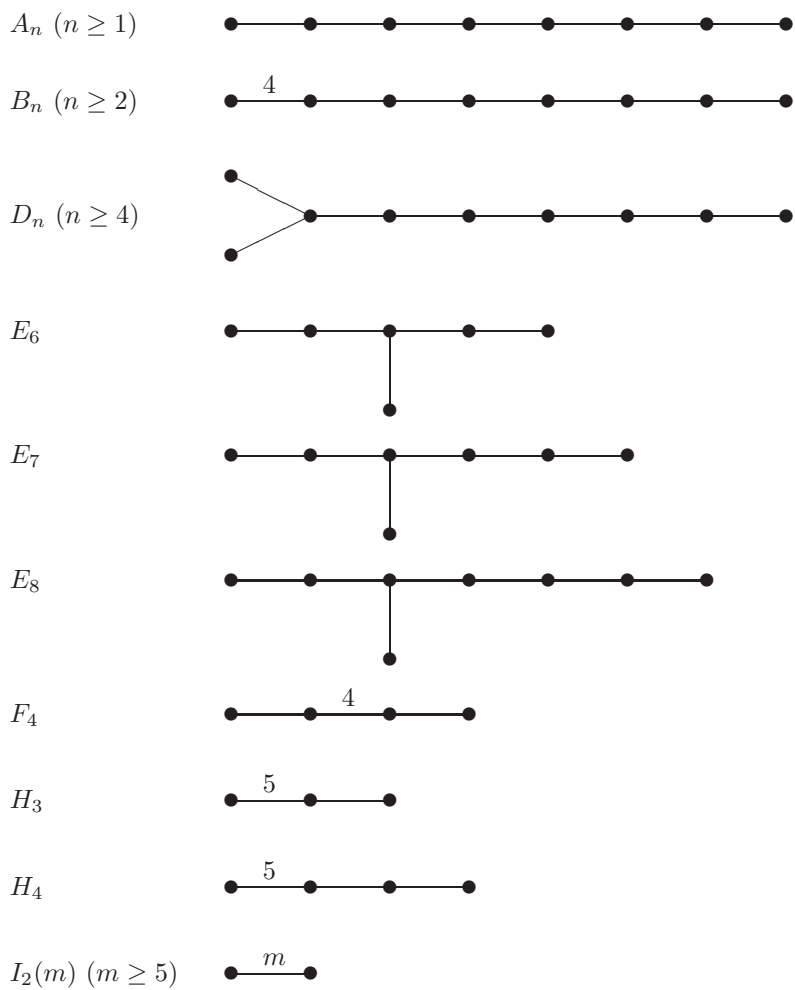
\begin{figure}[ht]
\vspace{-.2in} 
\[ 
\begin{array}{lcl} 
A_n\  (n\geq 1) && 
\setlength{\unitlength}{1.5pt} 
\begin{picture}(140,17)(0,-2) 
\put(0,-0.15){\line(1,0){140}} 
\multiput(0,0)(20,0){8}{\circle*{3}} 
\end{picture}\\
B_n\ (n\geq 2)
&& 
\setlength{\unitlength}{1.5pt} 
\begin{picture}(140,17)(0,-2) 
\put(0,-0.15){\line(1,0){140}}
\multiput(0,0)(20,0){8}{\circle*{3}} 
\put(8,2){$4$}
\end{picture} \\[.2in] 
D_n\ (n\geq 4)
&& 
\setlength{\unitlength}{1.5pt} 
\begin{picture}(140,17)(0,-2) 
\put(20,-0.15){\line(1,0){120}} 
\put(0,10){\line(2,-1){20}} 
\put(0,-10){\line(2,1){20}} 
\multiput(20,0)(20,0){7}{\circle*{3}} 
\put(0,10){\circle*{3}} 
\put(0,-10){\circle*{3}} 
\end{picture} 
\\[.2in] 
E_6 
&& 
\setlength{\unitlength}{1.5pt} 
\begin{picture}(140,17)(0,-2) 
\put(0,-0.15){\line(1,0){80}} 
\put(40,0){\line(0,-1){20}} 
\put(40,-20){\circle*{3}} 
\multiput(0,0)(20,0){5}{\circle*{3}} 
\end{picture} 
\\[.25in] 
E_7 
&& 
\setlength{\unitlength}{1.5pt} 
\begin{picture}(140,17)(0,-2) 
\put(0,-0.15){\line(1,0){100}} 
\put(40,0){\line(0,-1){20}} 
\put(40,-20){\circle*{3}} 
\multiput(0,0)(20,0){6}{\circle*{3}} 
\end{picture} 
\\[.25in] 
E_8 
&& 
\setlength{\unitlength}{1.5pt} 
\begin{picture}(140,17)(0,-2) 
\put(0,-0.15){\line(1,0){120}} 
\put(40,0){\line(0,-1){20}} 
\put(40,-20){\circle*{3}} 
\multiput(0,0)(20,0){7}{\circle*{3}} 
\end{picture} 
\\[.3in] 
F_4 
&& 
\setlength{\unitlength}{1.5pt} 
\begin{picture}(140,17)(0,-2) 
\put(0,-0.15){\line(1,0){60}} 
\multiput(0,0)(20,0){4}{\circle*{3}} 
\put(28,2){$4$}
\end{picture} 
\\[.05in] 
H_3
&&
\setlength{\unitlength}{1.5pt} 
\begin{picture}(140,17)(0,-2) 
\put(0,-0.15){\line(1,0){40}} 
\multiput(0,0)(20,0){3}{\circle*{3}} 
\put(8,2){$5$}
\end{picture} 
\\[.05in] 
H_4
&&
\setlength{\unitlength}{1.5pt} 
\begin{picture}(140,17)(0,-2) 
\put(0,-0.15){\line(1,0){60}} 
\multiput(0,0)(20,0){4}{\circle*{3}} 
\put(8,2){$5$}
\end{picture} 
\\[.05in] 
I_2(m)\  (m\geq 5) 
&&
\setlength{\unitlength}{1.5pt} 
\begin{picture}(140,17)(0,-2) 
\put(0,-0.15){\line(1,0){20}} 
\multiput(0,0)(20,0){2}{\circle*{3}} 
\put(8,2){$m$}
\end{picture} 
\end{array} 
\] 
\vspace{-.1in} 
\caption{Coxeter diagrams of finite irreducible Coxeter systems}
\label{Cox diagrams} 
\end{figure}

By Theorems~\ref{Cox ref root} and~\ref{th:refgp-is-coxeter}, 
the classification of finite Coxeter
groups is parallel to the classification of reflection groups and is 
essentially the same as the classification of root systems.
The difference is that the root systems $B_n$ and $C_n$ correspond to
the same Coxeter group~$B_n$. 
A Coxeter group is encoded by its \emph{Coxeter diagram}, a graph 
whose vertex set is $S$, with an edge $s$---$t$ whenever $m_{st}>2$.
If $m_{st}>3$, the edge is labeled by $m_{st}$.
Figure~\ref{Cox diagrams} shows the Coxeter diagrams of the finite
irreducible Coxeter systems, 
including the non-crystallographic Coxeter groups
$H_3$, $H_4$ and~$I_2(m)$.  The group $G_2$ appears as $I_2(6)$.
See~\cite{Humphreys} for more details.

\pagebreak

\subsection*{Regular polytopes}

By Theorem~\ref{regular}, the symmetry group of a regular polytope is
a reflection group. 
In fact, it is a Coxeter group whose Coxeter diagram 
is \emph{linear}:
the underlying graph is a path with no branching points.
This narrows down the possibilities, leading to the conclusion that
there are no other regular polytopes besides the ones described in 
Section~\ref{ref gp}.
In particular, there are no ``exceptional'' regular polytopes beyond 
dimension~4: 
only simplices, cubes, and crosspolytopes.
See~\cite{Coxeter}. 

\subsection*{Lie algebras}

The original motivation for the Cartan-Killing classification of root
systems came from Lie theory. 
Complex finite-dimensional simple Lie algebras correspond naturally,
and one-to-one, to finite irreducible crystallographic root systems.
There exist innumerable expositions of this classical subject; 
see, e.g.,~\cite{Fu-Ha}. 

\subsection*{Quivers of finite type
}

A \emph{quiver} is a directed graph; 
its {\em representation} assigns a complex vector
space to each vertex, and a linear map to each directed edge.
A quiver is \emph{of finite type} if it has only a finite
number of indecomposable representations (up to isomorphism);
a representation is \emph{indecomposable} if it cannot be obtained as
a nontrivial direct~sum. 
By Gabriel's Theorem, a quiver is of finite type if and only if its
underlying graph is a  
Dynkin diagram of type $A$, $D$ or~$E$. 
See \cite{reiten-dynkin} and references therein. 




\subsection*{Et cetera} 

And the list goes on: simple singularities, finite subgroups
of~$SU(2)$,
symmetric matrices with nonnegative integer entries and eigenvalues
between $-2$ and~$2$, \emph{etc}. 
For more, see \cite{Geck-Malle,hhsv, zuber}.
In Section~\ref{sec:clust-fintype}, we will present yet another
classification that is parallel to Cartan-Killing:
the classification of the cluster algebras of finite type.

\section{Reduced words and permutohedra} 

Each element $w\in W$ can be written as a product of elements of~$S$:
\[
w=s_{i_1}\cdots s_{i_\ell}\,. 
\]
A shortest factorization of this form 
(or the corresponding sequence of subscripts $(i_1,\dots,i_\ell)$) 
is called a \emph{reduced word} for~$w$;
the number of factors $\ell$ is called the
\emph{length} of~$w$. 

Any finite Coxeter group has a unique element $w_\circ$ of maximal length.
In the symmetric group $\mathcal{S}_{n+1}=A_{n}$, this is the permutation
$w_\circ$ that reverses the order of the elements of the set 
$\set{1,\ldots,n+1}$.

\begin{example}
\label{3412}
\rm
Let  $W=\mathcal{S}_4$ be the Coxeter group of type~$A_3$.
The standard choice of simple reflections yields
$S=\{s_1,s_2,s_3\}$, where  
$s_1$, $s_2$ and $s_3$
are the transpositions which interchange
$1$~with~$2$, $2$~with~$3$, and $3$~with~$4$, respectively.
(Cf.\ Example~\ref{An}.) 

The word 
$s_1s_2s_1s_3s_2s_3$ is a non-reduced word for the
permutation that interchanges $1$ with~$3$ and $2$ with~$4$. 
This permutation has two reduced words $s_2s_1s_3s_2$ and $s_2s_3s_1s_2$.

An example of a reduced word for $w_\circ$ is $s_1s_2s_1s_3s_2s_1$.
There are 16 such reduced words altogether.
(Cf.\ Example~\ref{permut23} and Theorem~\ref{numred}.) 
\end{example}

Recall from Section~\ref{ref gp} that we label the regions~$R_w$ of
the Coxeter arrangement by the elements of the reflection group~$W$,
so that $R_w$ is the image of $R_1$ under the action of~$w$.
More generally, $R_{uv}=u(R_v)$. 

\begin{lemma}
\label{lem:adj-regions}
In the Coxeter arrangement associated with a reflection group~$W$, 
regions $R_u$ and $R_v$ are adjacent
(that is, share a codimension~$1$ face)
if and only if $u^{-1}v$ is a simple reflection. 
\end{lemma}

Thus, moving to an adjacent region is encoded by multiplying \emph{on the
right} by a simple reflection; cf.\ Figure~\ref{I25}. 
(Warning:  this simple reflection is generally \emph{not} the same as the
reflection through the hyperplane separating the two adjacent regions.)
Consequently, reduced words for an element $w\in W$ correspond 
to equivalence classes of paths from $R_1$ to $R_w$ in the ambient
space of the Coxeter arrangement. 
More precisely, we consider the paths that cross hyperplanes of 
the arrangement one at a time, and cross each hyperplane at most once; 
two paths are equivalent if they cross the same hyperplanes
in the same order.

In order to make the correspondence between paths and reduced words 
more explicit, one can restrict the paths to the edges of the 
{\em $W$-permutohedron}, a convex polytope that we will now define.
Fix a point $x$ in the interior of~$R_1$.  
The $W$-permutohedron is the convex hull 
of the orbit of $x$ under the action of~$W$.
The name ``permutohedron'' comes from the fact that the vertices of an 
$A_{n}$-permutohedron are obtained by permuting the coordinates of a
generic point in~$\reals^{n+1}$. 

\begin{example}
\label{permut23}
{\rm
The $A_2$, $B_2$ and $G_2$ permutohedra are respectively a hexagon, an 
octagon and a dodecagon; under the right choices of~$x$, these
polygons are regular.
Figures~\ref{A3perm} and~\ref{B3perm} show the  permutohedra of types $A_3$ and~$B_3$.
Each of these realizations derives from a choice of $x\in R_1$ which
makes the permutohedron an Archimedean solid, so that in particular
its facets are all regular polygons.
The non-crystallographic $H_3$-permutohedron is also an Archimedean
solid\footnote{An \emph{Archimedean solid} is a non-regular 
polytope whose all facets are regular polygons, 
and whose symmetry group acts transitively on vertices. 
In dimension~$3$, 
there are $13$ Archimedean solids. 
The permutohedra of types $A_3$, $B_3$, and~$H_3$ are also known as
the \emph{truncated octahedron}, \emph{great rhombicuboctahedron},
and \emph{great rhombicosidodecahedron},
respectively. See, e.g.,~\cite{wolfram}. 
}.

\begin{figure}[ht]
\centerline{
        \epsfbox{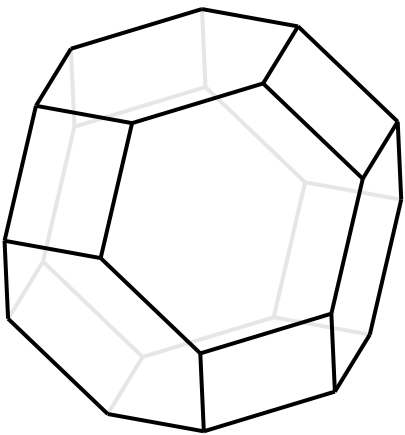}
}
\caption{The permutohedron of type $A_3$
}
\label{A3perm}
\end{figure}

\begin{figure}[ht]
\centerline{
\epsfbox{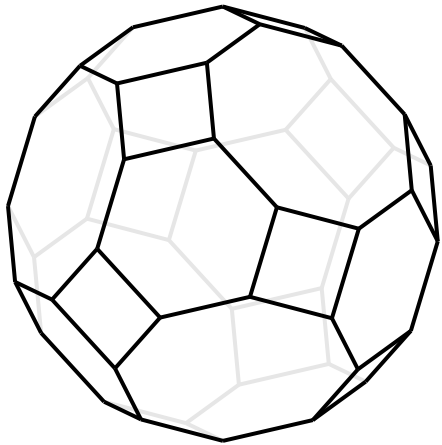}
}
\caption{The permutohedron of type $B_3$
}
\label{B3perm}
\end{figure}

In both pictures, the bottom vertex can be associated with
the identity element~$1\!\in\! W$, so that the top vertex is~$w_\circ$.
A reduced word for $w$ corresponds to a path along edges from $1$ to $w$ 
which moves up in a monotone fashion.
There are 16 such paths from $1$ to $w_\circ$ in the
$A_3$-permutohedron; cf.\ Example~\ref{3412}.  
}
\end{example}

The following beautiful formula is due to R.~Stanley~\cite{stanley-84}.

\begin{theorem}
\label{numred}
The number of reduced words for $w_\circ$ in the reflection group $A_n$ is 
\[\frac{\binom{n+1}{2}!}{1^n3^{n-1}5^{n-2}\cdots (2n-1)^1} .\]
\end{theorem}


\section{Coxeter element and Coxeter number}
\label{sec:Coxeter element}

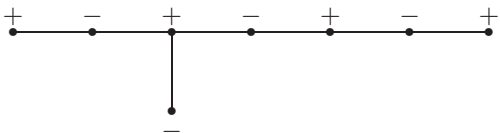
\begin{figure}[ht]
\begin{center}
\setlength{\unitlength}{1.5pt}
\begin{picture}(140,24)(-10,-18)
\put(0,0){\line(1,0){120}}
\put(40,0){\line(0,-1){20}}
\put(40,-20){\circle*{2}}
\multiput(0,0)(20,0){7}{\circle*{2}}
\put(0,4){\makebox(0,0){$+$}}
\put(40,-25){\makebox(0,0){$-$}}
\put(20,4){\makebox(0,0){$-$}}
\put(40,4){\makebox(0,0){$+$}}
\put(60,4){\makebox(0,0){$-$}}
\put(80,4){\makebox(0,0){$+$}}
\put(100,4){\makebox(0,0){$-$}}
\put(120,4){\makebox(0,0){$+$}}
\end{picture}
\end{center}
\caption{Bi-partition of the nodes of the Coxeter diagram of type~$E_8$}
\label{fig:I+I-}
\end{figure}

The underlying graph of the Coxeter diagram for a finite Coxeter group has no cycles.
Hence it is bipartite, i.e., we can write a disjoint union 
$I=I_+\cup I_-$ such that each of the sets $I_+$ and $I_-$ is totally
disconnected in the Coxeter diagram. 
An example is shown in Figure~\ref{fig:I+I-},
where the elements of $I_+$ and $I_-$ are marked by $+$ and~$-$,
respectively. 

The simple reflections associated with~$I_+$
(resp.,~$I_-$) commute pairwise.
Consequently, the following is well-defined: 
\[
c=\Biggl(\prod_{i\in I_+}s_i\Biggr) \Biggl(\prod_{i\in I_-}s_i\Biggr)\,. 
\]
The element $c\in W$ is called the {\em Coxeter element}\footnote{
\label{note:Coxeter element}
More broadly, one often calls 
the product of the elements in~$S$ (in any order) 
a Coxeter element, 
but for our present purposes the definition above will do.}.

\begin{example}\rm
In type $A_n$, let $I_-$ (resp.,~$I_+$) consist of the odd (resp.,
even) numbers in $I=[n]$. 
Then for example in $A_5=\mathcal{S}_6$, we have 
$c=s_2s_4s_1s_3s_5$. 
\end{example}

Thinking of $W$ as a reflection group, the Coxeter element $c$ is an 
interesting orthogonal transformation.
One important feature of $c$ is that it fixes a certain
two-dimensional plane~$L$  (as a set, not pointwise).
The action of $c$ on $L$ can be analyzed to determine the 
order of $c$ as an element of~$W$.
This order is called the {\em Coxeter number} of~$W$,
and is denoted by~$h$.

\begin{figure}[ht]
\centerline{
        \epsfbox{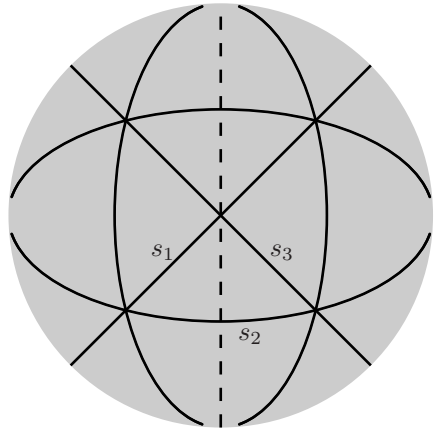}
        \begin{picture}(0,0)(160,0)
                \put(83,33){$s_2$}
                \put(95,65){$s_3$}
                \put(50,65){$s_1$}
        \end{picture}
}
\caption{The Coxeter arrangement~$A_3$ and the plane fixed by the
  Coxeter element
}
\label{A3plane}
\end{figure}

\begin{example}\rm
Figure~\ref{A3plane} shows the Coxeter 
arrangement of type~$A_3$ and the plane~$L$ fixed by the
  Coxeter element $c=s_2s_1s_3$ (dotted). 
The great circles represent the 
intersections of the six reflecting hyperplanes with a unit hemisphere.
The sphere is opaque, so only half of each circle is visible,
and appears either as a \emph{half} of an ellipse 
or as a straight line segment.
(The ``equator'' does not represent a hyperplane in the arrangement.) 
The restriction of $c$ onto~$L$ has order~$4$, so the Coxeter number for $A_3$
is~$h=4$.
\end{example}

\begin{figure}[ht]
\centerline{
        \epsfbox{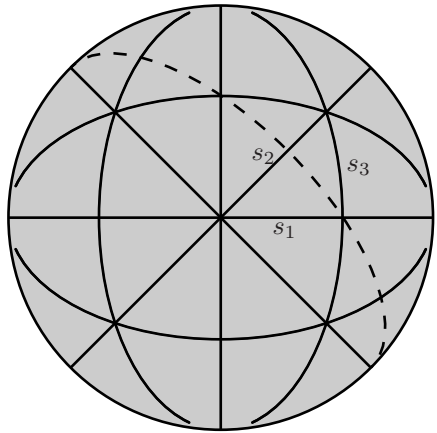}
        \begin{picture}(0,0)(160,0)
                \put(88,102){$s_2$}
                \put(96,74){$s_1$}
                \put(124,98){$s_3$}
        \end{picture}
}
\caption{The Coxeter arrangement $B_3$ and the plane fixed by the
  Coxeter element
}
\label{B3plane}
\end{figure}

\begin{example}\rm
Figure~\ref{B3plane} is a similar picture for $B_3$, illustrating that the 
Coxeter number for $B_3$ is $h=6$.
In this picture, the equator does represent a hyperplane in the 
arrangement.
\end{example}

\pagebreak

The action of $c$ on $L$ also leads to a determination of its 
eigenvalues, which all have the form 
$e^{2mi\pi/h}$, where $m$ is a positive integer less than~$h$. 
The $n$ values of $m$ which arise in this way 
are called the {\em exponents} of~$W$.
We denote the exponents by $e_1,\dots,e_n$. 
They pop up everywhere in the combinatorics of root systems and 
Coxeter groups.
For instance, the order (i.e., cardinality) of $W$ is expressed in terms of the
exponents by
\[
|W|=\prod_{i=1}^n (e_i+1)\,. 
\]
See Section~\ref{sec:numerology} for more examples. 

For a finite irreducible Coxeter group~$W$
Figure~\ref{data} tabulates some classical numerical invariants
associated to $W$ and the corresponding (not necessarily crystallographic) root system~$\Phi$.

\begin{figure}[ht]
\ \bigskip

\begin{tabular}{|c|c|c|c|c|}
\hline
type of $\Phi$&$|\Phi_+|$&$h$& $e_1,\dots,e_n$ &$|W|$\\\hline&&&&\\[-3.5mm]
\hline
$A_n$&$n(n+1)/2$&$n+1$&$1,2,\ldots,n$&$(n+1)!$\\\hline
$B_n, C_n$&$n^2$&$2n$&$1,3,5,\ldots,2n-1$&$2^nn!$\\\hline
$D_n$&$n(n-1)$&$2(n-1)$&$1,3,5,\ldots,2n-3,n-1$&$2^{n-1}n!$\\\hline
$E_6$&$36$&$12$&$1,4,5,7,8,11$&$2^73^45$\\\hline
$E_7$&$63$&$18$&$1,5,7,9,11,13,17$&$2^{10}3^45\cdot 7$\\\hline
$E_8$&$120$&$30$&$1,7,11,13,17,19,23,29$&$2^{14}3^55^27$\\\hline
$F_4$&$24$&$12$&$1,5,7,11$&$2^73^2$\\\hline
$G_2$&$6$&$6$&$1,5$&$2^23$\\\hline
$H_3$&$15$&$10$&$1,5,9$&$2^33\cdot 5$\\\hline
$H_4$&$60$&$30$&$1,11,19,29$&$2^63^25^2$\\\hline
$I_2(m)$&$m$&$m$&$1,m\!-\!1$&$2m$\\\hline
\end{tabular}

\medskip

\caption{
\hbox{Number of positive roots, Coxeter number, 
exponents, and the order of~$W$.} 
}
\label{data}
\end{figure}

\lecture{Associahedra and Mutations}
\label{lec3}

\section{Associahedron}

We start by discussing two classical problems of combinatorial
enumeration. 
\begin{itemize}
\item[(i)]
Count the number of \emph{bracketings}
 (\emph{parenthesizations}) of a non-associative product of
$n+2$ factors.
Note that we need $n$ pairs of brackets in order to 
make the product unambiguous. 
\item[(ii)]
Count the number of \emph{triangulations} of a convex
$(n\!+\!3)$-gon by diagonals. 
Note that each triangulation involves exactly $n$ diagonals. 
\end{itemize}

\begin{example}\rm
In the special cases $n=1,2,3$,
there are, respectively: 
\begin{itemize}
\item
$2$ bracketings $(ab)c$ and $a(bc)$ of a product of $3$ factors; 
\item
$5$ bracketings $((ab)c)d$, $(a(bc))d$, $a((bc)d)$, $(ab)(cd)$, 
and $a(b(cd))$ of a product of $4$ factors; 
\item
$14$ bracketings of a product of $5$ factors (check!).  
\end{itemize}
As to triangulations, there are: 
\begin{itemize}
\item
$2$ triangulations of a convex quadrilateral ($n=1$);
\item
$5$ triangulations of a pentagon ($n=2$, Figure~\ref{A2assoc_basic});
\item
$14$ triangulations of a hexagon ($n=3$, Figure~\ref{A3assoc}). 
\end{itemize}
\end{example}

\begin{theorem}
\label{th:two-catalan}
Both bracketings and triangulations described above are 
enumerated by the Catalan numbers 
$\frac{1}{n+2}\binom{2n+2}{n+1}$. 
\end{theorem}

There are a great many families of combinatorial objects enumerated by
the Catalan numbers;
more than a hundred of those are listed in \cite{ec2-catalan}. 
This list includes: ballot sequences; Young diagrams and tableaux
satisfying certain restrictions; noncrossing partitions; trees of
various kinds; Dyck paths; permutations avoiding patterns of
length~$3$; and much more. 
In Lecture~\ref{lec:num}, we will discuss several additional members
of the ``Catalan family,'' 
together with their analogues for arbitrary root systems.
(We will see that the ordinary Catalan numerology should be considered
as ``type~$A$.'') 

\medskip

\begin{figure}[ht]
\centerline{
\scalebox{0.666}{
        \epsfbox{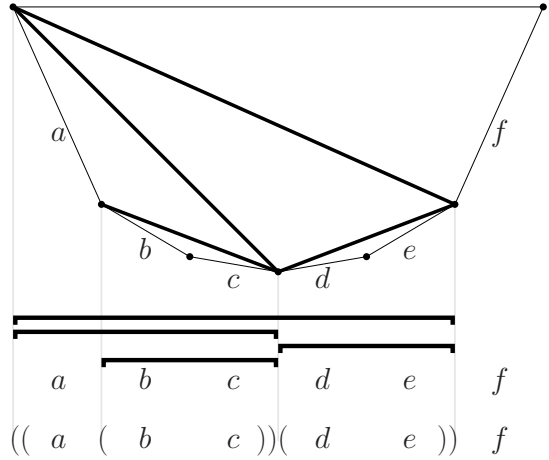}   
        \begin{picture}(0,0)(150,-250)
                \put(-133,-75){\huge$a$}
                \put(-83,-142){\huge$b$}
                \put(-33,-159){\huge$c$}
                \put(17,-159){\huge$d$} 
                \put(67,-142){\huge$e$} 
                \put(117,-75){\huge$f$} 
                \put(-133,-216){\huge$a$}
                \put(-83,-216){\huge$b$}
                \put(-33,-216){\huge$c$}
                \put(17,-216){\huge$d$} 
                \put(67,-216){\huge$e$} 
                \put(117,-216){\huge$f$}        
                \put(-133,-250){\huge$a$}
                \put(-83,-250){\huge$b$}
                \put(-33,-250){\huge$c$}
                \put(17,-250){\huge$d$} 
                \put(67,-250){\huge$e$} 
                \put(117,-250){\huge$f$}
                \put(-156.5,-250){\huge$(($}
                \put(-106.5,-250){\huge$($}
                \put(-16.25,-250){\huge$))($}
                \put(86.5,-250){\huge$))$}
        \end{picture}
}
}
\caption{The bijection between triangulations and bracketings.}
\label{triparen}
\end{figure}

A bijection between bracketings and triangulations is described 
in Figure~\ref{triparen}.

\begin{figure}[ht]
\centerline{
\scalebox{0.85}{
        \epsfbox{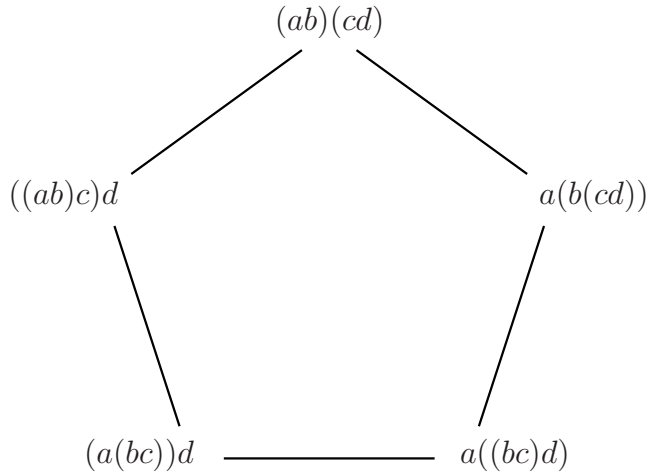}
        \begin{picture}(0,0)(104,-92)
                \setlength{\unitlength}{.75pt} 
                \put(-32,141){\LARGE$(ab)(cd)$}
                \put(125,38){\LARGE$a(b(cd))$}
                \put(78,-116){\LARGE$a((bc)d)$}
                \put(-145,-116){\LARGE$(a(bc))d$}
                \put(-190,38){\LARGE$((ab)c)d$}
        \end{picture}
}
}
\caption{Applying associativity to the bracketings of $abcd$.}
\label{A2paren}
\end{figure}

\begin{figure}[ht]
\centerline{
\scalebox{0.85}{
        \epsfbox{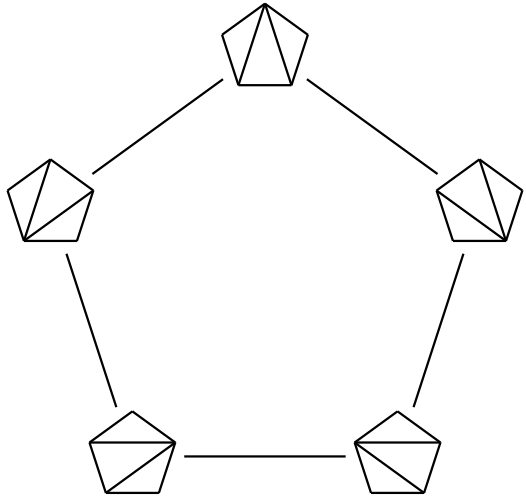}
}
}
\caption{The exchange graph for triangulations of a pentagon.}
\label{A2assoc_basic}
\end{figure}

\begin{figure}[ht!]
\centerline{
\scalebox{.87}{
\epsfbox{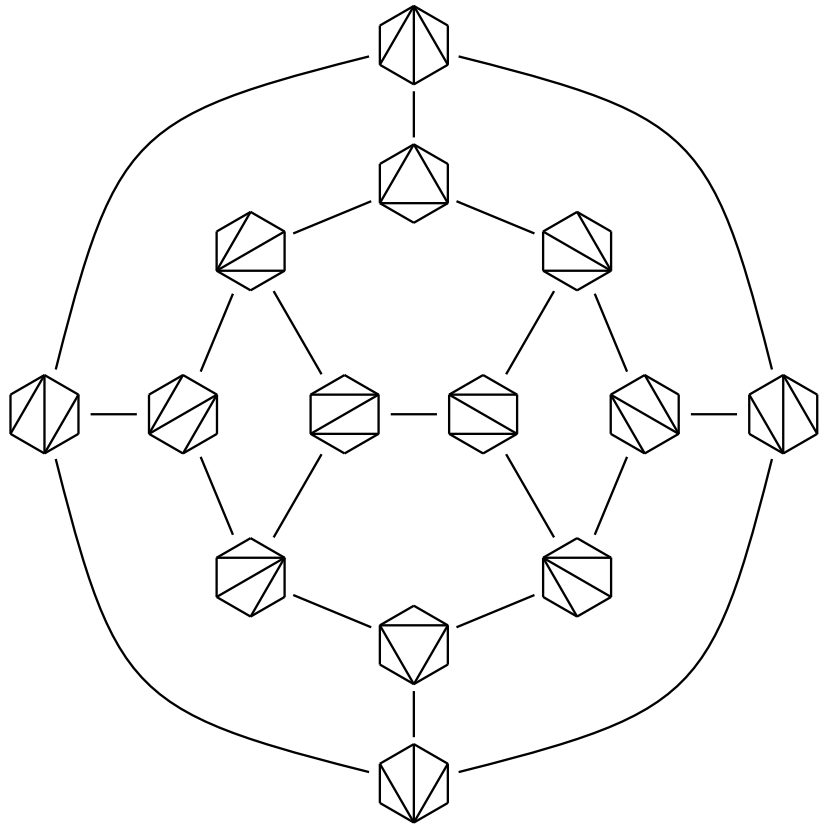}}
}
\caption{The exchange graph for triangulations of a hexagon.}
\label{A3assoc}
\end{figure}

For a fixed~$n$, the bracketings naturally form the set of vertices
of a graph whose edges correspond to applications of the 
associativity axiom.
Figure~\ref{A2paren} shows this graph for~$n=2$.

The bijection illustrated in Figure~\ref{triparen}
translates an 
application of the associativity axiom into a 
diagonal {\em flip} on the corresponding triangulation.
That is, one removes a diagonal to create a quadrilateral, then 
replaces the removed diagonal with the other diagonal of the
quadrilateral. 

We call the graph defined by diagonal flips the {\em exchange graph}.
The exchange graphs for $n=2$ and $n=3$ are shown in
Figures~\ref{A2assoc_basic} and~\ref{A3assoc}.

\bigskip

\begin{figure}[ht]
\centerline{\epsfbox{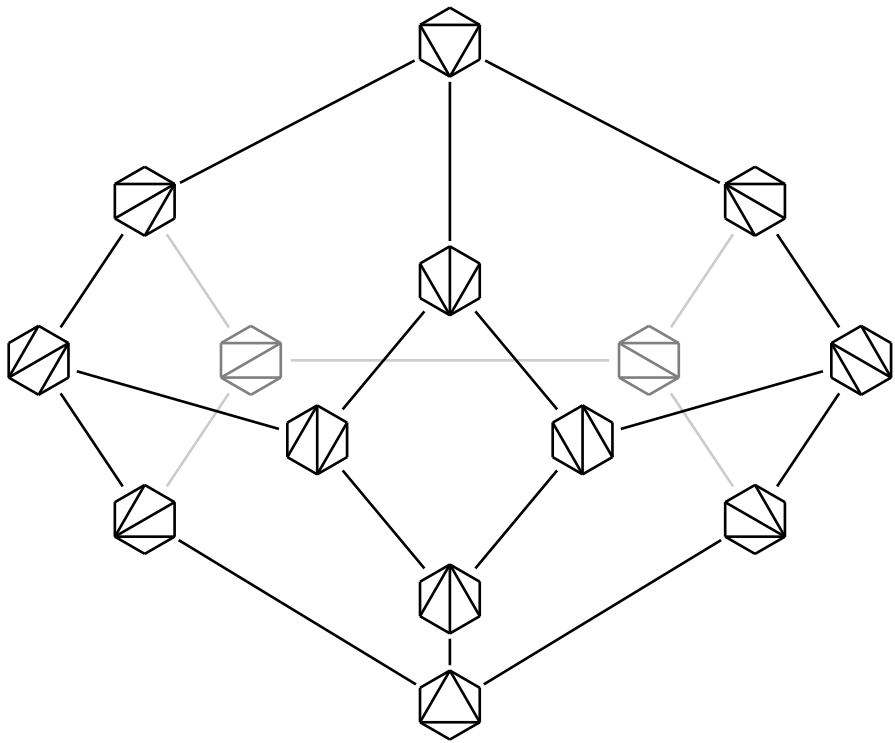}}
\caption{The 3-dimensional associahedron.}
\label{A3assoc_poly}
\end{figure}

The drawing of the exchange graph in Figure~\ref{A3assoc} 
fails to convey its crucial property: 
this exchange graph is the $1$-skeleton of a convex polytope, 
the \hbox{$3$-dimensional} \emph{associahedron}.  
(Sometimes it is also called the {\em Stasheff polytope}, 
after J.~Stasheff, who first defined it in~\cite{stasheff}.) 
Figure~\ref{A3assoc_poly} shows a polytopal realization of
this associahedron. 

In order to formally define the $n$-dimensional associahedron,
we start by describing the object which is dual to it, 
in the same sense in which the octahedron is dual to the cube,
and the dodecahedron is dual to the icosahedron. 

\begin{definition}[The dual complex of an associahedron]
\label{def:asso-dual}
{\rm
Consider the following 
simplicial complex:

\vspace{.05in}

\begin{center}

\begin{tabular}{|rl|}
\hline
&\\[-.13in]
vertices: & diagonals of a convex $(n\!+\!3)$-gon\\[.05in]
simplices: & partial triangulations of the $(n\!+\!3)$-gon\\
           &  (viewed as collections of non-crossing diagonals)\\[.05in]
maximal simplices: & triangulations of the $(n\!+\!3)$-gon\\
         & (collections of $n$ non-crossing diagonals).\\[.025in]
\hline
\end{tabular}

\vspace{.05in}

\end{center}
}\end{definition}

\begin{figure}[ht]
\centerline{
        \epsfbox{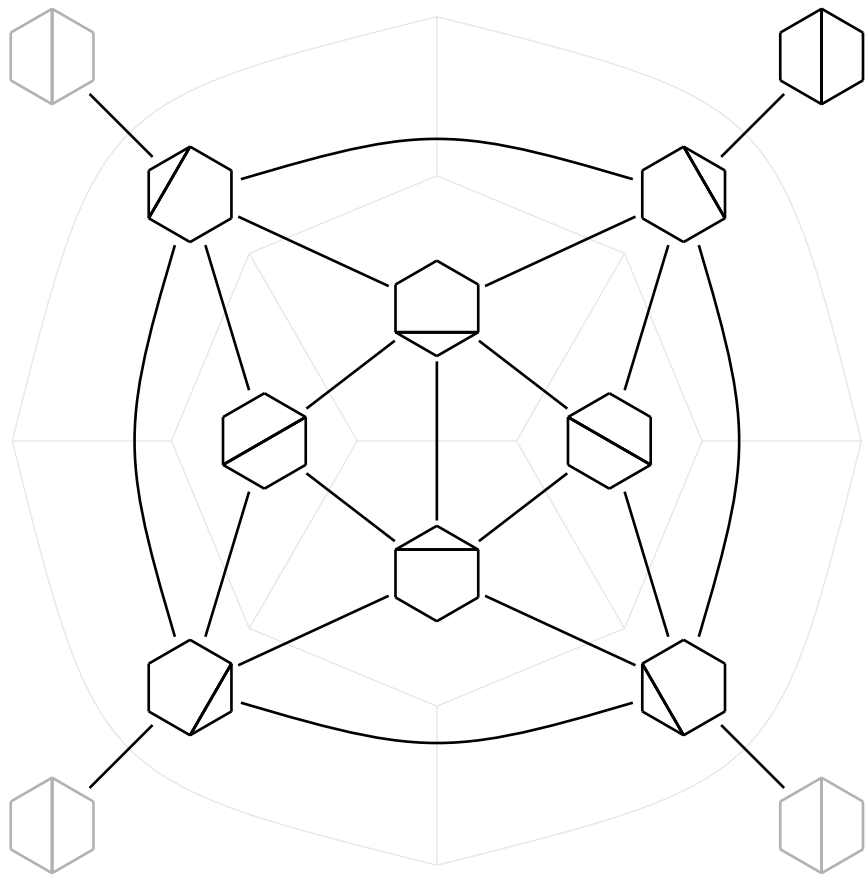}
        \begin{picture}(0,0)(230,-230)
        \end{picture}
}
\caption{
The simplicial complex dual to the $3$-dimensional associahedron.
}
\label{A3assoc_dual}
\end{figure}

Figure~\ref{A3assoc_dual} shows this simplicial complex for $n=3$, superimposed 
on a faint copy of the exchange graph.
Note that the facial structures of the $3$-dimensional associahedron and its dual
complex are indeed ``dual'' to each other:
two vertices of one polyhedron are adjacent if and only if the
corresponding faces of the other polyhedron share an edge. 

\pagebreak

It is not clear \emph{a priori} that these complexes are topological spheres.
But, as already mentioned, more is true.

\begin{theorem}
\label{dual-assoc}
The simplicial complex described in Definition~\ref{def:asso-dual} can
be realized as the boundary of an $n$-dimensional convex polytope. 
\end{theorem}

Theorem~\ref{dual-assoc} (or its equivalent reformulations) 
were proved independently by 
J.~Milnor, M.~Haiman, and C.~W.~Lee (first published proof~\cite{lee}).
This theorem also follows as a special case of
the very general theory of secondary polytopes
developed by I.~M.~Gelfand, M.~Kapranov and A.~Zelevinsky~\cite{gkz}.

\begin{definition}[The associahedron]
\label{def:asso}
{\rm
The $n$-dimensional \emph{associahedron} is the convex polytope 
(defined up to combinatorial equivalence) 
that is dual (or polar, see
\cite[Sec.~2.3]{Ziegler}) to the polytope of
Theorem~\ref{dual-assoc}. 
}\end{definition}

\begin{figure}[ht]
\centerline{
        \epsfbox{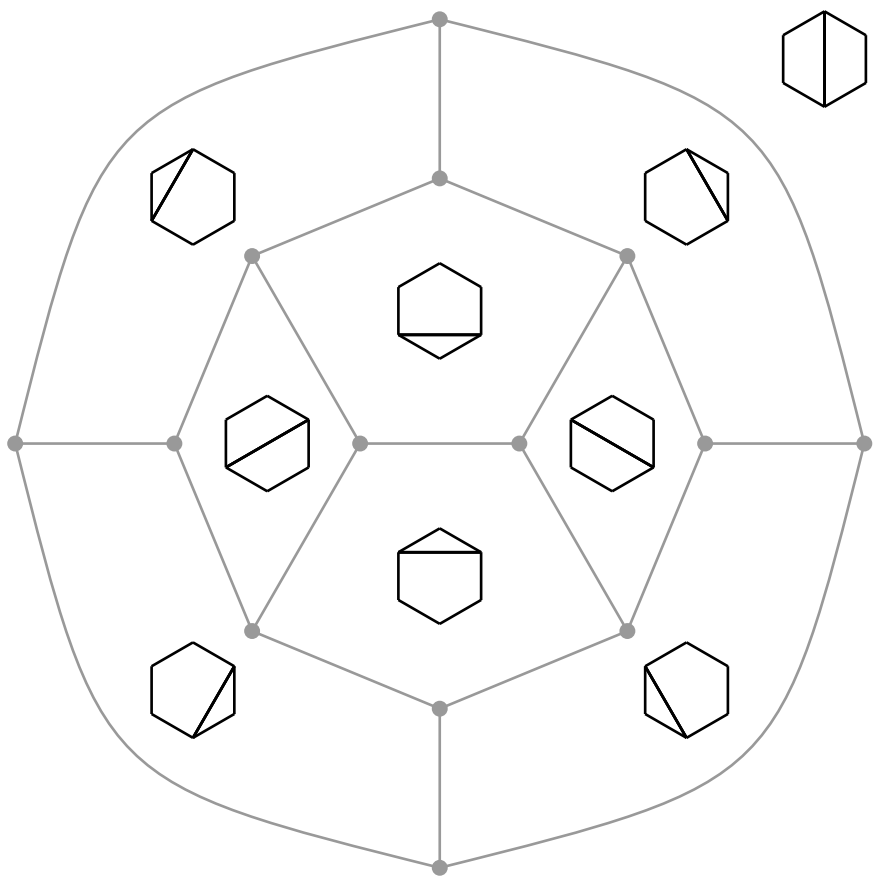}
        \begin{picture}(0,0)(230,-230)
        \end{picture}
}
\caption{Labeling the facets of the 
associahedron by diagonals 
}
\label{A3assoc_faces}
\end{figure}

The facial structure of an associahedron 
as a cell complex is dual to that of its polar: 
\begin{equation}
\label{eq:faces-asso}
\text{
\begin{tabular}{|rl|}
\hline
&\\[-.13in]
vertices: & triangulations\\[.05in]
faces: & partial triangulations\\[.05in]
facets:   & diagonals\\[.05in]
edges:  & diagonal flips\\[.025in]
\hline
\end{tabular}
}
\end{equation}
The labeling of the facets of an $n$-dimensional associahedron by the
diagonals of an $(n+3)$-gon is illustrated in 
Figure~\ref{A3assoc_faces} for the special case $n=3$ 
(compare to Figure~\ref{A3assoc}).

We note that we could have defined the associahedron directly, 
as a cell complex whose cell structure is described 
by~\eqref{eq:faces-asso}. 
(This would require resolving some technical issues that we would
rather avoid here.)  
The fact that these cell complexes are polytopal---i.e., 
the fact that a combinatorially defined associahedron 
can be realized as a convex polytope---is essentially equivalent to
Theorem~\ref{dual-assoc}. 

Associahedra play an important role in homotopy theory and the study of
operads~\cite{stasheff-notices}, 
in the analysis of real moduli/configuration spaces~\cite{devadoss-notices}, 
and other branches of mathematics. 
In these notes, we restrict our attention to the combinatorial
aspects of the associahedra. 

An $n$-dimensional polytope is called \emph{simple} if every vertex is
incident to exactly $n$ edges.
This is the case for the associahedron,
as every triangulation 
of an $(n\!+\!3)$-gon is adjacent to precisely $n$ others in the
exchange graph.

Much is known about the facial structure and enumerative invariants of the 
associahedron.
For example, each face is the direct product of smaller associahedra. 
The entries of the \emph{$h$-vector} of the associahedron are the Narayana numbers 
(see Section~\ref{sec:narayana}).
This allows one to calculate the number of faces of each dimension.

\section{Cyclohedron}

\begin{figure}[ht]
\centerline{\scalebox{0.6}{\epsfbox{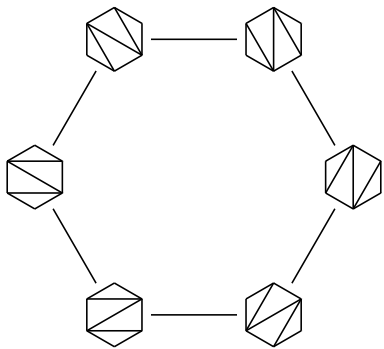}}}
\caption{The $2$-dimensional cyclohedron}
\label{B2assoc}
\end{figure}

\begin{figure}[ht]
\centerline{\scalebox{0.9}{\epsfbox{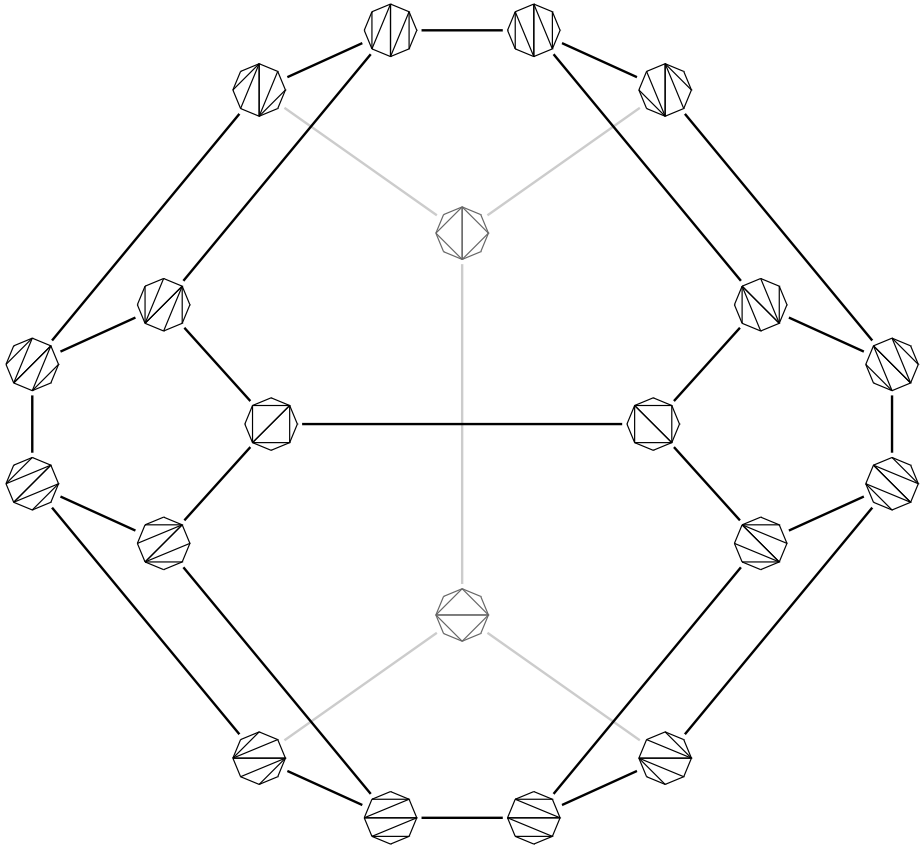}}}
\caption{The $3$-dimensional cyclohedron}
\label{B3assoc}
\end{figure}

The $n$-dimensional \emph{cyclohedron} (also known as the {\em Bott-Taubes 
polytope}~\cite{bott-taubes}) is constructed similarly to the associahedron 
using centrally-symmetric triangulations of a regular $(2n+2)$-gon.
Each edge of the cyclohedron represents either a diagonal flip involving two
diameters of the polygon, or a pair of two centrally-symmetric diagonal flips.
Figures~\ref{B2assoc} and~\ref{B3assoc} show the $2$- and $3$-dimensional 
cyclohedra respectively.
As these figures suggest, the cyclohedron is a convex polytope for any~$n$.  
Explicit polytopal realizations of cyclohedra were constructed by M.~Markl~\cite{markl}
and R.~Simion~\cite{simion-B}.
Each face of a cyclohedron is a product of smaller cyclohedra and
associahedra.

Further details about the combinatorics of cyclohedra, 
and about their appearance in the study of configuration spaces
can be found in~\cite{devadoss-cyclo}. 



The geometry of associahedra and cyclohedra
is related to the geometry of permutohedra, as the following theorem (due to Tonks~\cite{Tonks}) shows.

\begin{theorem}
\label{ABskeleta}
The $1$-skeleton of the $n$-dimensional associahedron (resp.,
cyclohedron) can be obtained from the  
$1$-skeleton of the permutohedron of type~$A_n$ (resp., type~$B_n$) 
by contraction of edges. 
\end{theorem}

\begin{figure}[ht!]
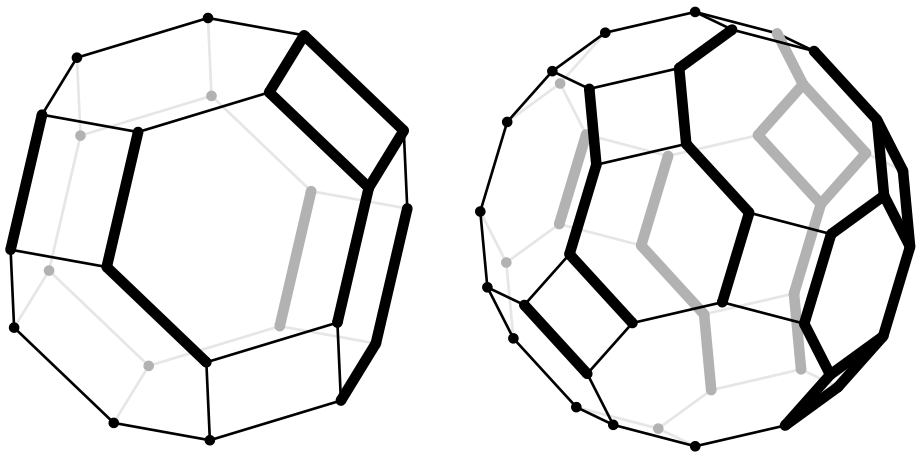

\centerline{
        \epsfbox{A3permequiv.ps}
\qquad \epsfbox{B3permequiv.ps}
}
\caption{Contracting edges of permutohedra of types $A_3$ and $B_3$ yields an 
associahedron and a cyclohedron
}
\label{permequiv-A3B3}
\end{figure}

Theorem~\ref{ABskeleta} is further discussed in 
Section~\ref{sec:lattice cong} in connection with 
Theorem~\ref{cluster refine}. 
For $n=3$, the theorem is illustrated in
Figure~\ref{permequiv-A3B3}. 
(Cf.\ Figures~\ref{A3perm} and~\ref{B3perm}.)

In light of Theorem~\ref{ABskeleta}, 
the cyclohedron can be viewed as a ``type~$B$ counterpart'' of the
associahedron (which is a ``type~$A$'' object). 

\ \bigskip

\pagebreak

\section{Matrix mutations}
\label{sec:matrix-mut}

Having looked closely at the associahedron and the cyclohedron, one is 
naturally led to wonder: 
are these two  
just a pair of isolated constructions, or is there a general framework 
that includes them as special cases? 
Given that the associahedra and the cyclohedra are related to the
root systems of types $A$ and~$B$, respectively, is there a 
classification of polytopes arising within this framework
that is similar to the Cartan-Killing classification of root systems?

As a first step towards answering these questions, 
we will develop the machinery of \emph{matrix mutations}, 
which encode the combinatorics of various models 
similar to the associahedron and the cyclohedron.
We begin our discussion of matrix mutations by continuing the example of the 
associahedron.

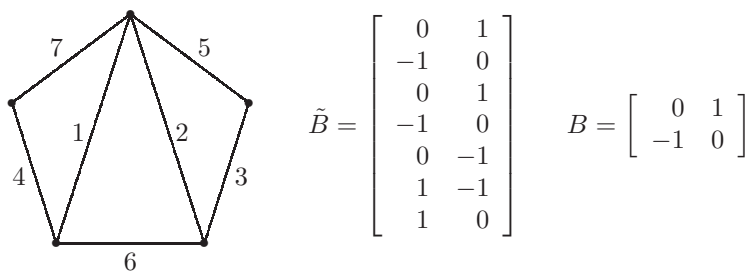
\begin{figure}[ht]
\begin{center}
\setlength{\unitlength}{2.8pt}
\begin{picture}(90,35)(0,-2)
  \put(6,0){\line(1,0){20}}
  \qbezier(6,0)(11,15.5)(16,31)
  \qbezier(26,0)(21,15.5)(16,31)
  \qbezier(6,0)(3,9.5)(0,19)
  \qbezier(26,0)(29,9.5)(32,19)
  \qbezier(0,19)(8,25)(16,31)
  \qbezier(32,19)(24,25)(16,31)

  \put(6,0){\circle*{1}}
  \put(26,0){\circle*{1}}
  \put(0,19){\circle*{1}}
  \put(32,19){\circle*{1}}
  \put(16,31){\circle*{1}}

\put(16,-2.5){\makebox(0,0){$6$}}
\put(1,9){\makebox(0,0){$4$}}
\put(6,26.5){\makebox(0,0){$7$}}
\put(26,26.5){\makebox(0,0){$5$}}
\put(31,9){\makebox(0,0){$3$}}

\put(9,15){\makebox(0,0){$1$}}
\put(23,15){\makebox(0,0){$2$}}

\put(40,15){{$
\tilde B=
\left[\begin{array}{rr}
0 & 1\\
-1& 0\\
0 & 1\\
-1& 0\\
0 &-1\\
1 &-1\\
1 & 0
\end{array}\right]
$}}

\put(75,15){{$
B=
\left[\begin{array}{rr}
0 & 1\\
-1& 0
\end{array}\right]
$}}

\end{picture}
\end{center}
\caption{Matrices $B$ and $\tilde B$ for a triangulation}
\label{B-tilde-example}
\end{figure}

Fix a triangulation $T$ of the $(n\!+\!3)$-gon.
Label the $n$ diagonals of $T$ arbitrarily by the numbers
$1,\dots,n$,
and label the $n+3$ sides of~$T$ 
by the numbers $n+1,\dots,2n+3$. 
The combinatorics of $T$ can be encoded by the (signed) \emph{edge-adjacency 
matrix}~\hbox{$\tilde B=(b_{ij})$}.  
This is the $(2n+3)\times n$ matrix 
whose entries 
are given by
\[
b_{ij}=
\begin{cases}
1 & \text{if $i$ and $j$ label two sides in some triangle of~$T$
  so that $j$ follows~$i$}\\
  & \text{\ \ in the clockwise traversal of the triangle's boundary;}\\
-1 & \text{if the same holds, with the counter-clockwise order;}\\
0 & \text{otherwise.} 
\end{cases}
\]
Note that the first index $i$ is a label for a side or a 
diagonal of the $(n\!+\!3)$-gon, while the second index~$j$ must label a
diagonal. 
The \emph{principal part} of~$\tilde B$ is an $n\times n$ submatrix 
$B=(b_{ij})_{i,j\in [n]}$ that encodes the signed adjacencies between
the diagonals of~$T$. 
An example is shown in Figure~\ref{B-tilde-example}. 

In the language of matrices~$\tilde B$ and~$B$, 
diagonal flips can be described as
certain transformations called matrix mutations. 

\begin{definition} 
\label{def:matrix mutation} 
{\rm 
Let 
$B = (b_{ij})
$ 
and $B' = (b'_{ij})
$ 
be 
integer matrices. 
We say that $B'$ is obtained from $B$ by a \emph{matrix mutation} in 
direction~$k$,  
and write $B' = \mu_k (B)$, if 
\begin{equation} 
\label{eq:mutation} 
b'_{ij} = 
\begin{cases} 
-b_{ij} & \text{if $k\in\{i,j\}$;} \\[.05in] 
b_{ij} + |b_{ik}| b_{kj} & \text{if  $k\notin\{i,j\}$ and $b_{ik}b_{kj}>0$;}\\
\ \ b_{ij} & \text{otherwise.}
\end{cases} 
\end{equation} 
}
\end{definition} 

It is easy to check that a matrix mutation is an involution,
i.e., $\mu_k(\mu_k(B))=B$. 

\begin{lemma}
\label{lem:mut-tri}
Assume that $\tilde B$ and $\tilde B'$ (resp., $B$ and~$B'$) 
are the edge-adjacency matrices 
(resp., their principal parts) 
for two triangulations $T$ and $T'$ obtained from each other by flipping the
diagonal numbered~$k$; 
the rest of the labels are the same in $T$ and~$T'$. 
Then $\tilde B'=\mu_k(\tilde B)$ (resp., $B'=\mu_k(B)$). 
\end{lemma}

\begin{figure}[ht]
\centerline{
        \epsfbox{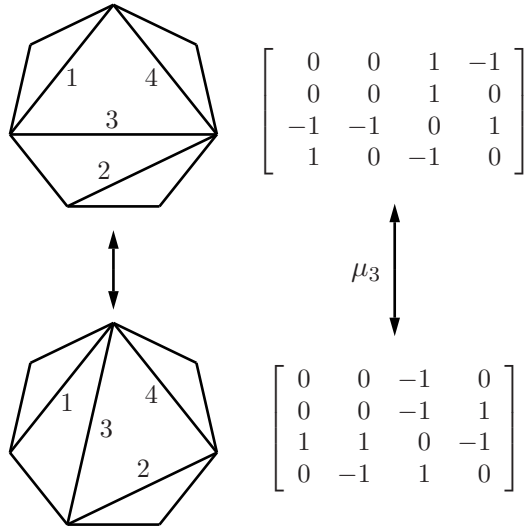}
        \begin{picture}(0,0)(100,-100)
        \put(-67,54){$3$}
        \put(-70,35){$2$}
        \put(-82,70){$1$}
        \put(-52,70){$4$}
        \put(-69,-63){$3$}
        \put(-55,-78){$2$}
        \put(-84,-53){$1$}
        \put(-52,-50){$4$}
        \put(-10,58){$\left[\begin{array}{rrrr}0&0&1&-1\\0&0&1&0\\
                                                -1&-1&0&1\\1&0&-1&0
                        \end{array}\right]$}
        \put(-6,-62){$\left[\begin{array}{rrrr}0&0&-1&0\\
                                                0&0&-1&1\\1&1&0&-1\\
                                                0&-1&1&0
                        \end{array}\right]$}
        \put(26,-2){\large$\mu_3$}
        \end{picture}
}
\ \vspace{-.3in}
\caption{A diagonal flip and the corresponding matrix mutation 
} 
\label{matrices}
\end{figure}

\begin{figure}[p!]
\centerline{
         \begin{picture}(0,0)(-134,-545)
        \put(-20,-26){$\left[\begin{array}{rr}0&1\\-1&0\end{array}\right]$}
        \put(-20,-126){$\left[\begin{array}{rr}0&-1\\1&0\end{array}\right]$}
        \put(-20,-226){$\left[\begin{array}{rr}0&1\\-1&0\end{array}\right]$}
        \put(-20,-326){$\left[\begin{array}{rr}0&-1\\1&0\end{array}\right]$}
        \put(-20,-426){$\left[\begin{array}{rr}0&1\\-1&0\end{array}\right]$}
        \put(-20,-526){$\left[\begin{array}{rr}0&-1\\1&0\end{array}\right]$}
        \put(-113,-26){$\scriptstyle{1}$}
        \put(-104,-26){$\scriptstyle{2}$}
        \put(-113,-126){$\scriptstyle{1}$}
        \put(-104,-126){$\scriptstyle{2}$}
        \put(-113,-226){$\scriptstyle{1}$}
        \put(-104,-213){$\scriptstyle{2}$}
        \put(-101,-335){$\scriptstyle{1}$}
        \put(-104,-313){$\scriptstyle{2}$}
        \put(-101,-435){$\scriptstyle{1}$}
        \put(-112,-426){$\scriptstyle{2}$}
        \put(-104,-526){$\scriptstyle{1}$}
        \put(-112,-526){$\scriptstyle{2}$}
         \end{picture}
       \epsfbox{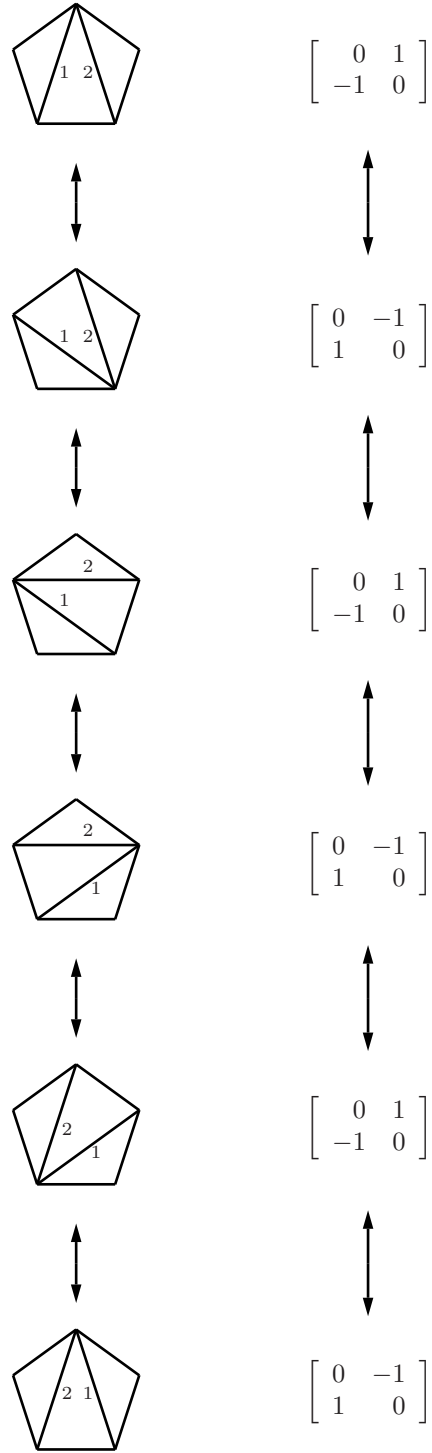}
}
\caption{Diagonal flips in a pentagon, and the corresponding
matrix mutations}
\label{flipmut}
\end{figure}

Lemma~\ref{lem:mut-tri} is illustrated in Figures~\ref{matrices}
and~\ref{flipmut}.
Note that the numbering of diagonals used in defining the
matrices~$\tilde B$ and~$B$ can change as we move along the exchange
graph. 
For instance, the sequence of $5$~flips shown in Figure~\ref{flipmut}
results in switching the labels of the two diagonals.

One can similarly define edge-adjacency matrices 
for centrally symmetric triangulations
(those matrices will have entries~$0$, $\pm1$, and~$\pm 2$), 
and verify that cyclohedral
flips translate precisely into matrix  mutations.

\section{Exchange relations}
\label{sec:exchange-rel}

We next introduce a set of algebraic (more
precisely, \emph{birational}) transformations that will go hand in hand 
with the matrix mutations. 
We start by explaining this construction in the case of an associahedron. 

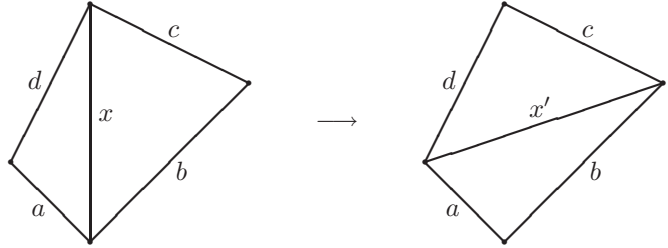
\begin{figure}[ht] 
\begin{center} 
\setlength{\unitlength}{1.5pt} 
\begin{picture}(60,60)(0,0) 
\thicklines 
  \put(0,20){\line(1,2){20}} 
  \put(0,20){\line(1,-1){20}} 
  \put(20,0){\line(0,1){60}} 
  \put(20,0){\line(1,1){40}} 
  \put(20,60){\line(2,-1){40}} 

  \put(20,0){\circle*{1}} 
  \put(20,60){\circle*{1}} 
  \put(0,20){\circle*{1}} 
  \put(60,40){\circle*{1}} 

\put(24,32){\makebox(0,0){$x$}} 

\put(7,8){\makebox(0,0){$a$}} 
\put(43,18){\makebox(0,0){$b$}} 
\put(41,53){\makebox(0,0){$c$}} 
\put(6,41){\makebox(0,0){$d$}} 
\end{picture} 
\begin{picture}(40,66)(0,0) 
\put(20,30){\makebox(0,0){$\longrightarrow$}} 
\end{picture} 
\begin{picture}(60,66)(0,0) 
\thicklines 
  \put(0,20){\line(1,2){20}} 
  \put(0,20){\line(1,-1){20}} 
  \put(0,20){\line(3,1){60}} 
  \put(20,0){\line(1,1){40}} 
  \put(20,60){\line(2,-1){40}} 

  \put(20,0){\circle*{1}} 
  \put(20,60){\circle*{1}} 
  \put(0,20){\circle*{1}} 
  \put(60,40){\circle*{1}} 

\put(29,34){\makebox(0,0){$x'$}}

\put(7,8){\makebox(0,0){$a$}} 
\put(43,18){\makebox(0,0){$b$}} 
\put(41,53){\makebox(0,0){$c$}} 
\put(6,41){\makebox(0,0){$d$}}

\end{picture} 

\end{center} 
\caption{A diagonal flip} 
\label{fig:type-a-exch} 
\end{figure}

Let us fix an arbitrary \emph{initial triangulation}~$T_\circ$ of a convex 
$(n+3)$-gon, 
and introduce a variable for each diagonal of this triangulation,
and also for each side of the $(n+3)$-gon. 
This gives $2n+3$ variables altogether. 
We are now going to associate a rational function in these $2n+3$
variables to \emph{every} diagonal of the $(n+3)$-gon. 
This will be done in a recursive fashion. 
Whenever we perform a diagonal flip as the one shown 
in Figure~\ref{fig:type-a-exch}, all but one rational functions associated to the
current triangulation remain unchanged: 
the rational function $x$ associated with the diagonal being removed gets
replaced by the rational function $x'$ associated with the ``new''
diagonal, where $x'$~is determined from the \emph{exchange relation} 
\begin{equation}
\label{eq:ptolemy}
x\,x'=a\,c+b\,d\,. 
\end{equation}

\begin{lemma}
\label{no mono}
The rational function $x_\gamma$ associated to each diagonal~$\gamma$ 
does not depend  on the particular choice of a sequence of flips that connects
the initial triangulation with another one containing~$\gamma$. 
\end{lemma}

Lemma~\ref{no mono} can be rephrased as saying that there 
are no ``monodromies'' associated with sequences of 
flips that begin and end at the same triangulation.

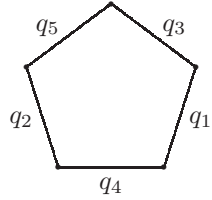
\begin{figure}[ht]
\begin{center}
\setlength{\unitlength}{2pt}
\begin{picture}(35,35)(0,-2)
  \put(6,0){\line(1,0){20}}
  \qbezier(6,0)(3,9.5)(0,19)
  \qbezier(26,0)(29,9.5)(32,19)
  \qbezier(0,19)(8,25)(16,31)
  \qbezier(32,19)(24,25)(16,31)

  \put(6,0){\circle*{1}}
  \put(26,0){\circle*{1}}
  \put(0,19){\circle*{1}}
  \put(32,19){\circle*{1}}
  \put(16,31){\circle*{1}}

\put(16,-3.5){\makebox(0,0){$q_4$}}
\put(-1,9){\makebox(0,0){$q_2$}}
\put(4,26.5){\makebox(0,0){$q_5$}}
\put(28,26.5){\makebox(0,0){$q_3$}}
\put(33,9){\makebox(0,0){$q_1$}}

\end{picture}
\end{center}
\caption{Labeling the sides of a pentagon}
\label{fig:pentagon-A2}
\end{figure}

To illustrate Lemma~\ref{no mono}, consider the triangulations of a
pentagon (i.e.,~$n=2$). 
We label the sides of the pentagon by the variables
$q_1,q_2,q_3,q_4,q_5$, as shown in Figure~\ref{fig:pentagon-A2}. 

\begin{figure}[ht]
\centerline{
        \epsfbox{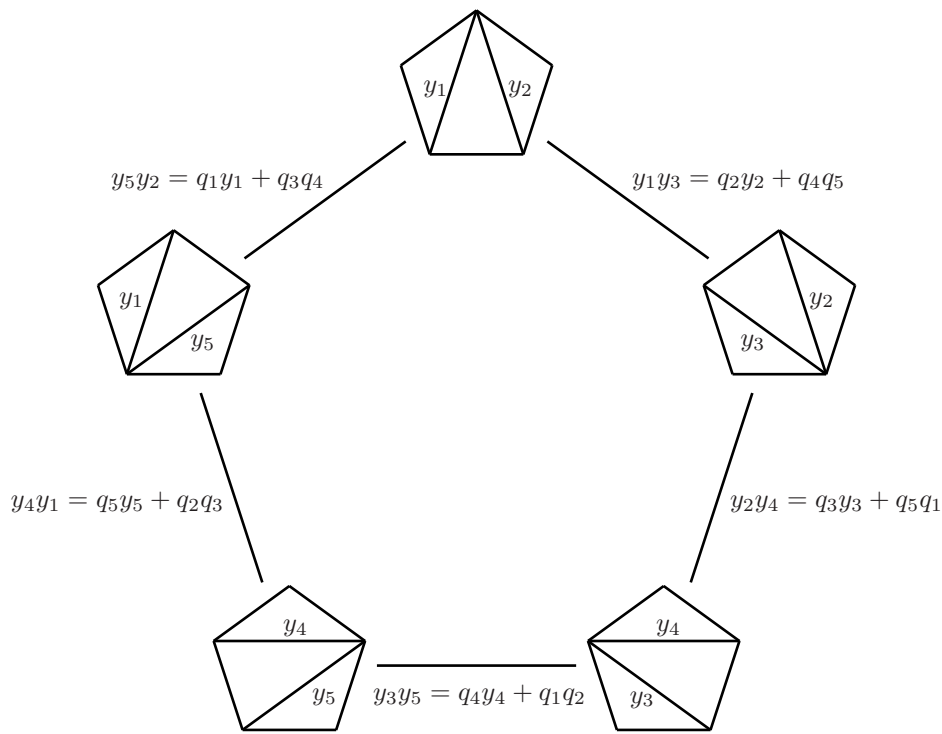}
        \begin{picture}(0,0)(144,-121)
                \put(-24,120){$y_1$}
                \put(8,120){$y_2$}
                \put(122,40){$y_2$}
                \put(96,24){$y_3$}
                \put(64,-84){$y_4$}
                \put(54,-111){$y_3$}
                \put(-77,-84){$y_4$}
                \put(-66,-110){$y_5$}
                \put(-139,40){$y_1$}
                \put(-112,24){$y_5$}
                \put(55,85){$y_1 y_3=q_2 y_2 + q_4 q_5$}
                \put(-142,85){$y_5 y_2=q_1 y_1 + q_3 q_4$}
                \put(-180,-37){$y_4 y_1=q_5 y_5 + q_2 q_3$}
                \put(-43,-110){$y_3 y_5=q_4 y_4 + q_1 q_2$}
                \put(92,-37){$y_2 y_4=q_3 y_3 + q_5 q_1$}
        \end{picture}
}
\caption{Exchange relations for the flips in a pentagon}
\label{A2assoc}
\end{figure}

We then label the diagonals incident to the top vertex by the
variables $y_1$ and~$y_2$.
Thus, our initial triangulation~$T_\circ$ appears at the top of
Figure~\ref{A2assoc}. 
The rational functions $y_3,y_4,y_5$ associated with the remaining
three diagonals are then computed from the exchange relations 
associated with the flips shown in 
Figure~\ref{A2assoc}. 

Starting at the top of Figure~\ref{A2assoc}
and moving clockwise, we recursively express $y_3,y_4,\dots$ in terms
of~$y_1,y_2$ and $q_1,\dots,q_5$:
\begin{align*}
y_3&=\frac{q_2 y_2 + q_4 q_5}{y_1}\,,\\
y_4&=\frac{q_3 y_3 + q_5 q_1}{y_2}
    =\frac{q_3q_2y_2+q_3q_4q_5+q_5q_1y_1}{y_1y_2}\,,\\
y_5&=\frac{q_4 y_4 + q_1 q_2}{y_3}
    =\cdots=\frac{q_3q_4+q_1y_1}{y_2} \quad \text{(check!)}, 
\end{align*}
and, finally, 
\begin{align*}
y_1&=\frac{q_5 y_5 + q_2 q_3}{y_4}=\cdots=y_1\,,\\
y_2&=\frac{q_1 y_1 + q_3 q_4}{y_5}=\cdots=y_2\,,
\end{align*}
recovering the original values and completing the cycle. 

We note that under the specialization $q_1=\cdots=q_5=1$,
the phenomenon we just observed is nothing else but the
$5$-periodicity of the pentagon recurrence,
which we have thus generalized.

Lemma~\ref{no mono} is a special case of a much more  general result from the
theory of \emph{cluster algebras}.
It can also be proved directly in at least two different ways
briefly sketched below; 
these proofs point at connections of this subject to other areas of
mathematics. 

\subsection*{Ptolemy's Theorem and hyperbolic geometry}

The classical Ptolemy's Theorem asserts that in an inscribed
quadrilateral, 
the sum of the products of the two pairs of
opposite sides equals the product of the two diagonals.
This relation looks exactly like the exchange
relation~\eqref{eq:ptolemy}. 
It suggests that one can prove Lemma~\ref{no mono} simply by
interpreting the rational function associated with each side or
diagonal as the Euclidean length of the corresponding segment. 
There is however a problem with this type of argument: the space of inscribed 
$(n+3)$-gons (up to congruence) is $(n+3)$-dimensional, 
whereas we need $2n+3$ independent variables in our setup. 

The problem can be resolved by passing from Euclidean to hyperbolic
geometry, where an analogue of Ptolemy's Theorem holds,
and where one can ``cook up'' the required additional degrees of
freedom. For much more on this topic, see~\cite{fock-goncharov1,gsv}.


\subsection*{Pl\"ucker coordinates on the Grassmannian
$\operatorname{Gr}(2,n+3)$} 

Take a $2\times (n\!+\!3)$ matrix $z=(z_{ij})$.
For any $k,l\in [n+3]$, $k<l$, let us denote by 
\[
P_{kl}=\det\begin{pmatrix}z_{1k}&z_{1l}\\z_{2k}&z_{2l}\end{pmatrix}
\] 
the $2\times 2$ minor of~$z$ that occupies columns $k$ and~$l$. 
These minors are the homogeneous \emph{Pl\"ucker coordinates} of the row span
of~$z$ as an element of the \emph{Grassmannian}
$\operatorname{Gr}(2,n+3)$ of all $2$-dimensional subspaces of an
$(n\!+\!3)$-space. See, e.g., \cite{Fu-Ha}. 

It is easy to check (the special case of) the
\emph{Grassmann-Pl\"ucker relations}:
\[
P_{ik}P_{jl}=P_{ij}P_{kl}+P_{il}P_{jk}\,.
\]
Once again, one recognizes the exchange
relation~\eqref{eq:ptolemy}. 
It is straightforward to construct, for a particular special
choice of initial triangulation~$T_\circ$, a matrix~$z$ for which the
values of the minors~$P_{kl}$ corresponding to the 
sides and diagonals of~$T_\circ$ are equal to the variables associated
with these segments. 
It then follows by induction that the rational function associated to
every diagonal is equal to the corresponding minor~$P_{kl}$,
implying Lemma~\ref{no mono}. 

\lecture{Cluster Algebras}
\label{lec:cluster}

Our next task is to create a general axiomatic theory of mutations
(``flips'') and exchanges, using the above examples as prototypes. 
This will lead us to the basic notions and results of the theory of
cluster algebras. 
Cluster algebras were introduced in \cite{ca1} as a
combinatorial/algebraic framework for the study of dual canonical
bases and related total positivity phenomena.
They since found applications in higher Teichm\"uller theory and
representation theory of quivers, among other fields.
All these motivations and applications will remain behind the scenes
in these lectures. 

Most of this lecture is based on \cite{ga, ca1, ca2}.
Sections~\ref{sec:polytopal} and~\ref{sec:3by3} are based on
\cite{gaPoly} and~\cite{ca3, tptp},
respectively. 

\section{Seeds and clusters} 


Consider a diagonal flip that transforms a triangulation~$T$ of a
convex $(n+3)$-gon into another triangulation~$T'$, as shown in
Figure~\ref{fig:type-a-exch}. 
The corresponding exchange relation~\eqref{eq:ptolemy} can be written
entirely in terms of the edge-adjacency matrix~$\tilde B$.
To be more precise, let us assume, as before, that the diagonals of~$T$ have
been labeled in some way by the numbers $1,\dots,n$, whereas the sides
of the $(n+3)$-gon have been assigned the labels going from $n+1$
through~$m=2n+3$. 
The labeling for~$T'$ is the same except for the one diagonal
(say, labeled~$k$) that is getting exchanged between $T$ and~$T'$.

This labeling of sides and diagonals of $T$ allows us to (temporarily) denote
the associated rational functions by $x_1,\dots,x_m$. 
For~$T'$, we get the same rational functions except that $x_k$
is replaced by~$x_k'$. 
Then the exchange relation under consideration takes the form 
\begin{equation}
\label{eq:exch-rel-gen}
x_k\, x'_k = \prod\limits_{
\begin{array}{c}\scriptstyle b_{ik}>0\\[-.05in] \scriptstyle 1\leq i\leq
  m\end{array}} 
x_i^{b_{ik}}+
\prod\limits_{\begin{array}{c}\scriptstyle b_{ik}<0\\[-.05in] \scriptstyle 1\leq i\leq
  m\end{array}} x_i^{-b_{ik}}\,. 
\end{equation}
In other words, the right-hand side of \eqref{eq:exch-rel-gen} is the
sum of two monomials whose exponents are the absolute values of the
entries in the $k$th column of~$\tilde B$,
while the sign of an entry determines which monomial the corresponding
term contributes~to. 

\begin{example}
{\rm
Let $T$ be the triangulation of a pentagon in
Figure~\ref{B-tilde-example}, with its edges labeled $1,\dots,7$ as
shown.
The exchange relations corresponding to flipping the diagonals $1$
and~$2$ are, respectively:
\begin{align*}
x_1 x_1' &= x_6 x_7 +x_2 x_4 \,,\\
x_2 x_2' &= x_1 x_3 +x_5 x_6 \,,
\end{align*}
in agreement with~\eqref{eq:exch-rel-gen}. 
}
\end{example}

To summarize, both the combinatorics of flips and the algebra of
exchange relations can be encoded entirely in terms of the
matrices~$\tilde B$ using, first, the machinery of matrix mutations
and, second, the ``birational dynamics'' governed
by~\eqref{eq:exch-rel-gen}. 
We shall now use this observation to lay out the axioms of a cluster
algebra. 
The formulation of these axioms will require some technical
preparation, which hopefully will make sense to the reader in light of
the examples discussed above.

A cluster algebra $\AA$ is a commutative ring 
contained in an  \emph{ambient field}~$\mathcal{F}$
isomorphic to the field of rational
functions in $m$ variables over~$\mathbb{Q}$.
(Think of the rational functions in the variables associated with the 
sides and diagonals of a fixed initial triangulation.)  

$\AA$ is generated
inside~$\mathcal{F}$ by a (possibly infinite) set of generators.
These generators are obtained from
an initial \emph{seed} via an iterative process of \emph{seed
  mutations} which follows a set of canonical rules.

A \emph{seed} in $\mathcal{F}$ is a pair $(\tilde \xx, \tilde B)$, where
\begin{itemize}
\item
$\tilde \xx\! =\! \{x_1, \dots, x_m\}$
is a set\footnote{
A diligent reader might object that we call $\tilde \xx$ a set rather
than a sequence. 
This is because we regard any two seeds obtained from each other by
simultaneous relabeling of the  elements~$x_i$ and the matrix entries
$b_{ij}$ as identical.
That is, for any permutation $w\in\mathcal{S}_m$ such that $w(i)=i$ for
$i>n$, 
we make no distinction between the seeds $(\tilde \xx, \tilde B)$
and $(w(\tilde\xx), w(\tilde B))$, where $w(\tilde\xx)=(x_{w(1)},\dots,x_{w(m)})$ and
$w(\tilde B)=(b_{w(i),w(j)})$.
}
of $m$ algebraically independent generators of~$\mathcal{F}$,
which is split into a disjoint union of an $n$-element \emph{cluster} 
$\xx = \{x_1, \dots, x_n\}$ 
and an $(m-n)$-element set of \emph{frozen variables} 
$\cc \! =\! \{x_{n+1}, \dots, x_m\}$; 
\item
$\tilde B\!=\! (b_{ij})$
is an $m\times n$ integer matrix of rank~$n$ 
whose \emph{principal part}~$B\!=\! (b_{ij})_{i,j\in [n]}$
is \emph{skew-symmetrizable}, 
i.e., 
there exists a diagonal matrix $D$ with positive diagonal entries such
that $DBD^{-1}$ is skew-symmetric. 
\end{itemize}
(Equivalently, there exist positive integers $d_1,\dots,d_n$ such that
$d_i b_{ij} = - d_j b_{ji}$ for all $i$ and~$j$.) 
The matrix $B$ is called the \emph{exchange matrix} of a seed.

A \emph{seed mutation} $\mu_k$ in direction~$k\in\{1, \dots, n\}$ transforms
a seed $(\tilde \xx,\tilde B)$ into another seed $(\tilde \xx',\tilde
B')$ defined as follows: 
\begin{itemize}
\item
$\tilde \xx' = \tilde \xx - \{x_k\} \cup \{x'_k\}$, where 
$x_k'$ is found from the exchange relation~\eqref{eq:exch-rel-gen};
\item
$\tilde B'=\mu_k(\tilde B)$, i.e., $\tilde B$ undergoes a matrix
mutation (hence so does~$B$). 
\end{itemize}
The following lemma justifies the definition of a seed mutation
by showing that $(\tilde \xx',\tilde B')$ is indeed a seed.  

\begin{lemma}
\label{lem:mut-inv}
Matrix mutations preserve the rank of a matrix.
If $B$ is skew-symmetrizable, then so is~$\mu_k(B)$,
with the same skew-symmetrizing matrix~$D$. 
\end{lemma}

Note that seed mutations do not change the frozen variables
$\cc\! =\! \{x_{n+1}, \dots, x_m\}$. 

\begin{example}
\label{example:typeA-1}
{\rm
Let $\xx$ and $\cc$ be the sets of variables associated with
the diagonals and sides, respectively, of some triangulation of
a convex~$(n+3)$-gon. 
(Thus $m=2n+3$.) 
Let $\tilde B$ be the sign-adjacency matrix of the triangulation. 
The mutations of seeds $(\tilde\xx,\tilde B)$ of this kind 
correspond to combining the exchange relations \eqref{eq:ptolemy}
with the matrix mutations associated with diagonal flips. 
}
\end{example}

\begin{figure}[ht]
\begin{center}
\setlength{\unitlength}{1.2pt}
\begin{picture}(200,110)(0,0)

\put(100,50){\makebox(0,0){\shortstack{initial\\ seed}}}
\put(100,50){\oval(36,22)}

\put(100,8){\makebox(0,0){seed}}
\put(100,8){\oval(32,12)}

\put(100,39){\line(0,-1){25}}

\put(115.2,58.2){\line(1,1){20.5}}
\put(84.8,58.2){\line(-1,1){20.5}}

\put(150,83){\makebox(0,0){seed}}
\put(150,83){\oval(32,12)}

\put(50,83){\makebox(0,0){seed}}
\put(50,83){\oval(32,12)}

\put(114.3,3.7){\line(1,-1){10}}
\put(85.7,3.7){\line(-1,-1){10}}

\put(150,89){\line(0,1){15}}
\put(50,89){\line(0,1){15}}

\put(164.3,78.7){\line(1,-1){10}}
\put(35.7,78.7){\line(-1,-1){10}}

\end{picture}
\end{center}
\caption{Seed mutations and the exchange graph} 
\label{fig:seed-mut}
\end{figure}
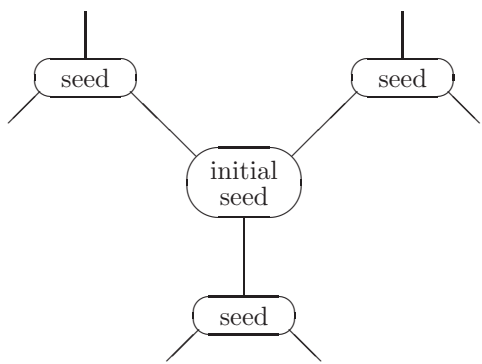

Seed mutations generate the \emph{mutation equivalence} relation on seeds:
$(\tilde \xx,\tilde B) \sim (\tilde \xx',\tilde B')$.
Let $\mathcal{S}$ be an equivalence class for this relation.
Thus, $\mathcal{S}$ is obtained by repeated 
mutations of an arbitrary initial seed in all possible directions.
This creates an \emph{exchange graph}. 
See Figure~\ref{fig:seed-mut}. 

\begin{figure}[ht]
\begin{center}
\setlength{\unitlength}{1.5pt}
\begin{picture}(200,110)(0,-5)

\put(100,50){\makebox(0,0){$x_1,x_2,x_3$}}
\put(100,50){\oval(36,12)}

\put(100,8){\makebox(0,0){$x_1',x_2,x_3$}}
\put(100,8){\oval(32,12)}

\put(100,44){\line(0,-1){30}}

\put(114.5,55.5){\line(1,1){22.3}}
\put(85.5,55.5){\line(-1,1){22.3}}

\put(150,83){\makebox(0,0){$x_1,x_2,x_3'$}}
\put(150,83){\oval(32,12)}

\put(50,83){\makebox(0,0){$x_1,x_2',x_3$}}
\put(50,83){\oval(32,12)}

\put(114.3,3.7){\line(1,-1){10}}
\put(85.7,3.7){\line(-1,-1){10}}

\put(150,89){\line(0,1){15}}
\put(50,89){\line(0,1){15}}

\put(164.3,78.7){\line(1,-1){10}}
\put(35.7,78.7){\line(-1,-1){10}}

\put(135,20){\shortstack[l]{}}

\end{picture}
\end{center}
\caption{Cluster variables} 
\label{fig:clust-var}
\end{figure}
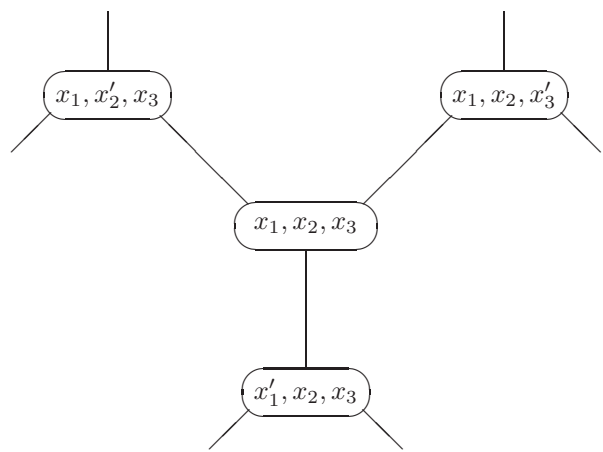

Let $\mathcal{X}=\mathcal{X}(\mathcal{S})$ be the union of all clusters
for all the seeds in~$\mathcal{S}$.
The elements of $\mathcal{X}$ are
called \emph{cluster variables}.
See Figure~\ref{fig:clust-var}.

The \emph{cluster algebra}\footnote{
Strictly speaking, this is a definition of a
\emph{skew-symmetrizable
cluster algebra of geometric type}. 
} 
$\AA=\AA(\mathcal{S})$
associated with $\mathcal{S}$ is generated inside~$\mathcal{F}$ 
by the cluster variables in  $\mathcal{X}$ together with the frozen
variables $x_{n+1},\dots,x_m$ and their inverses. 
(A variation of this definition includes 
cluster and frozen variables, 
but none of their inverses, in the generating set.) 
The integer $n$ is called the \emph{rank} of~$\AA$. 

\begin{theorem}
[The Laurent phenomenon~\cite{ca1}] 
\label{th:laurent-phenom}
Any cluster variable 
is expressed in terms of the variables $x_1,\dots,x_m$ of any
given seed as a Laurent polynomial with integer coefficients.
\end{theorem}

\begin{conjecture}[Nonnegativity conjecture~\cite{ca1}] 
\label{conj:nonneg}
Every coefficient in these Laurent polynomials 
is nonnegative.
\end{conjecture}

Conjecture~\ref{conj:nonneg} has been proved in a number of
special cases, including our main motivating example of an
associahedron, to which we return in Example~\ref {example:typeA-2}. 

\begin{example}
\label{example:typeA-2}
{\rm
In the case of Example~\ref{example:typeA-1},
the exchange graph on seeds is precisely the exchange graph on
triangulations illustrated in Figures~\ref{A2assoc_basic}
and~\ref{A3assoc}. 
The cluster algebra in this example is generated inside the ring of
rational functions in $2n+3$ independent variables 
by the rational functions 
associated with \emph{all}
diagonals and sides of the~$(n+3)$-gon. 
(Cf.~Lemma~\ref{no mono}.)
Here we use a variation of the definition of a cluster algebra where the
inverses of frozen 
variables are not included in the set of generators. 

This cluster algebra is canonically isomorphic to the homogeneous
coordinate ring of the Grassmannian $\operatorname{Gr}(2,n+3)$
with respect to its Pl\"ucker embedding. The cluster variables,
together with the frozen variables, form the set of all Pl\"ucker
coordinates on this Grassmannian. 
Theorem~\ref{th:laurent-phenom} and Conjecture~\ref{conj:nonneg} 
(proven in this special case) assert 
that any Pl\"ucker coordinate is written in terms of the Pl\"ucker
coordinates associated with a given triangulation as a Laurent
polynomial with nonnegative integer coefficients.
}\end{example}

\section{Finite type classification}
\label{sec:clust-fintype}

All results in this section were obtained in~\cite{ca2}. 

A cluster algebra is said to be
of \emph{finite
type} if it has finitely many distinct seeds.
Amazingly, the classification of the cluster algebras of finite type
turns out to be completely parallel
to the Cartan-Killing classification of (finite crystallographic) root systems.
Thus there is a cluster algebra of finite type for each Dynkin
diagram, or each Cartan matrix of finite type.
We shall now explain how.  

For a Cartan matrix $A\!=\!(a_{ij})$ of finite type,
we define a skew-symmetrizable  matrix $B(A)\!=\!(b_{ij})$ by
\[
b_{ij}=
\begin{cases}
\ \ 0 & \text{if $i=j$;} \\
\ a_{ij} & \text{if $i\neq j$ and $i\in I_+$;}\\
-a_{ij} & \text{if $i\neq j$ and $i\in I_-$,}
\end{cases}
\]
where $I_+$ and $I_-$ are defined as in
Section~\ref{sec:Coxeter element}. 
To illustrate, in type~$B_4$, we have (cf.\ Example~\ref{example:a4b4c4d4}):
\[
A=\left[\begin{array}{rrrrr}
2&-2&0&0\\
-1&2&-1&0\\
0&-1&2&-1\\
0&0&-1&2\\
\end{array}\right]\,,\qquad
B(A)=\left[\begin{array}{rrrrr}
0&-2&0&0\\
1&0&1&0\\
0&-1&0&-1\\
0&0&1&0\\
\end{array}\right]\,,
\]
under the convention $I_+=\{1,3\}$, $I_-=\{2,4\}$. 

\begin{theorem}[Finite type classification]
\label{th:fin-type-class}
A cluster algebra $\AA$ is of finite type
if and only if the exchange matrix at some
seed of~$\AA$ is of the form~$B(A)$, where $A$ is a Cartan matrix
of finite type.

The type of~$A$ (in the Cartan-Killing
nomenclature) is uniquely determined by the cluster
algebra~$\AA$, and is called the ``cluster type'' of~$\AA$.
\end{theorem}


We note that in deciding whether a cluster algebra is of finite type,
the bottom part of the matrix~$\tilde B$ plays no role whatsoever:
everything is determined by its principal part~$B$. 

In the special cases where a cluster algebra has rank~$n=2$,
is of finite type (that is, one of the types $A_2$, $B_2$, and~$G_2$),
and has no frozen variables (that is, $m=2$), 
Theorem~\ref{th:fin-type-class} brings us back to the recurrences of
Section~\ref{abel}. 
Indeed, these recurrences are precisely given by the exchange
relations in those cluster algebras. 
The periodicity of the corresponding sequences is simply a
reformulation of the ``finite type'' property for cluster algebras.

\begin{theorem}[Combinatorial criterion for finite type] 
\label{th:fin-type-characterizations}
A cluster algebra $\AA$ is of finite type if and only if 
the exchange matrix~$B=(b_{ij})$ for any seed of~$\AA$ satisfies the inequalities
$|b_{ij}b_{ji}|\leq 3$
for all $i,j\in\{1,\dots,n\}$.
\end{theorem}

To rephrase, a mutation equivalence class of skew-symmetrizable
$n\times n$ matrices defines a class of cluster algebras of finite
type
if and only if, for each matrix $B\!=\!(b_{ij})$ in this 
equivalence class, the inequality  $|b_{ij}b_{ji}|\leq 3$ holds
for all $i$ and~$j$.

Combining Theorems~\ref{th:fin-type-characterizations}
and~\ref{cartan thm} yields the following completely elementary
statement about integer matrices, no direct proof of which is
known\footnote{\emph{Note added in revision.}
According to A.~Zelevinsky, such a proof has been
recently found in his joint work with M.~Barot and C.~Geiss. 
}. 

\begin{corollary}
\label{cor:fintype-equivalence}
Let $\mathfrak{B}$ be a mutation equivalence class of
skew-symmetrizable integer matrices, with the skew-symmetrizing
matrix~$D$. (Cf.\ Lemma~\ref{lem:mut-inv}.) 
The following are equivalent:
\begin{itemize}
\item
any matrix $B\!=\!(b_{ij})\!\in\!\mathfrak{B}$ satisfies the inequalities
$|b_{ij}b_{ji}|\!\leq\! 3$, for all~$i$~and~$j$; 
\item
there exists a matrix $B\!=\!(b_{ij})\!\in\!\mathfrak{B}$
with the following property. Define $A=(a_{ij})$ by
\[
a_{ij}=\begin{cases}
-|b_{ij}| & \text{if $i\neq j$;} \\
\ \ 2 & \text{if $i=j$.}
\end{cases}
\]
Then $DAD^{-1}$ is positive definite. 
\end{itemize}
\end{corollary}


\medskip

Let $\Phi$ be an irreducible finite root system with Cartan matrix~$A$,
and let $\AA$ be a cluster algebra of the corresponding cluster type.
Theorem~\ref{th:fin-type-class} tells us that the set $\mathcal{X}$ 
of cluster variables is finite. 
A more detailed description of this set is provided by
Theorem~\ref{th:almost-positive} below. 

Let $\alpha_1,\dots,\alpha_n$ be the simple roots of~$\Phi$,
and let $\{x_1,\dots,x_n\}$ be the cluster at a seed
in~$\AA$ with the exchange matrix~$B(A)$.

Let $\Phi_{\geq -1}$ denote the set of roots in~$\Phi$ which are
either negative simple or positive.
Theorem~\ref{th:almost-positive} shows that the cluster
variables in~$\AA$ are naturally labeled by the roots in~$\Phi_{\geq
  -1}\,$.

\begin{theorem}
\label{th:almost-positive}
For any root
$\alpha=c_1\alpha_1+\cdots+c_n\alpha_n \in\Phi_{\geq -1}\,$,
there is a unique cluster variable $x[\alpha]$ given by
\begin{equation}
\label{eq:clvar-labeling}
x[\alpha]=\frac{P_\alpha(x_1,\dots,x_m)}{x_1^{c_1}\cdots x_n^{c_n}}\,,
\end{equation}
where $P_\alpha$ is a polynomial in $x_1,\dots,x_m$ with nonzero
constant term. 
The map $\alpha\mapsto x[\alpha]$ is a bijection between $\Phi_{\geq
  -1}$ and~$\mathcal{X}$. 
\end{theorem}

Note that the right-hand side of~\eqref{eq:clvar-labeling} is a
Laurent polynomial, in agreement with
Theorem~\ref{th:laurent-phenom}.

\section{Cluster complexes and generalized associahedra
}
\label{sec:ga}

This section is based on \cite{ga,ca2}, except for the last statement
in Theorem~\ref{th:cluster-complex-fintype}, which was proved in~\cite{gaPoly}.

\begin{figure}[ht]
\centerline{
\epsfbox{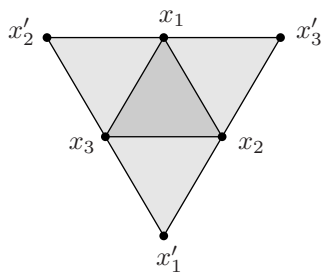}
\begin{picture}(0,0)(91.3,-13.2)
\put(42,81){$x_1$}
\put(72,33){$x_2$}
\put(8,33){$x_3$}
\put(41,-11){$x_1'$}
\put(-15,75){$x_2'$}
\put(94,75){$x_3'$}
\end{picture}
}
\caption{The cluster complex}
\label{fig:clust-complex}
\end{figure}

It can be shown that in a given cluster algebra of finite type,
each seed is uniquely determined by its cluster.
Consequently,
the combinatorics of exchanges 
is encoded by the \emph{cluster complex},  
a simplicial complex
(indeed, a pseudomanifold)
%
%
on the set of all cluster variables
whose maximal simplices (facets) are the clusters.
See Figure~\ref{fig:clust-complex}. 
By Theorem~\ref{th:almost-positive},
the cluster variables---hence the vertices of the cluster
complex---can be naturally labeled by the set~$\Phi_{\geq -1}$ of ``almost positive roots''
in the associated root system~$\Phi$.

This dual graph of the cluster complex is precisely the exchange graph 
of the cluster algebra. 

Theorem~\ref{th:cluster-complex-fintype} below shows 
that the cluster complex is always spherical, and moreover polytopal.

Recall that $Q_\reals$ denotes the $\reals$-span of~$\Phi$. 
The $\ZZ$-span of $\Phi$ is the \emph{root lattice}, denoted by~$Q$. 

\begin{theorem}
\label{th:cluster-complex-fintype}
The $n$ roots that label the
cluster variables in a given cluster
form a $\ZZ$-basis of
the root lattice~$Q$. 
The cones spanned by such $n$-tuples of roots form a complete
simplicial fan 
in the ambient real vector space~$Q_\reals$ (the ``cluster fan'').
This fan 
is the normal fan\footnote{
Let $P\subset V\cong \RR^n$ be an $n$-dimensional simple convex polytope. 
The \emph{support function} 
$F:V^*\to\RR$ 
of~$P$ is given by 
\begin{equation*} 
\label{eq:support-function} 
F(\gamma) = \max_{\mathbf{z} \in P} \langle\mathbf{z},\gamma\rangle . 
\end{equation*} 
The \emph{normal fan} $\mathcal{N}(P)$ is a complete simplicial fan in the 
dual space $V^*$ whose full-dimensional cones 
are the domains of linearity for~$F$.
More precisely, each vertex $\mathbf{z}$ of $P$ gives rise to the cone 
$\{\gamma \in V^*: F(\gamma) = 
\langle\mathbf{z},\gamma\rangle\}$
in~$\mathcal{N}(P)$. 
}
of a simple $n$-dimensional convex polytope 
in the dual space~$Q_\reals^*$. 
\end{theorem}

This polytope is called the \emph{generalized associahedron} of
the corresponding type.

Thus, the cluster complex of a cluster
algebra of finite type 
is canonically isomorphic to the dual simplicial complex of a
generalized associahedron of the corresponding type.
Conversely, the dual graph of the cluster complex is
the \hbox{$1$-skeleton} of the generalized associahedron.

\begin{figure}[ht]
\centerline{
\epsfbox{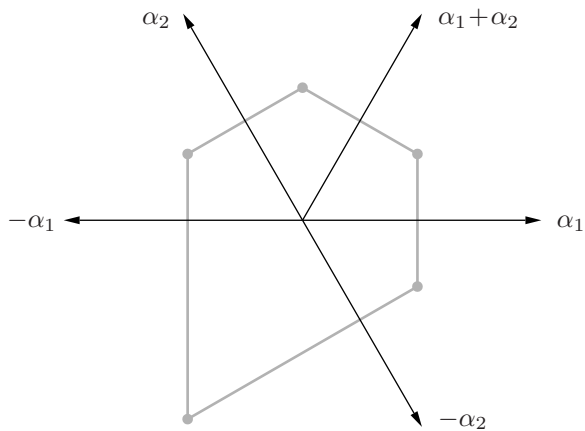}
\begin{picture}(0,0)(112.5,-84.6)
\put(97,-2){$\alpha_1$}
\put(-111,-2){$-\alpha_1$}
\put(52,75){$\alpha_1\!+\!\alpha_2$}
\put(-60,75){$\alpha_2$}
\put(52,-77){$-\alpha_2$}
\end{picture}
}
\caption{Associahedron of type~$A_2$ and its dual fan}
\label{fig:asso-A2}
\end{figure}

In type~$A_n$, 
this construction recovers the $n$-dimensional associahedron
(cf.\ Figure~\ref{fig:asso-A2}). 
The explanation involves an identification of the roots
in~$\Phi_{\geq -1}$ with diagonals of a convex $(n+3)$-gon
that will be discussed later in Example~\ref{example:compat-An}. 
In type~$B_n$, one obtains the $n$-dimensional cyclohedron. 
Thus the $n$-dimensional associahedron (resp., cyclohedron) is
dual to the cluster complex of an arbitrary cluster algebra of
type~$A_n$ (resp.,~$B_n$). 

Theorem~\ref{th:cluster-complex-fintype} leaves the following two
questions unanswered:
\begin{itemize}
\item
Which $n$-subsets of ``almost positive'' roots (``root clusters'')
label the clusters of the cluster algebra of finite type?
(An answer to this question would make explicit the combinatorics of a
generalized associahedron.) 
\item
What are the inequalities defining a generalized associahedron inside~$Q_\reals^*$? 
(We already know they are of the form $\langle
\mathbf{z},\alpha\rangle\leq \mathrm{const}$, for
$\alpha\in\Phi_{\geq-1}$.) 
\end{itemize}

We are now going to answer these questions, one after another.
The answer to the first question is facilitated by the following
property of a cluster complex.

\begin{theorem}
\label{th:clique}
The cluster complex is a clique complex for its \hbox{$1$-skeleton.}
In other words, a subset $S\subset \Phi_{\geq-1}$ is a simplex in the
cluster complex if and only if every \hbox{$2$-element} subset of~$S$
is a \hbox{$1$-simplex} in this complex. 
\end{theorem}

In type~$A_n$, Theorem~\ref{th:clique} reflects the basic property
of the dual complex of an associahedron: 
a collection of diagonals forms a simplex if and only
if any two of them do not cross.  


In order to describe the cluster complex, we therefore need only to
clarify which \emph{pairs} of roots label the edges of the cluster
complex. 
Thus, we need to define the root-theoretic analogue of the
notion of ``non-intersecting diagonals'' that lies at the heart of the
combinatorial construction of an associahedron.  

We will assume from now on that the root system $\Phi$ underlying a
cluster algebra~$\AA$ is irreducible.
(The general case can be obtained by taking direct products.) 
We retain the notation of Lecture~\ref{lec2}.
Thus, $n$ is the rank of~$\Phi$ (and~$\AA$); 
$I$~is the \hbox{$n$-element} indexing set, which is partitioned into
disconnected pieces $I_+$ and~$I_-$;
$W$~is the corresponding reflection group,
generated as a Coxeter group by the generators $s_i$, for $i\in I$;
$w_\circ$ is the element of maximal length in~$W$; 
$A=(a_{ij})$ is the Cartan matrix; 
$h$ is the Coxeter number.




\begin{definition}\rm
\label{def:tau}
Define involutions $\tau_\pm:\Phi_{\geq
  -1}\to\Phi_{\geq -1}$ by
\end{definition}
\begin{equation*} 
\label{eq:tau-pm-on-roots} 
\tau_\varepsilon(\alpha) = 
\begin{cases} 
\displaystyle 
\ \ \alpha & \text{if $\alpha = - \alpha_i\,$, for $i \in I_{- \varepsilon}$;}
\\[.1in] 
\displaystyle\prod_{i \in I_\varepsilon} s_i\,(\alpha) & \text{otherwise.} 
\end{cases} 
\end{equation*} 
For example, in type $A_2$, we get: 
\begin{equation*} 
\label{eq:A2-tau-tropical} 
\hspace{-.1in}
\begin{array}{ccc} 
-\alpha_1 & \stackrel{\textstyle\tau_+}{\longleftrightarrow} 
~\alpha_1~ \stackrel{\textstyle\tau_-}{\longleftrightarrow} 
~\alpha_1\,+\alpha_2~ 
\stackrel{\textstyle\tau_+}{\longleftrightarrow} ~\alpha_2~ 
\stackrel{\textstyle\tau_-}{\longleftrightarrow} 
& -\alpha_2\,\, \\
\circlearrowright & & \circlearrowright \\ \tau_- & & \tau_+ 
\end{array} 
\end{equation*} 
The product $\tau_-\tau_+$ can be viewed as a 
deformation of the Coxeter element. 
Hence, what is the counterpart of the Coxeter number? 

\begin{theorem}
\label{th:dihedral}
The order of $\tau_-  \tau_+$ is
$(h+2)/2$ if $w_\circ = -1$, and is $h+2$ otherwise.
Every $\langle \tau_-,\tau_+ \rangle$-orbit in  $\Phi_{\geq - 1}$ has a
nonempty intersection with $- \Pi$.
These intersections are precisely the $\langle - w_\circ \rangle$-orbits
in~$(- \Pi)$.
\end{theorem}


\begin{theorem}
\label{th:compat}
There is a unique binary relation (called ``compatibility'') on
$\Phi_{\geq -1}$ that has the following two properties: 
\begin{itemize}
\item
$\langle \tau_-,\tau_+ \rangle$-invariance: 
$\alpha$ and $\beta$ are compatible if and only if $\tau_\varepsilon \alpha$
and $\tau_\varepsilon\beta$ are, for $\varepsilon\in\{+,-\}$;
\item
a negative simple root $-\alpha_i$ is compatible with a root~$\beta$ 
if and only if the simple root expansion of~$\beta$ does not
involve~$\alpha_i$. 
\end{itemize}
This compatibility relation is symmetric.
%
The clique complex for the compatibility relation is canonically
isomorphic to the cluster complex. 
\end{theorem}

In other words (cf.\ Theorem~\ref{th:clique}), a subset of roots in
$\Phi_{\geq -1}$ forms a simplex in the cluster complex
if and only if every pair of roots in this subset is compatible.




\begin{example}
\label{example:compat-An}
{\rm
In type~$A_n$, the compatibility relation can be described in concrete
combinatorial terms using a particular identification of the roots in 
$\Phi_{\geq -1}$ with the diagonals of a 
regular $(n+3)$-gon. Under this identification, the 
roots in $- \Pi$ correspond to the diagonals on the ``snake'' 
shown in Figure~\ref{fig:octagon-snake}. 
Each positive root $\alpha_i + \alpha_{i+1} + \cdots + \alpha_j$
corresponds to the 
unique diagonal that crosses precisely the diagonals $- \alpha_i, 
- \alpha_{i+1}, \ldots,- \alpha_j$ from the snake (see 
Figure~\ref{fig:a2}). 
It is easy to check that the transformations $\tau_+$ and~$\tau_-$
act on the set of diagonals as if they were the reflections generating the
dihedral group of symmetries of the~$(n+3)$-gon. 
It then follows that two roots are compatible if and only if the
corresponding diagonals do not cross each other
(at an interior point).
}
\end{example}

\begin{figure}[ht] 
\begin{center} 
\setlength{\unitlength}{2.8pt} 
\begin{picture}(60,60)(0,0) 
\thicklines 
  \multiput(0,20)(60,0){2}{\line(0,1){20}} 
  \multiput(20,0)(0,60){2}{\line(1,0){20}} 
  \multiput(0,40)(40,-40){2}{\line(1,1){20}} 
  \multiput(20,0)(40,40){2}{\line(-1,1){20}} 

  \multiput(20,0)(20,0){2}{\circle*{1}} 
  \multiput(20,60)(20,0){2}{\circle*{1}} 
  \multiput(0,20)(0,20){2}{\circle*{1}} 
  \multiput(60,20)(0,20){2}{\circle*{1}} 

\thinlines \put(0,20){\line(1,0){60}} \put(0,40){\line(1,0){60}} 
\put(0,20){\line(2,-1){40}} \put(0,40){\line(3,-1){60}} 
\put(20,60){\line(2,-1){40}} 

\put(30,8){\makebox(0,0){$-\alpha_1$}} 
\put(30,22){\makebox(0,0){$-\alpha_2$}} 
\put(30,32){\makebox(0,0){$-\alpha_3$}} 
\put(30,42){\makebox(0,0){$-\alpha_4$}} 
\put(30,52){\makebox(0,0){$-\alpha_5$}} 


\end{picture} 
\end{center} 
\caption{The ``snake'' in type $A_5$} 
\label{fig:octagon-snake} 
\end{figure}
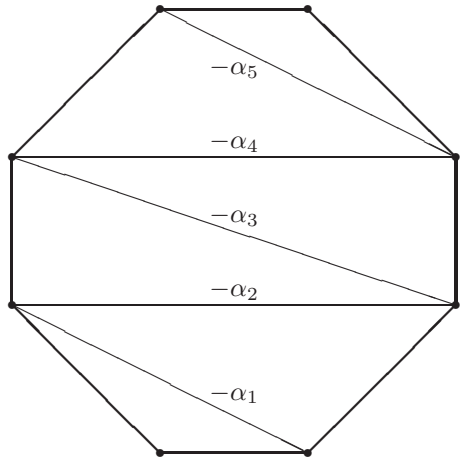 
\begin{figure}[ht] 
\begin{center} 
\setlength{\unitlength}{6pt} 
\begin{picture}(40,32)(-4,0) 
  \put(6,0){\line(1,0){20}} 
  \put(0,19){\line(1,0){32}} 
  \qbezier(6,0)(11,15.5)(16,31) 
  \qbezier(26,0)(21,15.5)(16,31) 
  \qbezier(6,0)(3,9.5)(0,19) 
  \qbezier(26,0)(29,9.5)(32,19) 
  \qbezier(0,19)(13,9.5)(26,0) 
  \qbezier(32,19)(19,9.5)(6,0) 
  \qbezier(0,19)(8,25)(16,31) 
  \qbezier(32,19)(24,25)(16,31) 

  \put(6,0){\circle*{.5}} 
  \put(26,0){\circle*{.5}} 
  \put(0,19){\circle*{.5}} 
  \put(32,19){\circle*{.5}} 
  \put(16,31){\circle*{.5}} 

\put(16,20.5){\makebox(0,0){$\alpha_1+\alpha_2$}} 
\put(23,16){\makebox(0,0){$-\alpha_2$}} 
\put(9,16){\makebox(0,0){$-\alpha_1$}} 
\put(12,8.7){\makebox(0,0){$\alpha_1$}} 
\put(20,8.7){\makebox(0,0){$\alpha_2$}} 

\end{picture} 
\end{center} 
\caption{Labeling of the diagonals in type $A_2$} \label{fig:a2} 
\end{figure}
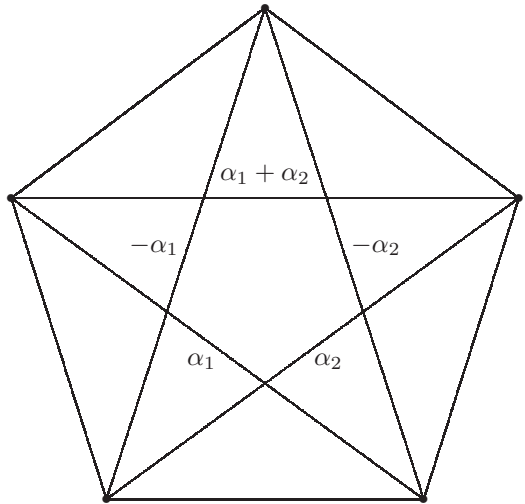


\section{Polytopal realizations of generalized associahedra}
\label{sec:polytopal}

We now demonstrate how to explicitly describe each generalized
associahedron by a set of linear inequalities. 

\begin{theorem} 
\label{th:tau-inv-support-functions} 
Suppose that a $(-w_\circ)$-invariant function $F: -\Pi \to \RR$ 
satisfies the inequalities 
\begin{equation*}
\label{eq:regular-dominant}
\sum_{i \in I} a_{ij} F(-\alpha_i) > 0 \quad \text{for all $j \in I$.} 
\end{equation*}
Let us extend $F$ (uniquely) to a $\langle \tau_-,\tau_+
\rangle$-invariant function on~$\Phi_{\geq -1}\,$. 
The generalized associahedron is then given in the dual
space~$Q_\reals^*$ by the linear inequalities 
\begin{equation*} 
\label{eq:inequalities-for-asso} 
\langle \mathbf{z}, \alpha \rangle\leq F(\alpha)\,,\ \text{for~all} \ 
\alpha\in\Phi_{\geq -1}\, . 
\end{equation*} 
\end{theorem} 

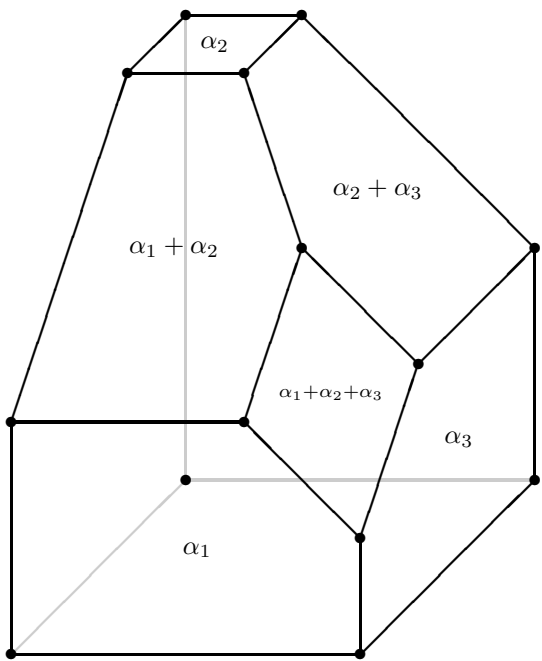
\begin{figure}[ht!]
\begin{center} 
\setlength{\unitlength}{2.2pt} 
\begin{picture}(90,110)(0,0) 
\thicklines

\verylight{
\put(30,30){\line(0,1){80}}
\put(30,30){\line(1,0){60}}
\put(0,0){\line(1,1){30}}
 }
\put(0,0){\line(1,0){60}} 
\put(0,0){\line(0,1){40}} 
\put(60,0){\line(0,1){20}} 
\put(60,0){\line(1,1){30}} 
\put(0,40){\line(1,0){40}} 
\put(0,40){\line(1,3){20}} 
\put(60,20){\line(-1,1){20}} 
\put(60,20){\line(1,3){10}} 
\put(90,30){\line(0,1){40}} 
\put(40,40){\line(1,3){10}} 
\put(70,50){\line(1,1){20}} 
\put(70,50){\line(-1,1){20}} 
\put(50,70){\line(-1,3){10}} 
\put(90,70){\line(-1,1){40}} 
\put(20,100){\line(1,0){20}} 
\put(20,100){\line(1,1){10}} 
\put(30,110){\line(1,0){20}} 
\put(40,100){\line(1,1){10}} 

\put(35,105){\makebox(0,0){$\alpha_2$}} 
\put(28,70){\makebox(0,0){$\alpha_1+\alpha_2$}} 
\put(63,80){\makebox(0,0){$\alpha_2+\alpha_3$}} 
\put(55,45){\makebox(0,0){$\scriptstyle \alpha_1+\alpha_2+\alpha_3$}} 
\put(32,18){\makebox(0,0){$\alpha_1$}} 
\put(77,37){\makebox(0,0){$\alpha_3$}}

\put(0,0){\circle*{2}} 
\put(60,0){\circle*{2}} 
\put(60,20){\circle*{2}} 
\put(30,30){\circle*{2}} 
\put(90,30){\circle*{2}} 
\put(0,40){\circle*{2}} 
\put(40,40){\circle*{2}} 
\put(70,50){\circle*{2}} 
\put(50,70){\circle*{2}} 
\put(90,70){\circle*{2}} 
\put(20,100){\circle*{2}} 
\put(40,100){\circle*{2}} 
\put(30,110){\circle*{2}} 
\put(50,110){\circle*{2}}

\end{picture} 
\end{center} 
\caption{Polytopal realization of 
the type $A_3$ associahedron}
\label{fig:A3asso-poly}
\end{figure}

An example of a function~$F$ satisfying the conditions in
Theorem~\ref{th:tau-inv-support-functions} is  obtained by setting 
$F(-\alpha_i)$ equal to the coefficient of the simple coroot~$\alpha_i^\vee$
  in the half-sum of all positive coroots.
(Coroots are the roots of the ``dual'' root system;
see \cite{Bourbaki, Humphreys}.) 
\begin{example}{\rm
In type~$A_3$, Theorem~\ref{th:tau-inv-support-functions} is illustrated in
Figure~\ref{fig:A3asso-poly}, which shows a $3$-dimensional
associahedron given by the inequalities 
\[ 
\begin{array}{rcl} 
\max(-z_1\,,\, 
-z_3\,,\, 
z_1\,,\, 
z_3\,,\, 
z_1+z_2\,,\, 
z_2+z_3 
)&\!\!\leq\!\!& 3/2\,,\\[.1in] 
\max(-z_2\,,\, 
z_2\,,\, 
z_1+z_2+z_3 
)&\!\!\leq\!\!& 2\,.
\end{array} 
\] 
}
\end{example}

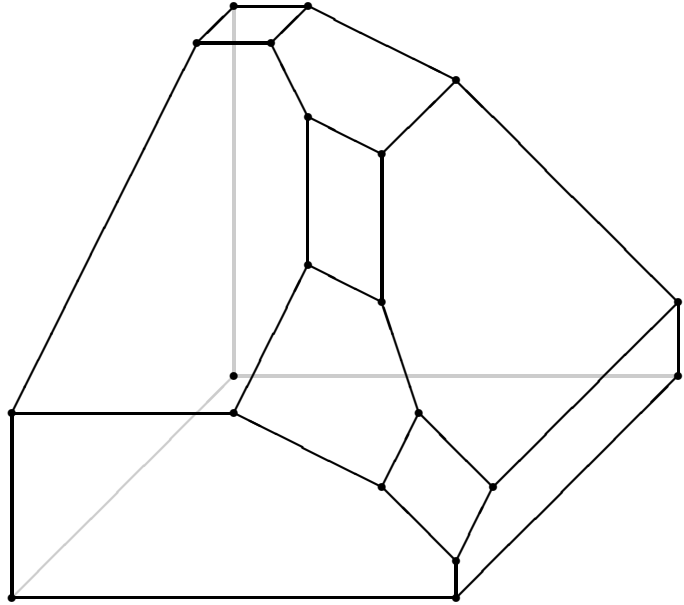
\begin{figure}[ht!]
\begin{center} 
\setlength{\unitlength}{1.4pt} 
\begin{picture}(180,170)(0,0) 
\thicklines

\verylight{
\put(0,0){\line(1,1){60}}
\put(60,60){\line(1,0){120}}
\put(60,60){\line(0,1){100}}
}

\put(0,0){\line(1,0){120}} 
\put(0,0){\line(0,1){50}} 
\put(120,0){\line(0,1){10}} 
\put(120,0){\line(1,1){60}} 
\put(0,50){\line(1,0){60}} 
\put(0,50){\line(1,2){50}} 
\put(120,10){\line(-1,1){20}} 
\put(120,10){\line(1,2){10}} 
\put(180,60){\line(0,1){20}} 
\put(60,50){\line(1,2){20}} 
\put(130,30){\line(1,1){50}} 
\put(130,30){\line(-1,1){20}} 
\put(100,30){\line(-2,1){40}} 
\put(100,30){\line(1,2){10}} 
\put(110,50){\line(-1,3){10}} 
\put(100,80){\line(-2,1){20}} 
\put(100,80){\line(0,1){40}} 
\put(80,90){\line(0,1){40}} 
\put(180,80){\line(-1,1){60}} 
\put(100,120){\line(-2,1){20}} 
\put(100,120){\line(1,1){20}} 
\put(80,130){\line(-1,2){10}} 
\put(50,150){\line(1,0){20}} 
\put(50,150){\line(1,1){10}} 
\put(70,150){\line(1,1){10}} 
\put(60,160){\line(1,0){20}} 
\put(80,160){\line(2,-1){40}} 


\put(0,0){\circle*{2}} 
\put(120,0){\circle*{2}} 
\put(120,10){\circle*{2}} 
\put(100,30){\circle*{2}} 
\put(130,30){\circle*{2}} 
\put(0,50){\circle*{2}} 
\put(60,50){\circle*{2}} 
\put(110,50){\circle*{2}} 
\put(60,60){\circle*{2}} 
\put(180,60){\circle*{2}} 
\put(100,80){\circle*{2}} 
\put(180,80){\circle*{2}} 
\put(80,90){\circle*{2}} 
\put(100,120){\circle*{2}} 
\put(80,130){\circle*{2}} 
\put(120,140){\circle*{2}} 
\put(50,150){\circle*{2}} 
\put(70,150){\circle*{2}} 
\put(60,160){\circle*{2}} 
\put(80,160){\circle*{2}}

\end{picture} 
\end{center} 
\caption{Polytopal realization of 
the type $C_3$ associahedron (cyclohedron)}
\label{fig:C3asso-poly}
\end{figure}

\begin{example}{\rm
In type~$C_3$, Theorem~\ref{th:tau-inv-support-functions} is illustrated in
Figure~\ref{fig:C3asso-poly}
that shows a $3$-dimensional cyclohedron given by the inequalities
\[ 
\!\!\!\!\!\!
\begin{array}{rcl} 
\max(-z_1\,,\, 
z_1\,,\, 
z_1+z_2\,,\, 
z_2+z_3 
)&\!\!\leq\!\!& 5/2\,,\\[.1in] 
\max(-z_2\,,\, 
z_2\,,\, 
z_1+z_2+z_3\,,\, 
z_1+2z_2+z_3 
)&\!\!\leq\!\!& 4\,,\\[.1in] 
\max(-z_3\,,\, 
z_3\,,\, 
2z_2+z_3\,,\, 
2z_1+2z_2+z_3 
)&\!\!\leq\!\!& 9/2\,.
\end{array} 
\] 
}
\end{example}

\pagebreak

\section{Double wiring diagrams and double Bruhat cells} 
\label{sec:3by3}

The goal of this section is to give a glimpse into how cluster
algebras come up in ``real life.''
We will present just one example:
the coordinate ring of the open double Bruhat cell
in~$GL_n(\complexes)$. 

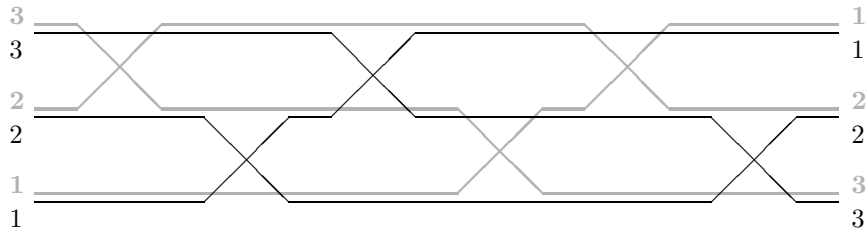
\begin{figure}[ht]
\setlength{\unitlength}{1.6pt} 
\begin{center}
\begin{picture}(180,45)(6,0)

\thicklines
\light{

  \put(0,2){\line(1,0){100}}
  \put(120,2){\line(1,0){70}}
  \put(0,22){\line(1,0){10}}
  \put(30,22){\line(1,0){70}}
  \put(120,22){\line(1,0){10}}
  \put(150,22){\line(1,0){40}}
  \put(0,42){\line(1,0){10}}
  \put(30,42){\line(1,0){100}}
  \put(150,42){\line(1,0){40}}

  \put(10,22){\line(1,1){20}}
  \put(100,2){\line(1,1){20}}
  \put(130,22){\line(1,1){20}}

  \put(10,42){\line(1,-1){20}}
  \put(100,22){\line(1,-1){20}}
  \put(130,42){\line(1,-1){20}}

  \put(193,2){$\mathbf{3}$}
  \put(193,22){$\mathbf{2}$}
  \put(193,42){$\mathbf{1}$}

  \put(-6,2){$\mathbf{1}$}
  \put(-6,22){$\mathbf{2}$}
  \put(-6,42){$\mathbf{3}$}
}
\dark{
\thinlines
  \put(0,0){\line(1,0){40}}
  \put(60,0){\line(1,0){100}}
  \put(180,0){\line(1,0){10}}
  \put(0,20){\line(1,0){40}}
  \put(60,20){\line(1,0){10}}
  \put(90,20){\line(1,0){70}}
  \put(180,20){\line(1,0){10}}
  \put(0,40){\line(1,0){70}}
  \put(90,40){\line(1,0){100}}

  \put(40,0){\line(1,1){20}}
  \put(70,20){\line(1,1){20}}
  \put(160,0){\line(1,1){20}}

  \put(40,20){\line(1,-1){20}}
  \put(70,40){\line(1,-1){20}}
  \put(160,20){\line(1,-1){20}}
  \put(193,-6){$3$}
  \put(193,14){$2$}
  \put(193,34){$1$}

  \put(-6,-6){$1$}
  \put(-6,14){$2$}
  \put(-6,34){$3$}
}
\end{picture}
\end{center}
\caption{Double wiring diagram}
\label{fig:double-wiring}
\end{figure}

We will need the notion of a \emph{double wiring diagram}
(of type~$(w_\circ,w_\circ)$), 
which is illustrated  in Figure~\ref{fig:double-wiring}.
Such a diagram consists of two families of $n$
piecewise-straight 
lines, each family colored with one of two colors.  
The crucial requirement is that each
pair of 
lines of like color intersect exactly once. 
The lines in a double wiring diagram are numbered separately within
each color, as shown in Figure~\ref{fig:double-wiring}. 

We note in passing that double wiring diagrams correspond naturally to
shuffles of two reduced words for the element~$w_\circ$ in the symmetric
group~$\mathcal{S}_n$. 

\begin{figure}[ht]
\setlength{\unitlength}{1.6pt} 
\begin{center}
\begin{picture}(180,45)(6,0)

\light{
\thicklines

  \put(0,2){\line(1,0){100}}
  \put(120,2){\line(1,0){70}}
  \put(0,22){\line(1,0){40}}
  \put(60,22){\line(1,0){40}}
  \put(120,22){\line(1,0){40}}
  \put(180,22){\line(1,0){10}}
  \put(0,42){\line(1,0){40}}
  \put(60,42){\line(1,0){100}}
  \put(180,42){\line(1,0){10}}

  \put(40,22){\line(1,1){20}}
  \put(100,2){\line(1,1){20}}
  \put(160,22){\line(1,1){20}}

  \put(40,42){\line(1,-1){20}}
  \put(100,22){\line(1,-1){20}}
  \put(160,42){\line(1,-1){20}}

  \put(193,2){$\mathbf{1}$}
  \put(193,22){$\mathbf{2}$}
  \put(193,42){$\mathbf{3}$}

  \put(-6,2){$\mathbf{3}$}
  \put(-6,22){$\mathbf{2}$}
  \put(-6,42){$\mathbf{1}$}
}
\thinlines
\dark{
  \put(0,0){\line(1,0){10}}
  \put(30,0){\line(1,0){100}}
  \put(150,0){\line(1,0){40}}
  \put(0,20){\line(1,0){10}}
  \put(30,20){\line(1,0){40}}
  \put(90,20){\line(1,0){40}}
  \put(150,20){\line(1,0){40}}
  \put(0,40){\line(1,0){70}}
  \put(90,40){\line(1,0){100}}

  \put(10,0){\line(1,1){20}}
  \put(70,20){\line(1,1){20}}
  \put(130,0){\line(1,1){20}}

  \put(10,20){\line(1,-1){20}}
  \put(70,40){\line(1,-1){20}}
  \put(130,20){\line(1,-1){20}}
  \put(193,-6){${3}$}
  \put(193,14){${2}$}
  \put(193,34){${1}$}

  \put(-6,-6){${1}$}
  \put(-6,14){${2}$}
  \put(-6,34){${3}$}
}
\end{picture}
\end{center}
\caption{An isotopic double wiring diagram}
\label{fig:isotopic-double-wiring}
\end{figure}
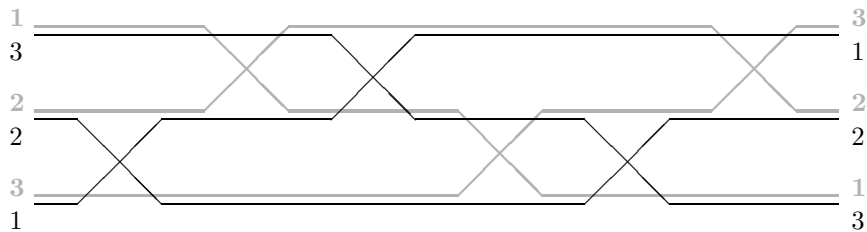

 From now on, we will not distinguish between double
wiring diagrams that are \emph{isotopic}, i.e., have the same
``topology.''
For example, the diagrams in Figures~\ref{fig:double-wiring}
and~\ref{fig:isotopic-double-wiring} are isotopic to each other. 
The diagram in Figure~\ref{fig:isotopic-double-wiring} is obtained
from Figure~\ref{fig:double-wiring}
by sliding the two leftmost crossings past each other,
and also doing the same for the two rightmost crossings. 

The following lemma is a direct corollary of a theorem of G.~Ringel
(1956). It can also be obtained from the type~$A$ version of a
classical result by J.~Tits (1969) concerning the word problem
in Coxeter groups. 

\begin{lemma}
\label{lem:ringel-tits}
Any two (isotopy classes of) double wiring diagrams can be
transformed into each other by a sequence of local ``moves'' of
three different kinds, shown in Figure~\ref{fig:moves}. 
(Each of these local moves only changes a small portion of a 
double wiring diagram, leaving the rest of it intact.) 
\end{lemma}

The reader is asked to ignore, for now, the labels
$A,B,\dots,Z$ in Figure~\ref{fig:moves}. 


\begin{figure}[ht]
\setlength{\unitlength}{1pt} 
\begin{center}
\begin{picture}(120,45)(10,0)
\thicklines
\light{
  \put(0,2){\line(1,0){100}}
  \put(0,22){\line(1,0){100}}
  \put(0,42){\line(1,0){100}}
}
\thinlines
\dark{
  \put(0,0){\line(1,0){10}}
  \put(30,0){\line(1,0){40}}
  \put(50,0){\line(1,0){10}}
  \put(90,0){\line(1,0){10}}
  \put(0,20){\line(1,0){10}}
  \put(30,20){\line(1,0){10}}
  \put(60,20){\line(1,0){10}}
  \put(90,20){\line(1,0){10}}
  \put(0,40){\line(1,0){40}}
  \put(60,40){\line(1,0){40}}

  \put(10,0){\line(1,1){20}}
  \put(40,20){\line(1,1){20}}
  \put(70,0){\line(1,1){20}}

  \put(10,20){\line(1,-1){20}}
  \put(40,40){\line(1,-1){20}}
  \put(70,20){\line(1,-1){20}}
}
  \put(115,24){\vector(1,0){10}}
  \put(125,16){\vector(-1,0){10}}



  \put(2,8){$A$}
  \put(46,8){$Y$}
  \put(89,8){$D$}

  \put(16,28){$B$}
  \put(76,28){$C$}

\end{picture}
\begin{picture}(100,40)(-5,0)

\thicklines
\light{
  \put(0,2){\line(1,0){100}}
  \put(0,22){\line(1,0){100}}
  \put(0,42){\line(1,0){100}}
}
\thinlines
\dark{

  \put(0,40){\line(1,0){10}}
  \put(30,40){\line(1,0){40}}
  \put(50,40){\line(1,0){10}}
  \put(90,40){\line(1,0){10}}
  \put(0,20){\line(1,0){10}}
  \put(30,20){\line(1,0){10}}
  \put(60,20){\line(1,0){10}}
  \put(90,20){\line(1,0){10}}
  \put(0,0){\line(1,0){40}}
  \put(60,0){\line(1,0){40}}

  \put(10,20){\line(1,1){20}}
  \put(40,0){\line(1,1){20}}
  \put(70,20){\line(1,1){20}}

  \put(10,40){\line(1,-1){20}}
  \put(40,20){\line(1,-1){20}}
  \put(70,40){\line(1,-1){20}}
}
  \put(2,28){$B$}
  \put(46,28){$Z$}
  \put(89,28){$C$}

  \put(16,8){$A$}
  \put(76,8){$D$}

\end{picture}
\end{center}
\begin{center}
\begin{picture}(120,60)(10,0)
\thicklines
\light{

  \put(0,2){\line(1,0){10}}
  \put(30,2){\line(1,0){40}}
  \put(50,2){\line(1,0){10}}
  \put(90,2){\line(1,0){10}}
  \put(0,22){\line(1,0){10}}
  \put(30,22){\line(1,0){10}}
  \put(60,22){\line(1,0){10}}
  \put(90,22){\line(1,0){10}}
  \put(0,42){\line(1,0){40}}
  \put(60,42){\line(1,0){40}}

  \put(10,2){\line(1,1){20}}
  \put(40,22){\line(1,1){20}}
  \put(70,2){\line(1,1){20}}

  \put(10,22){\line(1,-1){20}}
  \put(40,42){\line(1,-1){20}}
  \put(70,22){\line(1,-1){20}}
}
  \put(2,8){$A$}
  \put(46,8){$Y$}
  \put(89,8){$D$}

  \put(16,28){$B$}
  \put(76,28){$C$}

\thinlines
  \put(115,24){\vector(1,0){10}}
  \put(125,16){\vector(-1,0){10}}
\dark{
  \put(0,0){\line(1,0){100}}
  \put(0,20){\line(1,0){100}}
  \put(0,40){\line(1,0){100}}
}

\end{picture}
\begin{picture}(100,40)(-5,0)
\thicklines
\light{

  \put(0,42){\line(1,0){10}}
  \put(30,42){\line(1,0){40}}
  \put(50,42){\line(1,0){10}}
  \put(90,42){\line(1,0){10}}
  \put(0,22){\line(1,0){10}}
  \put(30,22){\line(1,0){10}}
  \put(60,22){\line(1,0){10}}
  \put(90,22){\line(1,0){10}}
  \put(0,2){\line(1,0){40}}
  \put(60,2){\line(1,0){40}}

  \put(10,22){\line(1,1){20}}
  \put(40,2){\line(1,1){20}}
  \put(70,22){\line(1,1){20}}

  \put(10,42){\line(1,-1){20}}
  \put(40,22){\line(1,-1){20}}
  \put(70,42){\line(1,-1){20}}
}

\thinlines
\dark{
  \put(0,0){\line(1,0){100}}
  \put(0,20){\line(1,0){100}}
  \put(0,40){\line(1,0){100}}
}
  \put(2,28){$B$}
  \put(46,28){$Z$}
  \put(89,28){$C$}

  \put(16,8){$A$}
  \put(76,8){$D$}

\end{picture}
\end{center}
\begin{center}
\begin{picture}(100,60)(5,-10)

\thicklines
\light{

  \put(0,2){\line(1,0){40}}
  \put(60,2){\line(1,0){10}}

  \put(0,22){\line(1,0){40}}
  \put(60,22){\line(1,0){10}}

  \put(40,2){\line(1,1){20}}
  \put(40,22){\line(1,-1){20}}
}
\thinlines
\dark{

  \put(0,0){\line(1,0){10}}
  \put(30,0){\line(1,0){40}}

  \put(0,20){\line(1,0){10}}
  \put(30,20){\line(1,0){40}}

  \put(10,0){\line(1,1){20}}
  \put(10,20){\line(1,-1){20}}
}
  \put(2,8){$A$}
  \put(32,8){$Y$}
  \put(62,8){$C$}

  \put(32,28){$B$}
  \put(32,-12){$D$}

\thinlines
  \put(100,15){\vector(1,0){10}}
  \put(110,7){\vector(-1,0){10}}

\end{picture}
\begin{picture}(100,60)(-30,-10)

\thicklines
\light{

  \put(0,2){\line(1,0){10}}
  \put(30,2){\line(1,0){40}}

  \put(0,22){\line(1,0){10}}
  \put(30,22){\line(1,0){40}}

  \put(10,2){\line(1,1){20}}
  \put(10,22){\line(1,-1){20}}
}
\thinlines
\dark{
  \put(0,0){\line(1,0){40}}
  \put(60,0){\line(1,0){10}}

  \put(0,20){\line(1,0){40}}
  \put(60,20){\line(1,0){10}}

  \put(40,0){\line(1,1){20}}
  \put(40,20){\line(1,-1){20}}
}
  \put(2,8){$A$}
  \put(32,8){$Z$}
  \put(62,8){$C$}

  \put(32,28){$B$}
  \put(32,-12){$D$}

\end{picture}
\end{center}
\caption{Local ``moves''}
\label{fig:moves}
\end{figure}


To illustrate Lemma~\ref{lem:ringel-tits},
the double wiring diagram in  Figure~\ref{fig:double-wiring} 
allows $4$~different local moves,
all of which are of the kind shown at the bottom 
of Figure~\ref{fig:moves}. 
Two of these moves can be performed by 
first passing to the isotopic
Figure~\ref{fig:isotopic-double-wiring}. 
To make each of the other two moves,
slide the two innermost crossings
in  Figure~\ref{fig:double-wiring} past each other;
this will create two patterns of the form shown 
at the bottom of Figure~\ref{fig:moves}.

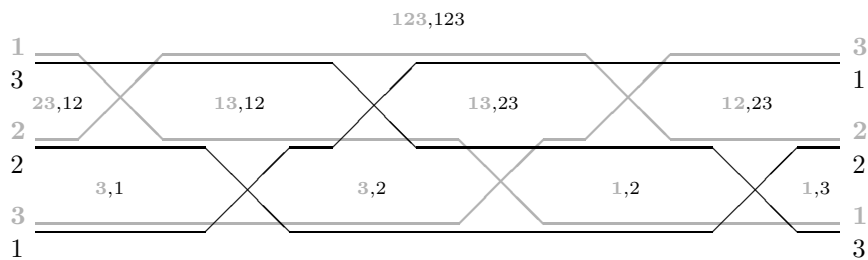
\begin{figure}[ht]
\setlength{\unitlength}{1.6pt} 
\begin{center}
\begin{picture}(180,50)(6,0)

\light{
\thicklines

  \put(0,2){\line(1,0){100}}
  \put(120,2){\line(1,0){70}}
  \put(0,22){\line(1,0){10}}
  \put(30,22){\line(1,0){70}}
  \put(120,22){\line(1,0){10}}
  \put(150,22){\line(1,0){40}}
  \put(0,42){\line(1,0){10}}
  \put(30,42){\line(1,0){100}}
  \put(150,42){\line(1,0){40}}

  \put(10,22){\line(1,1){20}}
  \put(100,2){\line(1,1){20}}
  \put(130,22){\line(1,1){20}}

  \put(10,42){\line(1,-1){20}}
  \put(100,22){\line(1,-1){20}}
  \put(130,42){\line(1,-1){20}}

  \put(193,2){$\mathbf{{1}}$}
  \put(193,22){$\mathbf{{2}}$}
  \put(193,42){$\mathbf{{3}}$}

  \put(-6,2){$\mathbf{{3}}$}
  \put(-6,22){$\mathbf{{2}}$}
  \put(-6,42){$\mathbf{{1}}$}


  \put(14,10){$_{\light{\mathbf{3}},{\dark{1}}}$}
  \put(76,10){$_{\light{\mathbf{3}},{\dark{2}}}$}
  \put(136,10){$_{\light{\mathbf{1}},{\dark{2}}}$}
  \put(181,10){$_{\light{\mathbf{1}},{\dark{3}}}$}

  \put(-1,30){$_{\light{\mathbf{23}},{\dark{12}}}$}
  \put(42,30){$_{\light{\mathbf{13}},{\dark{12}}}$}
  \put(102,30){$_{\light{\mathbf{13}},{\dark{23}}}$}
  \put(162,30){$_{\light{\mathbf{12}},{\dark{23}}}$}

  \put(84,50){$_{\light{\mathbf{123}},{\dark{123}}}$}
}
\thinlines
\dark{

  \put(0,0){\line(1,0){40}}
  \put(60,0){\line(1,0){100}}
  \put(180,0){\line(1,0){10}}
  \put(0,20){\line(1,0){40}}
  \put(60,20){\line(1,0){10}}
  \put(90,20){\line(1,0){70}}
  \put(180,20){\line(1,0){10}}
  \put(0,40){\line(1,0){70}}
  \put(90,40){\line(1,0){100}}

  \put(40,0){\line(1,1){20}}
  \put(70,20){\line(1,1){20}}
  \put(160,0){\line(1,1){20}}

  \put(40,20){\line(1,-1){20}}
  \put(70,40){\line(1,-1){20}}
  \put(160,20){\line(1,-1){20}}

  \put(193,-6){${3}$}
  \put(193,14){${2}$}
  \put(193,34){${1}$}

  \put(-6,-6){${1}$}
  \put(-6,14){${2}$}
  \put(-6,34){${3}$}
}
\end{picture}
\end{center}
\caption{Chamber minors}
\label{fig:chamber-sets}
\end{figure}

A \emph{chamber} of a double wiring diagram is a connected component
of the complement to the union of the lines, with the exception of the
``crumbs'' made of narrow horizontal isthmuses and small triangular
regions; the large component at the very bottom is not included either. 
With these conventions, there are exactly $n^2$ chambers altogether
(e.g., $9$~chambers in Figure~\ref{fig:double-wiring}). 
We then assign to every chamber a pair of subsets
of the set $[1,n]=\{1,\dots,n\}$: 
each subset indicates which lines of the corresponding
color pass \emph{below} that chamber; 
see Figure~\ref{fig:chamber-sets}. 

Suppose we are given an $n\times n$ matrix~$x=(x_{ij})$. 
For any subsets $I,J\subset\{1,\dots,n\}$ of equal cardinality, 
we denote by $\Delta_{I,J}(x)$ the corresponding minor of~$x$, 
that is, the determinant of the submatrix of~$x$ 
occupying the rows and columns specified by the sets $I$ and~$J$. 
Then each chamber of a double wiring diagram is naturally associated
with a \emph{chamber minor} $\Delta_{I,J}$ (viewed as a function on
the general linear group~$GL_n(\complexes)$), 
where $I$ and $J$ are the sets written into that chamber. 
%

We note that two double wiring diagrams have the same associated
collections of chamber minors if and only if they are isotopic. 

Let $\mathcal{F}$ denote the field of rational functions
on~$GL_n(\complexes)$, i.e., the field of rational functions with
complex coefficients in the
matrix entries~$x_{ij}$ (viewed as indeterminates). 

\begin{lemma}
\label{lem:ch-min-gen}
The $n^2$ chamber minors of an arbitrary double wiring diagram 
form a set of algebraically independent generators of the
field~$\mathcal{F}$. 
\end{lemma}

Notice that each local move in Figure~\ref{fig:moves}
exchanges a single chamber minor~$Y$
(associated with a \emph{bounded}, or interior, chamber) with another
chamber minor~$Z$, and keeps all other chamber minors in place. 
We can therefore define, by analogy with triangulations, a graph of
exchanges whose vertices correspond to (isotopy classes of) double
wiring diagrams, and whose edges correspond to the moves in
Figure~\ref{fig:moves}. 

\begin{figure}[ht]
\setlength{\unitlength}{1.6pt} 
\begin{center}
\begin{picture}(180,55)(6,0)

\light{
\thicklines

  \put(0,2){\line(1,0){10}}
  \put(30,2){\line(1,0){40}}
  \put(90,2){\line(1,0){100}}
  \put(0,22){\line(1,0){10}}
  \put(30,22){\line(1,0){10}}
  \put(60,22){\line(1,0){10}}
  \put(90,22){\line(1,0){100}}
  \put(0,42){\line(1,0){40}}
  \put(60,42){\line(1,0){130}}

  \put(10,2){\line(1,1){20}}
  \put(40,22){\line(1,1){20}}
  \put(70,2){\line(1,1){20}}

  \put(10,22){\line(1,-1){20}}
  \put(40,42){\line(1,-1){20}}
  \put(70,22){\line(1,-1){20}}

  \put(193,2){$\mathbf{1}$}
  \put(193,22){$\mathbf{2}$}
  \put(193,42){$\mathbf{3}$}

  \put(-6,2){$\mathbf{3}$}
  \put(-6,22){$\mathbf{2}$}
  \put(-6,42){$\mathbf{1}$}

  \put(2,10){$_{\light{\mathbf{3}},{\dark{1}}}$}
  \put(46,10){$_{\light{\mathbf{2}},{\dark{1}}}$}
  \put(91,10){$_{\light{\mathbf{1}},{\dark{1}}}$}
  \put(136,10){$_{\light{\mathbf{1}},{\dark{2}}}$}
  \put(181,10){$_{\light{\mathbf{1}},{\dark{3}}}$}

  \put(12,30){$_{\light{\mathbf{23}},{\dark{12}}}$}
  \put(87,30){$_{\light{\mathbf{12}},{\dark{12}}}$}
  \put(162,30){$_{\light{\mathbf{12}},{\dark{23}}}$}

  \put(84,50){$_{\light{\mathbf{123}},{\dark{123}}}$}
}
\thinlines
\dark{

  \put(0,0){\line(1,0){100}}
  \put(120,0){\line(1,0){40}}
  \put(180,0){\line(1,0){10}}
  \put(0,20){\line(1,0){100}}
  \put(120,20){\line(1,0){10}}
  \put(150,20){\line(1,0){10}}
  \put(180,20){\line(1,0){10}}
  \put(0,40){\line(1,0){130}}
  \put(150,40){\line(1,0){40}}

  \put(100,0){\line(1,1){20}}
  \put(130,20){\line(1,1){20}}
  \put(160,0){\line(1,1){20}}

  \put(100,20){\line(1,-1){20}}
  \put(130,40){\line(1,-1){20}}
  \put(160,20){\line(1,-1){20}}

  \put(193,-6){${3}$}
  \put(193,14){${2}$}
  \put(193,34){${1}$}

  \put(-6,-6){${1}$}
  \put(-6,14){${2}$}
  \put(-6,34){${3}$}
}
\end{picture}
\end{center}
\caption{A double wiring diagram allowing $3$ local moves}
\label{fig:lexmindiagram}
\end{figure}
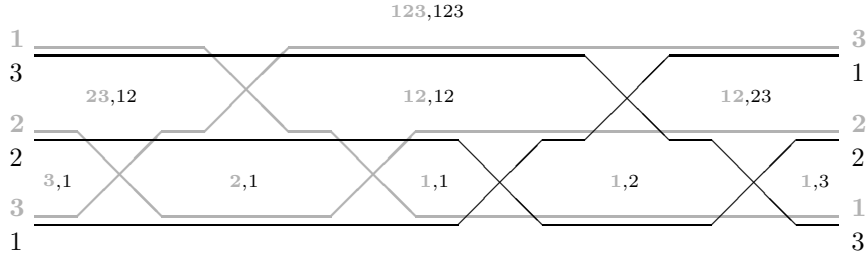

\begin{example}
\label{example:gl3-34}
{\rm
For $n=3$, there are $34$ non-isotopic double wiring diagrams. 
The corresponding $34$-vertex graph of exchanges can be found in
\cite[Figure~10]{tptp}. 
It has $18$ vertices of degree~$4$, and $16$ vertices of degree~$3$. 
They correspond, respectively, to the double wiring diagrams
that allow $4$ local moves (as the diagram in
Figure~\ref{fig:chamber-sets}) 
and those allowing only $3$ local moves (as the diagram in
Figure~\ref{fig:lexmindiagram}). 
}
\end{example}

\begin{lemma}
\label{lem:3-term}
Whenever two double wiring diagrams differ by a single local move of
one of the three types shown in Figure~\ref{fig:moves},
the chamber minors 
appearing there satisfy the identity $AC+BD=YZ$. 
\end{lemma}

Lemmas~\ref{lem:ch-min-gen} and~\ref{lem:3-term} suggest the existence
of a cluster algebra 
structure associated with $n\times n$ matrices. 
We next present one of several versions of this structure, 
leaving out most of the technical details. 
The ambient field for our cluster algebra is the field $\mathcal{F}$ 
of rational functions 
on~$GL_n(\complexes)$ introduced above. 
Each double wiring diagram provides us with a seed
whose cluster variables are the $(n-1)^2$ chamber minors associated with the
bounded chambers; 
the frozen variables are the $2n-1$ chamber minors associated with the
unbounded chambers at the edges of the diagram. 
It remains to define the matrices~$\tilde B$. 

Take any double wiring diagram in which every bounded chamber can be
``flipped'' (such a diagram can be constructed for any~$n$). 
Comparing the corresponding exchange relations $AC+BD=YZ$
with~\eqref{eq:exch-rel-gen}, determine the matrix entries of~$\tilde
B$. It can be shown that exchanges associated with the local moves on 
double wiring diagrams are compatible with the cluster algebra axioms.
Furthermore, applying these axioms uncovers hitherto hidden clusters 
which do not correspond to any wiring diagrams.
Each variable in these clusters is a regular function on~$GL_n(\complexes)$
(a polynomial in the matrix entries). 
The resulting cluster algebra coincides with the coordinate ring of
the \emph{open double Bruhat cell} $G^{w_\circ,w_\circ}$
in~$GL_n(\complexes)$.
We refer to~\cite{ca3} for further details. 

\begin{example}
{\rm
The open double Bruhat cell $G^{w_\circ,w_\circ}\subset GL_3(\complexes)$
consists of all complex $3\times 3$ matrices $x=(x_{ij})$
whose minors
\begin{equation}
\label{eq:sl3-w0w0-nonvanishing}
x_{13},
{\ \ }
\Bigl|\!\Bigl|\begin{array}{cc}
x_{12} & x_{13} \\
x_{22} & x_{23}
\end{array}\Bigr|\!\Bigr|,
{\ \ }
x_{31},
{\ \ }
\Bigl|\!\Bigl|\begin{array}{cc}
x_{21} & x_{22} \\
x_{31} & x_{32}
\end{array}\Bigr|\!\Bigr|,
{\ \ }
\det(x) 
\end{equation}
are nonzero.
(These $5$ minors correspond to the unbounded chambers of any double
wiring diagram for~$GL_3(\complexes)$.) 
The coordinate ring $\complexes[G^{w_\circ,w_\circ}]$ turns out to be
a cluster algebra of type~$D_4$
over the ground ring generated by the minors in~\eqref{eq:sl3-w0w0-nonvanishing} and
their inverses.
Thus, the ring of rational functions on~$GL_3$ exhibits some quite
unexpected symmetries of type~$D_4$.  

This cluster algebra has $16$~cluster variables, corresponding to the $16$ roots
in~$\Phi_{\geq -1}$.
These variables are:

\begin{itemize}
\item 
$14$ (among the $19$ total) minors of~$x$, namely,
all except those listed in~\eqref{eq:sl3-w0w0-nonvanishing};

\item 
two ``hidden'' variables: 
$x_{12}x_{21}x_{33}-x_{12}x_{23}x_{31}-x_{13}x_{21}x_{32}+x_{13}x_{22}x_{31}$
and
$x_{11}x_{23}x_{32}-x_{12}x_{23}x_{31}-x_{13}x_{21}x_{32}+x_{13}x_{22}x_{31}$.
\end{itemize}
These $16$ variables form $50$~clusters of size~$4$, one for each of
the $50$ vertices of the type~$D_4$ associahedron.
}
\end{example}

For any $n\geq 4$, the construction described above produces a cluster
algebra of infinite type.

\lecture{Enumerative Problems}
\label{lec:num}

\section{Catalan combinatorics of arbitrary type 
}
\label{sec:numerology}

Let $\Phi$ be a finite irreducible crystallographic root system of rank~$n$,
and $W$ the corresponding reflection group.  
We retain the root-theoretic notation used in Lectures~\ref{lec2}
and~\ref{lec:cluster}. 
In particular, $e_1, \dots, e_n$ are the exponents of~$\Phi$,
and $h$ is the Coxeter number.

The number of vertices of an $n$-dimensional associahedron
(or, equivalently, the number of clusters in a cluster algebra of
type~$A_n$) is the Catalan number $\frac{1}{n+2}
\binom{2n+2}{n+1}$.
It is natural to ask similar enumerative questions for other
Cartan-Killing types. 

\begin{theorem}[\cite{ga}]
\label{th:N-thru-exponents}
The number of clusters in a cluster algebra of finite type
associated with a root system~$\Phi$
(or, equivalently, the number of vertices of the corresponding
generalized associahedron) is equal to
\begin{equation}
\label{eq:N-thru-exponents}
N(\Phi) \stackrel{\rm def}{=} \prod_{i=1}^n \frac{e_i + h + 1}{e_i +
  1} \,. 
\end{equation}
\end{theorem}

Figure~\ref{W-Cat} shows the values of $N(\Phi)$ for all~$\Phi$.
Recall that the exponents of root systems are tabulated in
Figure~\ref{data}.

\begin{figure}[ht]
\begin{center}
\begin{tabular}{|c|c|c|c|c|c|c|c|c|}
\hline
$A_n$ &
$B_n,C_n
$ & $D_n$ & $E_6$ & $E_7$ & $E_8$ & $F_4$  & $G_2$ \\
\hline
&&&&&&&\\[-.1in]
$\textstyle\frac{1}{n+2} \binom{2n+2}{n+1}$&
$\textstyle\binom{2n}{n}$ 
&   $\textstyle\frac{3n-2}{n}\binom{2n-2}{n-1}$ 
&
$\textstyle 833$
&
$\textstyle 4160$
&
$\textstyle 25080$
&
$\textstyle 105$
& $\textstyle 8$ \\[.05in]
\hline
\end{tabular}
\end{center}
\caption{The numbers $N(\Phi)$ 
} 
\label{W-Cat}
\end{figure}

As the numbers $N(\Phi)$ given by \eqref{eq:N-thru-exponents}
can be thought of as generalizations of the
Catalan numbers to an arbitrary Cartan-Killing type, 
it comes as no surprise that they count a host of various 
combinatorial objects related to the root system~$\Phi$. 
Below in this section, 
we briefly describe several families of objects counted by~$N(\Phi)$.
We refer the reader to the introductory sections
of~\cite{Ath,Ath-Rei,Ath-TAMS,chapoton-enum,panyushev} 
for the history of research in this area, 
for further details and references, and for numerous generalizations and connections. 

The numbers $N(\Phi)$ 
seem to have first appeared in D.~Djokovi\'c's
work \cite{djokovic}  on enumeration of conjugacy classes of elements of finite
\hbox{order in Lie groups.}

\subsection*{Antichains in the root poset (non-nesting
partitions)}  

\begin{figure}[ht!]
\centerline{
        \epsfbox{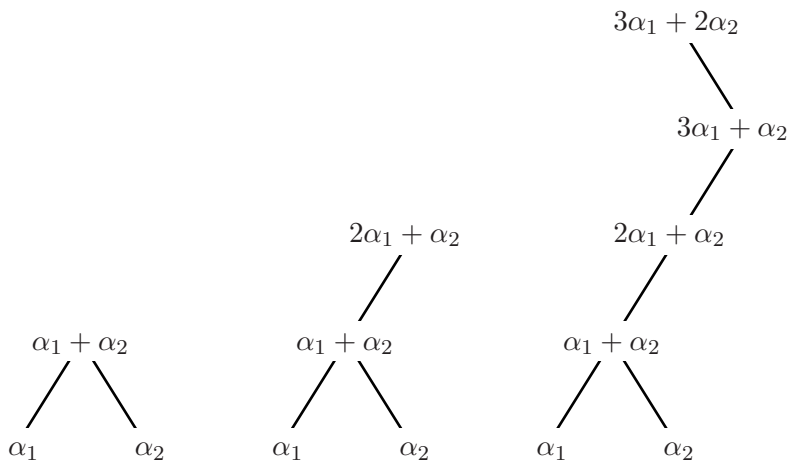}
        \begin{picture}(0,0)(264,-5)
                \put(-32,-3){\large$\alpha_1$}
                \put(16,-3){\large$\alpha_2$}
                \put(-23,37){\large$\alpha_1+\alpha_2$}
                \put(68,-3){\large$\alpha_1$}
                \put(116,-3){\large$\alpha_2$}
                \put(77,37){\large$\alpha_1+\alpha_2$}
                \put(97,78){\large$2\alpha_1+\alpha_2$}
                \put(168,-3){\large$\alpha_1$}
                \put(216,-3){\large$\alpha_2$}
                \put(178,37){\large$\alpha_1+\alpha_2$}
                \put(197,78){\large$2\alpha_1+\alpha_2$}
                \put(221,118){\large$3\alpha_1+\alpha_2$}
                \put(197,158){\large$3\alpha_1+2\alpha_2$}
        \end{picture}
}
\caption{The root posets of types $A_2$, $B_2$ and $G_2$.}
\label{rootposet}
\end{figure}

\begin{figure}[ht!]
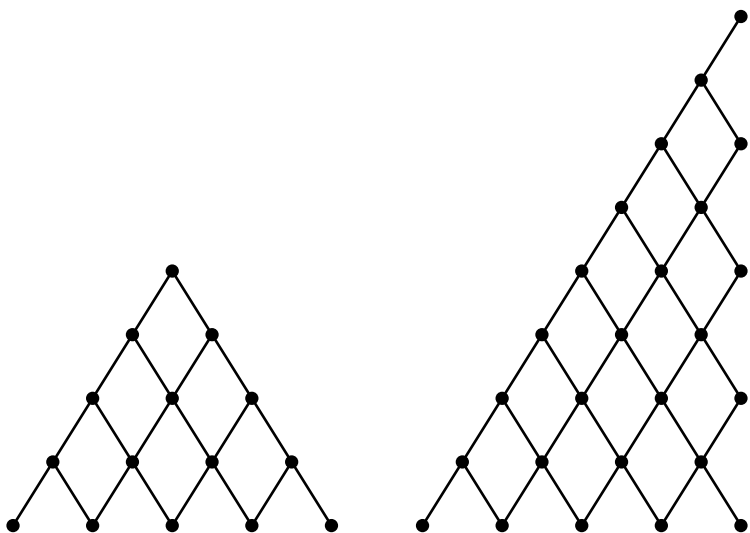

\centerline{
        \epsfbox{A5rootposet.ps}
        \hspace{.3 in}
        \epsfbox{B5rootposet.ps}
}
\caption{The root posets of types $A_5$ and $B_5$.}
\label{rootposetAB5}
\end{figure}

The {\em root poset} of $\Phi$ is the partial order
on the set of positive roots $\Phi_+$ such that $\beta\le\gamma$ if and only if
$\gamma-\beta$ is a nonnegative (integer) linear combination of simple
roots. 
See Figures~\ref{rootposet} and~\ref{rootposetAB5}.

\begin{theorem}[\cite{cellini-papi,Reiner-noncrossing,shi}]
The number of antichains (i.e., sets of pairwise non-comparable
elements) in the root poset of $\Phi$  is equal to~$N(\Phi)$. 
\end{theorem}


\subsection*{Positive regions of the Shi arrangement}

The \emph{Shi arrangement} is the arrangement of affine hyperplanes
defined by the equations 
\[
\begin{array}{r}
\br{\beta,x}=0\\
\br{\beta,x}=1
\end{array}
\quad 
\text{for all $\beta\in\Phi_+$.}
\]
(Thus, the number of hyperplanes in the Shi arrangement is equal to
the number of roots in the root system~$\Phi$.) 
The \emph{positive regions} of this arrangement are the regions contained
in the {\em positive cone}, which consists of the points~$x$ 
such that $\br{\beta,x}>0$ for any $\beta\in\Phi_+\,$.

\begin{theorem}[\cite{shi}] 
The number of positive regions in the Shi arrangement is equal to~$N(\Phi)$. 
\end{theorem}

\begin{figure}[ht!]
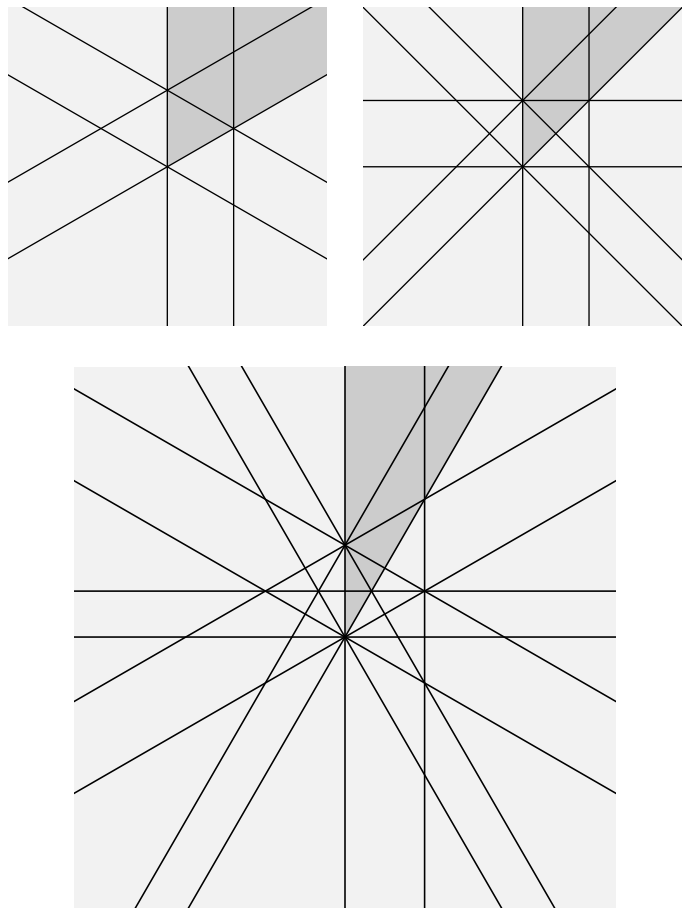

\centerline{
{\epsfbox{A2shi.ps}}
\hspace{.1 in}
{\epsfbox{B2shi.ps}}
}
\vspace{.2 in}
\centerline{\scalebox{1.2}{\epsfbox{G2shi.ps}}}
\caption{The Shi arrangements of types $A_2$, $B_2$ and
$G_2$.  The positive cone is shaded.}
\label{shi}
\end{figure}

Figure~\ref{shi} shows the Shi arrangements of types $A_2$, $B_2$ and
$G_2$, oriented so as to agree with the root systems as drawn in
Figure~\ref{rank2}. 

\subsection*{$W$-orbits in a discrete torus 
}  

The reflection group $W$ acts on the root lattice~$Q=\ZZ\Phi$, 
hence on the ``discrete torus''
$Q/(h+1)Q$ obtained as a quotient of~$Q$ by its subgroup~$(h+1)Q$.

\begin{theorem}[\cite{haiman}]
The number of $W$-orbits in $Q/(h+1)Q$ is equal to~$N(\Phi)$. 
\end{theorem}

\begin{figure}[ht!]
\centerline{\scalebox{0.8}{\epsfbox{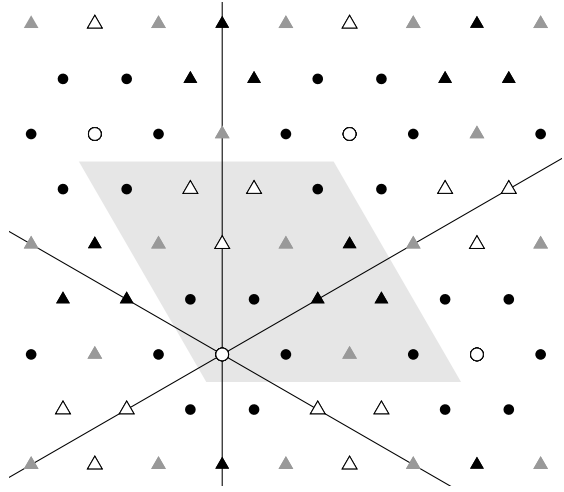}}}
\caption{$A_2$-orbits in $Q/4Q$. 
Each orbit is labeled by a different symbol.
}
\label{A2lattice}
\end{figure}

\begin{figure}[ht!]
\centerline{\scalebox{0.8}{\epsfbox{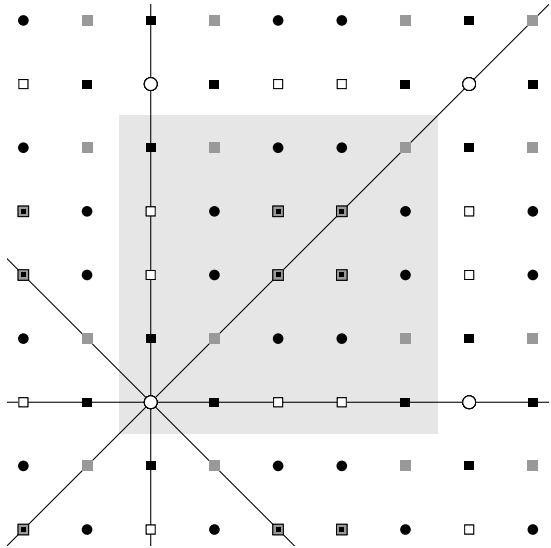}}}
\caption{$B_2$-orbits in $Q/5Q$. 
}
\label{B2lattice}
\end{figure}

Figures~\ref{A2lattice} and~\ref{B2lattice} illustrate these orbits in
types~$A_2$ and~$B_2$, where $h=3$ and $h=4$, respectively. 
Each figure shows the reflection lines of the Coxeter arrangement; 
the shaded region is a fundamental domain for 
the \hbox{translations in~$(h+1)Q$.}

\subsection*{Non-crossing partitions} 

The classical \emph{non-crossing partitions} introduced by Kreweras are
(unordered) partitions of the set $[n+1]=\{1,\dots,n+1\}$ 
into non-empty subsets called \emph{blocks} 
which satisfy the following ``non-crossing'' condition: 
\begin{itemize}
\item
there does not exist an ordered quadruple $(a<b<c<d)$ 
such that the two-element sets $\set{a,c}$ and $\set{b,d}$
are contained in different blocks.
\end{itemize}

\begin{figure}[ht!]
\centerline{
{\epsfbox{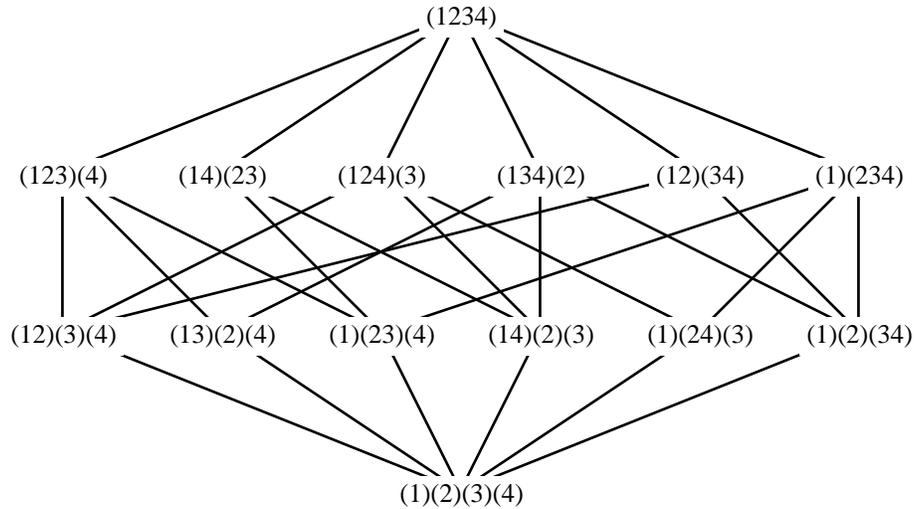}}}
\caption{The non-crossing partition lattice of type $A_3$}
\label{A3nc}
\end{figure}

\begin{figure}[ht!]
\centerline{\scalebox{0.8}
{\epsfbox{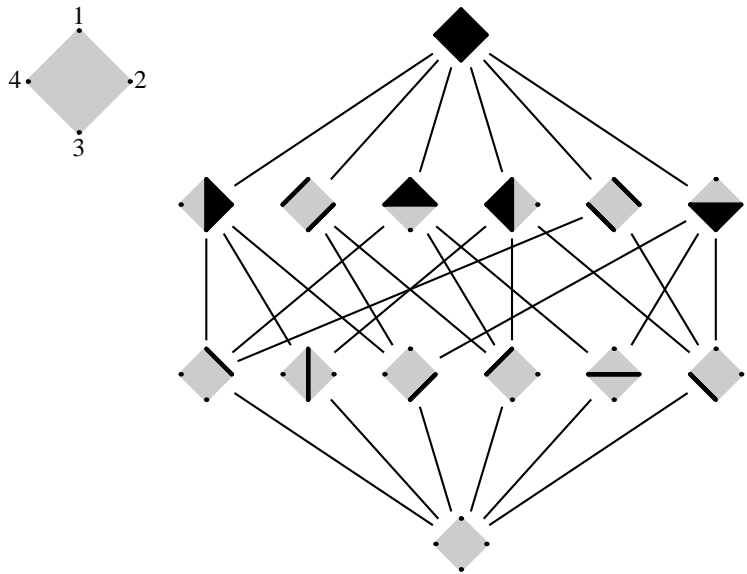}}}
\caption{Planar representation of non-crossing partitions
}
\label{A3nc_pics}
\end{figure}

Figure~\ref{A3nc} shows the $14$ non-crossing partitions for~$n=3$,
partially ordered by refinement.
Such partial order is in fact a lattice for any~$n$;
the number of non-crossing partitions is a Catalan number. 

An alternative way of representing non-crossing partitions 
is shown in Figure~\ref{A3nc_pics}.
Place the elements of $[n+1]$ around a circle.
Then the non-crossing partitions are those set partitions in which the
convex hulls of blocks do not intersect.

We will now explain how this construction arises as a type-$A$
special case of a general construction 
valid for \emph{any} (possibly infinite) Coxeter system~$(W,S)$.

A \emph{reflection} in a Coxeter group $W$ is an element
conjugate to a generator~$s\in S$. 
Any element $w\in W$ can be written as a product of reflections.  
Let $L(w)$ denote the length (i.e., number of factors) of a shortest
such factorization. 
We then partially order $W$ by setting $u\preceq uv$ whenever 
$L(uv)=L(u)+L(v)$, i.e., whenever concatenating shortest
factorizations for $u$ and $v$ gives a shortest factorization
for~$uv$. 
Equivalently, $w$ covers $u$ in this partial order if and only if
$L(w)=L(u)+1$ and there is a reflection $t$ such that $w=ut$. 


Let $c$ be a product (in an arbitrary order) of the generators
in~$S$. 
Thus, $c$ is a \emph{Coxeter element} in~$W$, in the broader sense of
the notion 
alluded to in a footnote 
in Section~\ref{sec:Coxeter element}. 
The \emph{non-crossing partition lattice} for~$W$
(see~\cite{bessis, brady-watt}) 
is the interval $[1,c]$ in the partial order $(W,\preceq)$ defined
above. 
It is a classical result that all Coxeter elements are conjugate to
each other.  
Since the set of all reflections is fixed under
conjugation, it follows that different choices of $c$ yield isomorphic posets. 
(These posets are lattices, which is a non-trivial theorem.) 

The following theorem was obtained in~\cite{bessis, picantin}. 
A version for the classical types $ABCD$ appeared earlier
in~\cite{Reiner-noncrossing}. 

\begin{theorem}
\label{th:non-crossing-W}
Let $\,W$ be the reflection group associated with a finite root
system~$\Phi$. 
Then the non-crossing partition lattice
for~$W$ has $N(\Phi)$ elements. 
\end{theorem}

In type~$A_n$, the general construction presented above
recovers the ordinary non-crossing
partition lattice. To realize why, look again at Figure~\ref{A3nc},
and interpret each element of the poset as a permutation in
$\mathcal{S}_4$ written in cycle notation.

The non-crossing partition lattice of type~$B_n$ can also be given a
direct combinatorial description.
Let us take the ordinary lattice of non-crossing partitions of a $2n$-element
set in its representation illustrated in Figure~\ref{A3nc_pics}.  
Then consider the sublattice consisting of those partitions whose planar 
representations are centrally symmetric.
The result (for $n=3$) is shown in Figure~\ref{B3nc_pics}.

\begin{figure}[ht!]
\centerline{{\epsfbox{B3nc_pics.ps}}}
\caption{The non-crossing partitions of type $B_3$.
}
\label{B3nc_pics}
\end{figure}

\section{Generalized Narayana Numbers}
\label{sec:narayana}

For any enumerative problem whose answer is a Catalan number,
replacing a simple count by a generating function with respect to some
combinatorial statistic results in a $q$-analogue of a Catalan
number. 
There are at least three such $q$-analogues that routinely pop up in
various contexts. 
One is obtained from the usual formula $\frac{1}{n+2}\binom{2n+2}{n+1}$
by replacing $n+2$ and $\binom{2n+2}{n+1}$ with their standard 
             $q$-analogues.
A different answer is obtained while counting order ideals in the root
poset of type~$A_n$ 
by the cardinality of an ideal.
For more on these $q$-analogues, see~\cite{Fur-Hof,Gar-Hai,ec2-q-catalan}.

We will focus on a third $q$-analogue that is related to the Narayana
numbers, defined by the formula $\frac{1}{n+1}\binom{n+1}{k}\binom{n+1}{k+1}$. 
The Narayana numbers form a triangle shown on the right in
Figure~\ref{fig:pascal-narayana}. 
Thus, the numbers in each row of this triangle are obtained by looking
at the corresponding row of Pascal's triangle on the left, 
computing products of consecutive pairs of entries, and dividing them 
by~$n+1$.

\begin{figure}[ht]
\[\begin{array}{ccccccccccc}
&&&&&1\\
&&&&1&&1\\
&&&1&&2&&1\\
&&1&&3&&3&&1\\
&1&&4&&6&&4&&1\\
1&&5&&10&&10&&5&&1
\end{array}
\hspace{34pt}
\begin{array}{ccccccccc}
\\
&&&&1\\
&&&1&&1\\
&&1&&3&&1\\
&1&&6&&6&&1\\
1&&10&&20&&10&&1
\end{array}\]
\caption{The Pascal triangle and the Narayana numbers} 
\label{fig:pascal-narayana}
\end{figure}

Remarkably, the row sums in the triangle of Narayana numbers are the Catalan 
numbers:
\[
\sum_{k=0}^n\frac{1}{n+1}\binom{n+1}{k}\binom{n+1}{k+1} = 
\frac{1}{n+2}\binom{2n+2}{n+1}\,. 
\]
This suggests introducing a $q$-analogue of the Catalan numbers given by 
\begin{equation}
\label{eq:narayana-gf}
\sum_{k=0}^n\frac{1}{n+1}\binom{n+1}{k}\binom{n+1}{k+1}q^k\,.
\end{equation}
We will now explain the connection between this $q$-analogue and the
classical (type~$A$) associahedron.
This connection will lead us to an extension of the definition to other root
systems. 

We will need the notions of the $f$-vector and $h$-vector 
of an $(n\!-\!1)$-dimensional simplicial complex.
The \emph{$f$-vector} is $(f_{-1},f_0,\ldots,f_{n-1})$ where $f_i$ denotes 
the number of $i$-dimensional faces.
The unique ``$(-1)$-dimensional'' face is the empty face.
The \emph{$h$-vector} $(h_0,h_1,\ldots,h_n)$ is determined from
the $f$-vector by the ``reverse Pascal's triangle'' recursion which we illustrate
by an example.

\begin{figure}[ht]
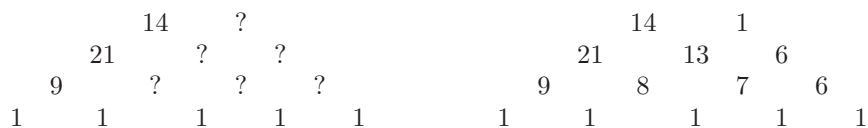

\[\begin{array}{ccccccccc}
&&&14&&?\\
&&21&&?&&?\\
&9&&?&&?&&?\\
1&&1&&1&&1&&1
\end{array}
\qquad\qquad
\begin{array}{ccccccccc}
&&&14&&1\\
&&21&&13&&6\\
&9&&8&&7&&6\\
1&&1&&1&&1&&1\\
\end{array}
\]
\caption{Computing the $h$-vector} 
\label{fig:h-from-f}
\end{figure}

\begin{example}\rm
\label{A3h}
The $f$-vector of the simplicial complex dual to the
associahedron of type~$A_3$ is
$(1,9,21,14)$. (See Figure~\ref{A3assoc_dual}.) 
To calculate the $h$-vector, we place the $f$-vector and 
a row of $1$'s in a triangular array
as shown in Figure~\ref{fig:h-from-f} on the left, 
with most of the entries as yet undetermined.
The remaining entries are then filled in by applying the following
rule: 
each entry is the difference between the entry preceding it in its row
and the entry directly southwest of it. 
Thus, we get $9-1=8$, $21-8=13$, \emph{etc.}
Finally, we obtain the $h$-vector $(1,6,6,1)$ by reading the rightmost
entries in every row. 
Notice that these are exactly the Narayana numbers appearing in the
third row in Figure~\ref{fig:pascal-narayana}. 
\end{example}

\begin{lemma}
\label{lem:asso-narayana}
The components of the $h$-vector of the simplicial complex 
dual to an $n$-dimensional associahedron 
are the Narayana numbers 
$\frac{1}{n+1}\binom{n+1}{k}\binom{n+1}{k+1}$. 
\end{lemma}

Motivated by Lemma~\ref{lem:asso-narayana}, 
we define the (generalized) Narayana numbers $N_k(\Phi)$ ($k=0,\dots,n$) 
for an arbitrary root system $\Phi$ as the
entries of the $h$-vector of the simplicial complex 
dual to the corresponding generalized
associahedron.

\begin{example}\rm
\label{B3h}
The $f$-vector of the simplicial complex 
dual to the $3$-dimensional cyclohedron
(the associahedron of type~$B_3$) is $(1,12,30,20)$. 
The corresponding $h$-vector is $(1,9,9,1)$.
In general, the Narayana numbers of type $B_n$ are the squares of
entries of Pascal's triangle: $N_k(B_n)=\binom{n}{k}^2$.
\end{example}

It is easy to see that 
the entries of an $h$-vector always add up to~$f_{n-1}$, 
the number of top-dimensional faces in the simplicial complex. 
Thus, \hbox{$\sum_k N_k(\Phi)=N(\Phi)$.} 
Consequently, the generating function for the Narayana numbers of
type~$\Phi$ 
\[
N(\Phi,q)=\textstyle\sum_{k=0}^n N_k(\Phi) q^k
\]
provides a $q$-analogue of~$N(\Phi)$ 
which generalizes~\eqref{eq:narayana-gf}. 
These generating functions for the finite crystallographic root systems 
are tabulated in Figure~\ref{fig:narayana}.

\begin{figure}[ht]
\begin{eqnarray*}
N(A_n,q)&=&\sum_{k=0}^n\frac{1}{n+1}\binom{n+1}{k}\binom{n+1}{k+1}q^k\\
N(B_n,q)&=&\sum_{k=0}^n{\binom{n}{k}\!}^2q^k\\
N(D_n,q)&=&1+q^n+
        \sum_{k=1}^{n-1}
                \left[{\binom{n}{k}\!}^2-
                        \frac{n}{n-1}\binom{n-1}{k-1}\binom{n-1}{k}\right]q^k\\
N(E_6,q)&=&1+36q+204q^2+351q^3+204q^4+36q^5+q^6\\
N(E_7,q)&=&1+63q+546q^2+1470q^3+1470q^4+546q^5+63q^6+q^7\\
N(E_8,q)&=&1+120q+1540q^2+6120q^3+9518q^4
        \\&&\qquad\qquad\qquad\qquad\qquad+6120q^5+1540q^6+120q^7+q^8\\
N(F_4,q)&=&1+24q+55q^2+24q^3+q^4\\
N(G_2,q)&=&1+6q+q^2
\end{eqnarray*}
\caption{Generating functions for generalized Narayana numbers}
\label{fig:narayana}
\end{figure}

The Narayana numbers provide refined counts for the various interpretations of 
$N(\Phi)$ given in Section~\ref{sec:numerology}. 
These enumerative results are listed in
Theorem~\ref{th:phi-narayana-combinat} below;
we elaborate on the items in the theorem in subsequent comments. 

Theorem~\ref{th:phi-narayana-combinat} is a combination of results in
\cite{Ath-TAMS,ga,panyushev,reiner-welker,sommers};
see \cite{Ath-TAMS} for a historical overview, and for further
generalizations. 

\pagebreak

\begin{theorem}
\label{th:phi-narayana-combinat}
The following numbers are equal to each other, and to~$N_k(\Phi)$: 
\begin{itemize}
\item[(i)]
the $k$th component of the $h$-vector for the dual complex of a
generalized associahedron of type~$\Phi$; 
\item[(ii)]
the number of elements of rank~$k$ in the
non-crossing partition lattice for~$W$; 
\item[(iii)]
the number of antichains of size $k$ in the root 
poset for~$\Phi$; 
\item[(iv)]
the number of $W$-orbits in $Q/(h+1)Q$ consisting of elements 
whose stabilizer has rank~$k$; 
\item[(v)]
the components of the $h$-vector for the dual cell complex 
of the positive part of the Shi arrangement. 
\end{itemize}
\end{theorem}

\begin{remark}[Comments on Theorem~\ref{th:phi-narayana-combinat}]
{\rm{\ }

(i) This was our definition of $N_k(\Phi)$. 
\smallskip

(ii) 
The lattice of non-crossing partitions of type~$\Phi$ is graded, 
and $N_k(\Phi)$ is the number
of elements of rank $k$.
\smallskip

(iii) 
The $h$-vector of any simplicial polytope satisfies the 
\emph{Dehn-Sommer\-ville equations} $h_i=h_{d-i}$. 
Thus interpretation (i) implies that
$N_k(\Phi)=N_{n-k}(\Phi)$. 
This symmetry of the 
Narayana numbers is also apparent in the interpretation~(ii) 
because the non-crossing partition lattices are self-dual. 
However, this symmetry is not at all obvious in the interpretations
(iii)--(v). 
In particular, no direct combinatorial explanation is known for why 
the number of antichains of size $k$ in the root poset is the same as
the number of antichains of size~$n-k$.  
\smallskip

(iv)
The stabilizer of an element in $Q/(h+1)Q$  is a reflection
subgroup of~$W$.
The stabilizers of elements in the same $W$-orbit are conjugate,
and therefore have the same rank. 
$N_k(\Phi)$ is the number of orbits in which the stabilizers have rank~$k$.
For example, in type $A_2$ there is $1$ orbit whose stabilizer has
rank~$2$ (the unfilled circle in
Figure~\ref{A2lattice}), $3$ orbits whose stabilizers have rank~$1$ 
(each symbolized by a triangle) and
$1$ orbit whose stabilizers have rank~$0$ (the filled circles),
in agreement with $N(A_2,q)=1+3q+q^2$. 
\smallskip

(v)
The positive regions of the Shi arrangement can be used to define a ``dual''
cell complex.
The vertices of this complex correspond to the positive regions of~the
Shi arrangement. 
The faces of the complex correspond to those faces of the closures of
these regions that are not contained in the boundary of the positive
cone. 
Accordingly, the maximal faces correspond to the vertices of the
arrangement which lie  in the interior of the positive cone.
See Figure~\ref{A2shidual}. 
Amazingly, this cell complex has the same $f$-vector (hence the same 
$h$-vector) as the corresponding associahedron.
In the example of Figure~\ref{A2shidual}, we get $5$~vertices,
$5$~faces, and $1$~two-dimensional face, matching the numbers for the
pentagon (the type~$A_2$ associahedron). 

\begin{figure}[ht]
\centerline{\scalebox{0.8}{\epsfbox{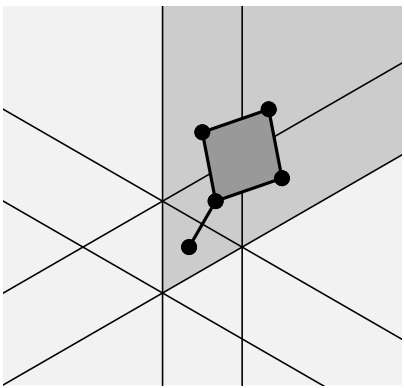}}}
\caption{The dual complex for the positive part of the Shi arrangement
  of type~$A_2$.}
\label{A2shidual}
\end{figure}
}

\end{remark}

\section{Non-crystallographic types}
\label{sec:non-crys assoc}

\begin{figure}[ht]
\centerline{
\scalebox{1.2}{\epsfbox{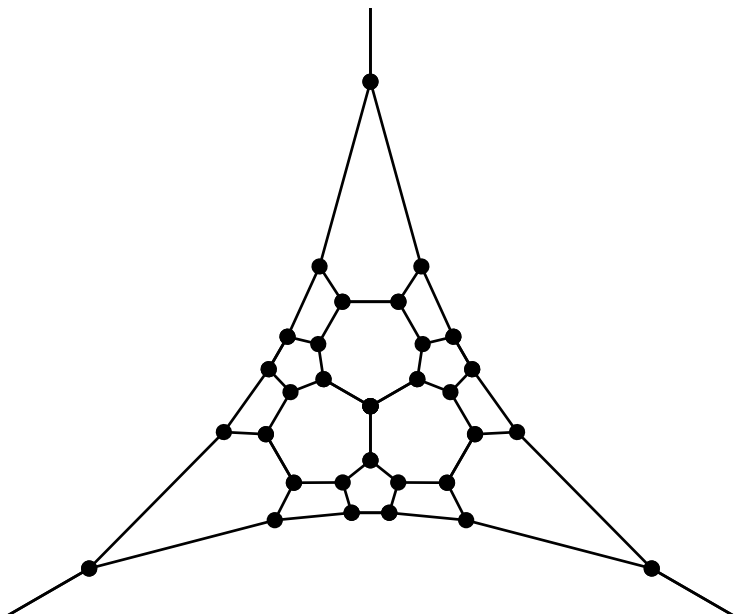}}}   
\caption{The 
associahedron of type~$H_3$}
\label{H3assoc}
\end{figure}

The construction of generalized associahedra via Definition~\ref{def:tau} and 
Theorems~\ref{th:dihedral} and~\ref{th:compat} can be carried out verbatim
for the non-crystallographic root systems $I_2(m)$, $H_3$ and $H_4$.
(However, the last sentence of Theorem~\ref{th:compat} must be ignored, 
since no ``cluster complex'' exists for non-crystallographic root systems.)
The associahedron of type $I_2(m)$ is an $(m+2)$-gon. 
The $1$-skeleton of the associahedron for $H_3$ is shown in Figure~\ref{H3assoc}.
(The vertex at infinity completes the three unbounded regions to heptagons.)

The analogue of Theorem~\ref{th:N-thru-exponents} holds true in types 
$H_3$, $H_4$, and~$I_2(m)$: 
the number of vertices of a generalized associahedron  
is equal to~$N(\Phi)$. 
The latter number is still given
by~\eqref{eq:N-thru-exponents}, with the exponents
taken from Figure~\ref{data}. 
Figure~\ref{W-Cat-noncrystal} shows these values of $N(\Phi)$ explicitly. 

\begin{figure}[ht]
\begin{center}
\begin{tabular}{|c|c|c|}
\hline
$H_3$ & $H_4$ & $I_2(m)$ \\
\hline
&&\\[-.1in]
$\textstyle 32$
& $\textstyle 280$
& $\textstyle m+2$ \\[.05in]
\hline
\end{tabular}
\end{center}
\caption{The numbers $N(\Phi)$ in non-crystallographic cases} 
\label{W-Cat-noncrystal}
\end{figure}

The corresponding $h$-vectors (``Narayana numbers'') are given by
\begin{eqnarray*}
N(I_2(m),q)&=&1+mq+q^2,\\
N(H_3,q)&=&1+15q+15q^2+q^3,\\
N(H_4,q)&=&1+60q+158q^2+60q^3+q^4.
\end{eqnarray*}

The construction of the non-crossing partition lattice does not require a
crystallographic Coxeter group. Theorem~\ref{th:non-crossing-W} and
Theorem~\ref{th:phi-narayana-combinat}(ii) remain valid for the finite
non-crystallographic root systems.
At present, the other manifestations of $N(\Phi)$ and $N_k(\Phi)$ 
presented in Sections~\ref{sec:numerology} and~\ref{sec:narayana} 
(including Parts (iii)--(v) of
Theorem~\ref{th:phi-narayana-combinat}) 
do not appear to extend to the non-crystallographic cases.

\section{Lattice congruences and the weak order}
\label{sec:lattice cong}

This section is based on~\cite{cambrian}. 
Its main goal is to establish a relationship between two fans
associated with a root system~$\Phi$ and the corresponding reflection
group~$W$: 
\begin{itemize}
\item
the \emph{Coxeter fan} created by (the regions of) the Coxeter
arrangement, and 
\item
the \emph{cluster fan} described in
Theorem~\ref{th:cluster-complex-fintype}. 
\end{itemize}
These fans are the normal fans of a permutahedron and an
associahedron of the corresponding type, respectively. 

Let $\omega_i$ denote the \emph{fundamental weight}~\cite{Bourbaki} 
corresponding to~$\alpha_i$. 
For $i\!\in \!I$, we~set 
\[
\ep(i)=\begin{cases}
+1 & \text{if $i\in I_+\,$,} \\
-1 & \text{if $i\in I_-\,$.} 
\end{cases}
\]

\begin{theorem}
\label{cluster refine}
The linear automorphism $Q_\RR\to Q_\RR$ defined by 
$\alpha_i\mapsto\ep(i)\omega_i$ 
moves the cluster fan to a fan refined by the Coxeter fan. 
\end{theorem}

The gluing of maximal cones of the Coxeter fan corresponds to contraction of
edges in the $1$-skeleton of a permutahedron. 
By Theorem~\ref{cluster refine}, this can be done in such a way that 
the result of the contraction is the $1$-skeleton of a generalized
associahedron. 
We have thus extended Theorem~\ref{ABskeleta} to all types.

\medskip

The statement of Theorem~\ref{cluster refine} does not specify which
regions of the Coxeter arrangement should be combined together to produce the
maximal cones of the transformed cluster fan. 
We next present a lattice-theoretic construction that, conjecturally,
answers this question. 


The \emph{weak order} on $W$ is the partial order in which $u\le v$ if and 
only if some reduced word for $u$ occurs as an initial segment of a reduced
word for~$v$.
In particular, $v$ covers $u$ in the weak order if and only if 
$u^{-1}v$ is a simple reflection, and the
length of $v$ is greater than the length of~$u$ (necessarily by~$1$). 
Lemma~\ref{lem:adj-regions} (see also the paragraph that follows it)
implies that the Hasse diagram of the weak order can be identified
with the $1$-skeleton of a $W$-permutohedron. 

\begin{theorem}[\cite{bjorner-orderings}]
The weak order on a finite Coxeter group is a lattice.
\end{theorem}

\begin{figure}[ht!]
\centerline{\scalebox{.8}{\epsfbox{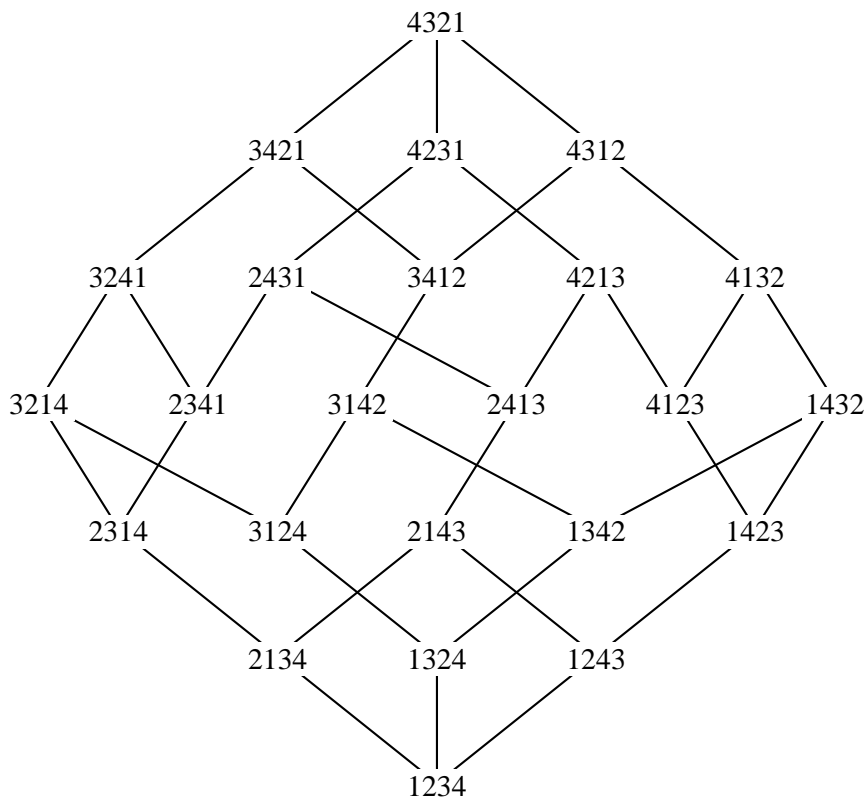}}}
\caption{The weak order of type $A_3$}
\label{weakA3}
\end{figure}

\begin{example}{\rm 
The weak order of type~$A_n$ can be described in the language of permutations of
$[n\!+\!1]$, written in one-line notation.
Permutation $v=(v_1,\dots,v_{n+1})$ covers 
$u\!=\!(u_1,\dots,u_{n+1})$ if $v$ is obtained from~$u$ 
by exchanging two 
entries $u_i$ and $u_{i+1}$ with 
$u_i\!<\!u_{i+1}\,$. 
Figure~\ref{weakA3} shows the weak order on~$A_3$. (Cf.\
Figure~\ref{A3perm}.) 
}
\end{example}

A \emph{congruence} on a lattice is an equivalence relation which respects the
meet and join operations. 
A (bipartite) \emph{Cambrian congruence} on the weak order of~$W$ 
is defined as the (unique) coarsest congruence~``$\equiv$'' such that,
for each edge $(s,t)$ in the Coxeter diagram, with $t\in I_-$, we have 
\[
t\equiv tsts\cdots \quad \text{($m_{st}-1$ factors).}
\]



\begin{figure}[ht]
\centerline{\scalebox{.8}{\epsfbox{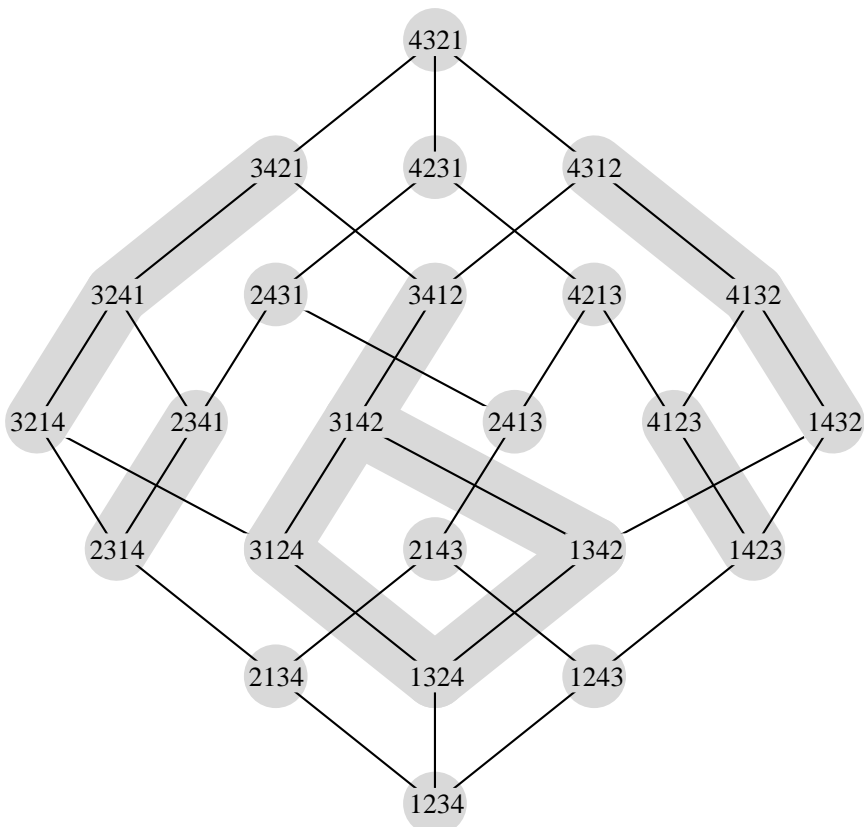}}}
\caption{A bipartite Cambrian congruence of type~$A_3$}
\label{weakS4cong}
\end{figure}

\begin{example}\rm
Figure~\ref{weakS4cong} shows the bipartite Cambrian congruence for 
$W$ of type~$A_3$, i.e., the coarsest congruence on the weak order of the
symmetric group~$\mathcal{S}_4$
such that $1324\equiv 3124$ and $1324\equiv 1342$.
The congruence classes are shaded.
\end{example}

\begin{conjecture}
\label{conj:cluster-cambrian}
Two regions $R_u$ and $R_v$ of the Coxeter arrangement are contained in the same 
maximal cone of the transformed cluster fan (see Theorem~\ref{cluster
  refine})
if and only if $u\equiv v$ under the bipartite Cambrian congruence. 
\end{conjecture}

Conjecture~\ref{conj:cluster-cambrian} has been proved in types~$A_n$
and~$B_n$.  The proof makes explicit the combinatorics of the Cambrian
congruence and connects it to constructions given by
Billera and Sturmfels~\cite{Iterated} (type~$A$) and
Reiner~\cite{Equivariant} (type~$B$). The conjecture implies in particular
that the Hasse diagram of the quotient of the weak order by the Cambrian
congruence (called the \emph{Cambrian lattice}) is isomorphic to the
$1$-skeleton of the generalized associahedron. 

\begin{figure}[ht]
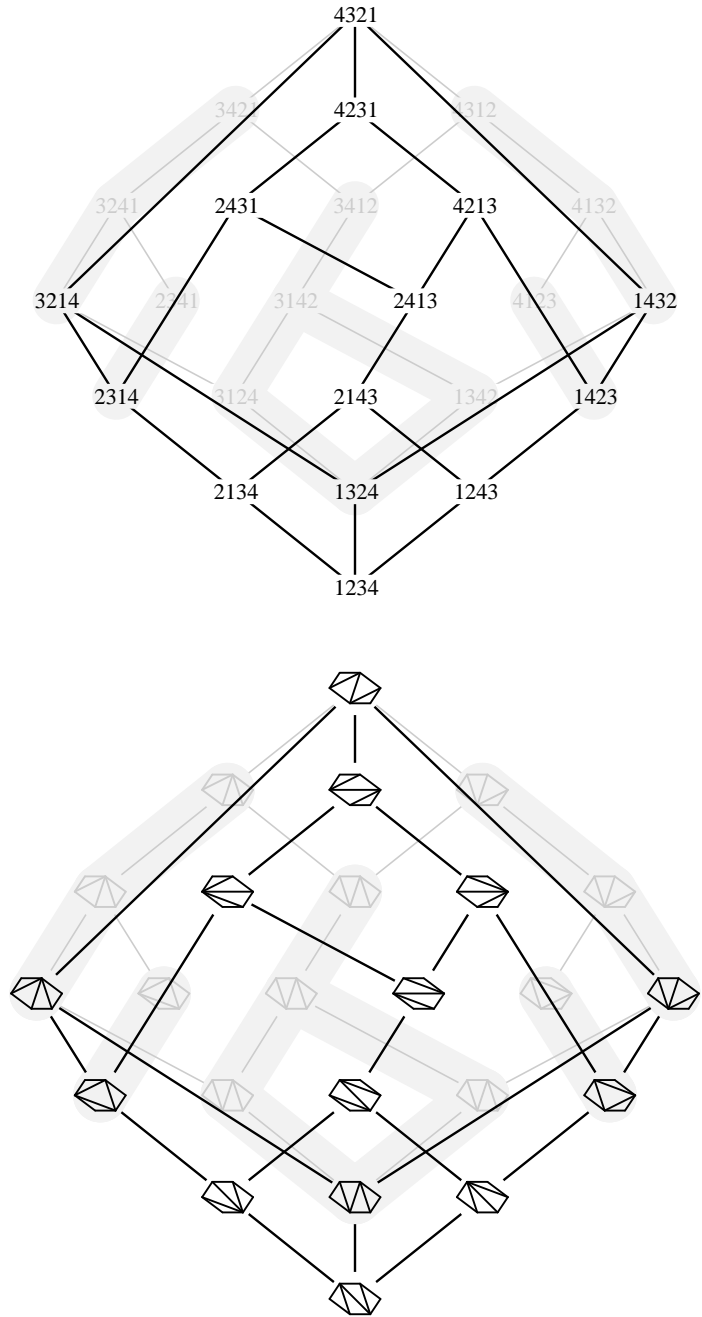

\vspace{.1in}
\centerline{\scalebox{.6}{\epsfbox{weakS4bottoms.ps}}}
%
\vspace{.3in}
\centerline{\scalebox{.6}{\epsfbox{permtri_bottoms.ps}}}
\caption{A bipartite Cambrian lattice of type $A_3$}
\label{permtri_bottoms}
\end{figure}

More concretely, the
Cambrian lattice is obtained as the induced subposet of the weak order
formed by taking the (unique) smallest element in each (Cambrian)
congruence class; see Figure~\ref{permtri_bottoms}. 
We omit the description of the bijection used to translate the top
picture in Figure~\ref{permtri_bottoms} (the Cambrian lattice labeled  by
permutations) into the bottom one (the associahedron labeled by
triangulations).


\newcommand{\journalname}[1]{\textrm{#1}}
\newcommand{\booktitle}[1]{\textrm{#1}}

\end{document}

%% file: pcmslmod-modified.tex
\makeatletter

\newif\iffirstlecture\firstlecturefalse

\newcommand{\lectureseries}{\firstlecturetrue
              \secdef\@lectureseries\@slectureseries} 

\newcommand{\@lectureseries}[2][default]{\chapter*{#2}%
              \gdef\thelectureseries{#1}} 

\newcommand{\@slectureseries}[1]{\chapter*{#1}}

\renewcommand{\auth}{\secdef\@auth\@sauth}

\newcommand{\@auth}[2][default]{\vspace{-1pc}{\raggedleft
        \Large\bf\noindent
        #2\endgraf
        \vspace*{2pc}
        }
        \def\@author{#1}%
}

\newcommand{\@sauth}[1]{\vspace{-1pc}{\raggedleft
        \Large\bf\noindent
        #1\endgraf
        \vspace*{2pc}
        }
        \def\@author{#1}%
}

\def\lecture#1{\global\Lecturetrue\global\Monographfalse
\iffirstlecture\else\chapter*{}\fi\firstlecturefalse
  \global\let\sectionmark\@gobble 
  \addtocounter{lecture}1\relax
  \refstepcounter{chapter}%
\gdef\thelecturename{#1\unskip}
  {\Large\bfseries
   \raggedleft
   \@xp\uppercase\@xp{\thelecturelabel} {\LARGE\thelecturenum}\\
   \vspace*{3pt}%
   \thelecturename
   \endgraf}%
  \let\@secnumber=\thelecturenum
  \@xp\lecturemark\@xp{\thelecturename}%
  \addcontentsline{toc}{chapter}{%
    \thelecturelabel\ \thelecturenum.\ \thelecturename}%
  \vspace{34\p@}\noindent}
  
\def\lecture{\global\Lecturetrue\global\Monographfalse
\iffirstlecture\else\chapter*{}\fi%
  \global\let\sectionmark\@gobble 
\secdef\@lecture\@slecture}

\def\@lecture[#1]#2{%
  \addtocounter{lecture}1\relax
  \refstepcounter{chapter}%
\gdef\thelecturename{#1\unskip}\firstlecturefalse
  {\Large\bfseries
   \raggedleft
   \@xp\uppercase\@xp{\thelecturelabel} {\LARGE\thelecturenum}\\
   \vspace*{3pt}%
    #2\unskip
   \endgraf}%
  \let\@secnumber=\thelecturenum
  \@xp\lecturemark\@xp{\thelecturename}%
  \addcontentsline{toc}{chapter}{%
    \thelecturelabel\ \thelecturenum.\ #2}%
  \vspace{34\p@}\noindent}
  
\def\slecturerunhead#1#2#3{%
    \let\@tempa\chaptername
    \uppercasenonmath{\@tempa}%
    \def\@tempb{#3\unskip}%
    \uppercasenonmath{\@tempb}%
    {\normalfont\@tempb}
    }
\def\slecturemark{
    \@secmark\markright\slecturerunhead\chaptername}%

\def\@slecture#1{%
\iffirstlecture
\gdef\thelecturename{#1\unskip}\firstlecturefalse
  {\Large\bfseries
\noindent\thelecturename
   \endgraf}%
  \let\@secnumber=\thelecturenum
  \@xp\slecturemark\@xp{\thelecturename}%
  \addcontentsline{toc}{chapter}{%
    \thelecturename}%
 \vspace{-6\p@}\noindent
\else
\gdef\thelecturename{#1\unskip}\firstlecturefalse
  {\Large\bfseries
   \raggedleft
   \@xp\uppercase\@xp{\thelecturename}
   \endgraf}%
  \let\@secnumber=\thelecturenum
  \@xp\slecturemark\@xp{\thelecturename}%
  \addcontentsline{toc}{chapter}{%
    \thelecturename}%
  \vspace{34\p@}\noindent
\fi}

\ifLecture
  \def\chapterrunhead#1#2#3{%
    \let\@tempa\@author
    \uppercasenonmath{\@tempa}%
    \uppercasenonmath{\thelectureseries}%
    \textmd{\@tempa, \thelectureseries}
    }
  \def\lecturerunhead#1#2#3{%
    \let\@tempa\chaptername
    \uppercasenonmath{\@tempa}%
    \def\@tempb{#3\unskip}%
    \uppercasenonmath{\@tempb}%
    \textmd{\@tempa\ #2. \@tempb}
    }
\else
  \let\chapterrunhead\partrunhead
\fi

\newif\ifBibliographyIsASection\BibliographyIsASectionfalse

  \def\bibliomark{
    \@secmark\markright\bibliorunhead\chaptername}%

  \def\bibliorunhead#1#2#3{%
    \let\@tempa\chaptername
    \uppercasenonmath{\@tempa}%
    \def\@tempb{#3\unskip}%
    \uppercasenonmath{\@tempb}%
    \textmd{\@tempb}
    }

\def\thebibliography#1{%
  \ifBibliographyIsASection
    \section*\refname
    \if@backmatter
      \markboth{\refname}{\refname}%
    \fi
  \else
\chapter*{}
  {\Large\bfseries
   \raggedleft
   \@xp\uppercase\@xp{\bibname} \\
   \endgraf}%
  \let\@secnumber=\thelecturenum
  \@xp\bibliomark\@xp{\bibname}%
  \addcontentsline{toc}{chapter}{%
    \bibname}%
  \vspace{34\p@}\noindent
  \fi
  \normalsize\labelsep .5em\relax
  \list{\@arabic\c@enumi.}{\settowidth\labelwidth{\@biblabel{#1}}%
  \leftmargin\labelwidth
  \advance\leftmargin\labelsep
	\usecounter{enumi}}\sloppy
  \clubpenalty9999 \widowpenalty\clubpenalty  \sfcode`\.\@m}

  \def\indexmark{
    \@secmark\markright\indexrunhead\chaptername}%

  \def\indexrunhead#1#2#3{%
    \let\@tempa\chaptername
    \uppercasenonmath{\@tempa}%
    \def\@tempb{#3\unskip}%
    \uppercasenonmath{\@tempb}%
    \textmd{\@tempb}
    }

\ifLecture
\def\theindex{\cleardoublepage
\@restonecoltrue\if@twocolumn\@restonecolfalse\fi
\columnseprule \z@ \columnsep 35pt
\def\indexchap{\@startsection
		{chapter}{1}{\z@}{8pc}{34pt}%
		{\raggedleft
		\Large\bfseries}}%
 \twocolumn[\indexchap[{\indexname}]{\@xp\uppercase\@xp{\indexname}}]
  \@xp\indexmark\@xp{\indexname}%
	\thispagestyle{plain}\let\item\@idxitem\parindent\z@
	 \footnotesize\parskip\z@ plus .3pt\relax\let\item\@idxitem}
\fi

\def\@makefntext{\noindent\@makefnmark}

\def\setaddress{%
  {\let\@makefnmark\relax  \let\@thefnmark\relax
        \nobreak
        \addressnum@=\z@
        \loop\ifnum\addressnum@<\addresscount@\advance\addressnum@\@ne
           \footnote{
           \csname @address\number\addressnum@\endcsname
           \csname @curraddr\number\addressnum@\endcsname
           \csname @email\number\addressnum@\endcsname}\repeat
  \ifx\@empty\@date\else \@footnotetext{\@setdate}\fi
  \ifx\@empty\@subjclass\else \@footnotetext{\@setsubjclass}\fi
  \ifx\@empty\@keywords\else \@footnotetext{\@setkeywords}\fi
  \ifx\@empty\thankses\else \@footnotetext{%
    \def\par{\let\par\@par}\@setthanks}\fi
    }%
  \@setcopyright
}

\def\@tmpevenhead{\relax}

\def\cleardoublepage{\clearpage\if@twoside \ifodd\c@page\else
    \let\@tmpevenhead\@evenhead \let\@evenhead\relax\hbox{}\eject 
    \let\@evenhead\@tmpevenhead\if@twocolumn\hbox{}\newpage\fi\fi\fi}

\def\@setcopyright{%
  \let\copyrightyear\currentyear             
  \insert\copyins{\hsize\textwidth
    \parfillskip\z@ \leftskip\z@\@plus.9\textwidth
    \fontsize{6}{7\p@}\normalfont\upshape
    \everypar{}%
    \vskip-\skip\copyins \nointerlineskip
    \noindent\vrule\@width\z@\@height\skip\copyins
    \copyright\copyrightyear\ 
Sergey Fomin and Nathan Reading\par
    \kern\z@}%
}

\renewcommand{\@auth}[2][default]{{\raggedleft
        \begingroup
  \fontsize{\@xivpt}{18}\bfseries
  #2\par \endgroup
        \vspace*{2pc}
        }
        \def\@author{#1}%
}

\renewcommand{\@sauth}[1]{{\raggedleft
        \begingroup
  \fontsize{\@xivpt}{18}\bfseries
  #1\par \endgroup
        \vspace*{2pc}
        }
        \def\@author{#1}%
}

\def\@lecture[#1]#2{%
  \addtocounter{lecture}1\relax
  \refstepcounter{chapter}%
\gdef\thelecturename{#1\unskip}\firstlecturefalse
  {\Large\bfseries
   \raggedleft
   \@xp\uppercase\@xp{\thelecturelabel} {\LARGE\thelecturenum}\\
   \vspace*{3pt}%
    #2\unskip
   \endgraf}%
  \let\@secnumber=\thelecturenum
  \@xp\lecturemark\@xp{\thelecturename}%
  \addcontentsline{toc}{chapter}{%
    \thelecturelabel\ \thelecturenum.\ #2}%
  \vspace{10\p@}\noindent}
  
\def\@slecture#1{%
\iffirstlecture
\gdef\thelecturename{#1\unskip}\firstlecturefalse
  {\Large\bfseries
\noindent\thelecturename
   \endgraf}%
  \let\@secnumber=\thelecturenum
  \@xp\slecturemark\@xp{\thelecturename}%
  \addcontentsline{toc}{chapter}{%
    \thelecturename}%
 \vspace{-6\p@}\noindent
\else
\gdef\thelecturename{#1\unskip}\firstlecturefalse
  {\Large\bfseries
   \raggedleft
   \@xp\uppercase\@xp{\thelecturename}
   \endgraf}%
  \let\@secnumber=\thelecturenum
  \@xp\slecturemark\@xp{\thelecturename}%
  \addcontentsline{toc}{chapter}{%
    \thelecturename}%
  \vspace{10\p@}\noindent
\fi}

\makeatother